\documentclass[12pt,A4]{amsart}%
\usepackage{amsfonts}
\usepackage{amsmath}
\usepackage{amssymb}
\usepackage{graphicx}
\usepackage{enumerate}
\usepackage{multicol}
\usepackage{color}
\usepackage{mathtools}
\usepackage{datetime2}
\usepackage{pstricks, pst-node, pst-text, pst-3d}
\usepackage[arrow, matrix, curve, xdvi, dvips]{xy}
\usepackage{ulem}

\definecolor{ufogreen}{rgb}{0.24, 0.80, 0.44}
\definecolor{uablue}{rgb}{0.0, 0.2, 0.67}
\definecolor{purpleheart}{rgb}{0.41, 0.21, 0.61}
\definecolor{pinegreen}{rgb}{0.0, 0.47, 0.44}

\definecolor{mgreen}{rgb}{0.1, 0.7, 0.1}

\newtheorem{theorem}{Theorem}
\newtheorem{corollary}{Corollary}
\newtheorem{lemma}{Lemma}
\newtheorem{definition}{Definition}
\newtheorem{proposition}{Proposition}
\newtheorem{remark}{Remark}

\DeclareMathOperator{\pr}{pr}

\DeclareMathOperator{\tr}{trace}
\newcommand{\be}{\begin{equation}}
\newcommand{\ee}{\end{equation}}
\newcommand{\bean}{\begin{equation*}}
\newcommand{\ean}{\end{equation*}}
\newcommand{\ba}{\begin{array}}
\newcommand{\ea}{\end{array}}
\newcommand{\bpm}{\begin{pmatrix}}
\newcommand{\epm}{\end{pmatrix}}
\newcommand{\Ad}{\text{\rm Ad}}

\newcommand{\ad}{\text{\rm ad}}

\newcommand{\End}{\text{\rm End}}

\newcommand{\fg}{\mathfrak {g}}
\newcommand{\fh}{\mathfrak {h}}

\newcommand{\fp}{\mathfrak {p}}

\newcommand{\M}{\widehat{M}}

\newcommand{\field}[1]{\mathbb{#1}}
\newcommand{\liealg}[1]{\mathfrak{#1}}
\newcommand{\Stiefel}[2]{\mathrm{St}_{{#1}, {#2}}}
\newcommand{\at}[2][\big]{#1\vert_{#2}}
\newcommand{\Secinfty}{\Gamma^{\infty}}

\DeclareMathOperator{\D}{d}
\newcommand{\matR}[2]{\field{R}^{{#1} \times {#2}}}
\newcommand{\id}{\mathrm{id}}
\newcommand{\Stab}{\mathrm{Stab}}
\newcommand{\tensor}{\otimes}
\DeclareMathOperator{\vect}{vec}
\DeclareMathOperator{\Span}{span}
\newcommand{\e}{\mathrm{e}}

\newcommand{\exIso}{T}

\usepackage{hyperref}

\begin{document}

\title{Rolling Stiefel manifolds equipped with $\alpha$-metrics}
%    author one information
\author{M. Schlarb}
\address{Institute of Mathematics, Julius-Maximilians-Universit\"at W\"urzburg, 
	\newline \hspace*{0,4 cm}W\"urzburg, Germany}
\email{markus.schlarb@mathematik.uni-wuerzburg.de}
\thanks{.}

\author{K. H\"uper}
\address{Institute of Mathematics, Julius-Maximilians-Universit\"at W\"urzburg, 
 \newline \hspace*{0,4 cm}W\"urzburg, Germany}
\email{hueper@mathematik.uni-wuerzburg.de}
\thanks{.}

%    author two information
\author{I. Markina}
\address{Department of Mathematics, University of Bergen, 
 \newline \hspace*{0,4 cm}P.O.~Box 7803, Bergen N-5020, Norway}
\email{irina.markina@uib.no}
\thanks{.}

%    author three information
\author{F. Silva Leite}
\address{Institute of Systems and Robotics, Uni\-versity of Coimbra,
 \newline \hspace*{0,4 cm}P\'olo II, 3030-290 Coimbra, Portugal \hspace*{0,4 cm} and
 \newline \hspace*{0,4 cm}Department of Mathematics, University of Coimbra,
 \newline \hspace*{0,4 cm}Largo D. Dinis, 3000-143 Coimbra, Portugal}
\email{fleite@mat.uc.pt}

\begin{abstract}
	We discuss the rolling, without slip and without twist, of Stiefel manifolds
	equipped with $\alpha$-metrics, from an intrinsic and an extrinsic point of view.
	We, however,  start with a more general perspective, namely	by investigating intrinsic rolling of normal na\-turally
	reductive homogeneous spaces.
	This gives evidence why a seemingly
	straightforward generalization of intrinsic rolling of
	symmetric spaces to normal naturally reductive
	homogeneous spaces is not possible, in general.
	For a given control curve, we derive a system of explicit time-variant  ODEs
	 whose solutions describe the desired rolling.
	These findings are applied to obtain the
	intrinsic rolling of Stiefel manifolds, which is then extended
	to an extrinsic one.
	Moreover, explicit solutions of the kinematic equations
	are obtained provided that the 
	development curve is the projection of a not necessarily horizontal
	one-parameter subgroup.
	In addition, our results are put into perspective with examples of rolling
	 Stiefel manifolds known from the literature.\\
	 
	 \noindent Keywords: Intrinsic rolling, Extrinsic rolling, Stiefel manifolds, normal naturally
	reductive homogeneous spaces, Covariant derivatives, Parallel vector fields, Kinematic equations.\\
	
	\noindent MSC Classification: 53B21,
	53C30, 53C25, 37J60, 58D19.
	 
\end{abstract}

\date {\DTMnow}
\maketitle
%%%%%%%%%%%%%%%%%%%%%%%%%%%%%%%%%
%\tableofcontents

%%%%%%%%%%%  Section 1  %%%%%%%%%%%%

\section{Introduction}

In recent years there has been increasing interest in so-called
rolling maps of differentiable manifolds. Researchers have taken
different points of view to study the differential geometry behind
these constructions. From our point of view it seems to be natural to
distinguish between two approaches, the intrinsic and the extrinsic
one. The first viewpoint does not require any embedding space to study
rolling maps, whereas the second needs one. At first glance the
intrinsic approach seems to be of more pure mathematical flavour,
simply because intrinsic properties stay in the foreground and any
influence of an embedding space, which might a priori not be known or
even considered to be artificial, will be ignored. In some sense, in
that framework choosing coordinates is a no-go. On the other hand,
however, the extrinsic approach might be considered to be of more
applied character,
mainly because some of the related applications actually stem from
rolling rigid or convex
bodies in the geometric mechanic sense and/or from closely related
questions of geometric
control. 
Although there is an overlap of both approaches, i.e. interpretations of mathematical results of rolling without slip or twist have partially been discussed from both sides, definitions usually differ, including assumptions and consequences. We want to emphasize that by extrinsic we do not mean working with coordinates in the sense of charts. The access to an embedding vector space often opens nevertheless the path to a coordinate free approach, similar to treating the standard sphere $S^n$ embedded into $\mathbb{R}^{n+1}$.

The purpose of this paper is at least threefold. Firstly, we put both approaches, intrinsic and extrinsic, into perspective, clarifying the sometimes subtle differences and discuss their consequences. In particular, we claim that the role of the no-twist conditions become more clarified. Secondly, we study a sufficiently rich class of manifolds, namely the rolling of normal natural reductive homogeneous spaces. An essentially constructive procedure to generalize the rolling of symmetric spaces is presented here for the first time. Thirdly, rolling Stiefel manifolds, serves as our role model, as it is well known that although spheres and orthogonal groups within the set of real Stiefel manifolds are symmetric spaces, all the others are not. We also put all our results into perspective by comparing them to partial results scattered in the literature.

Central to our treatise is the derivation of the so-called kinematic
equations, i.e. a set of ODEs to be considered under certain
nonholonomic constraints. Certainly, the rich theory behind
differential geometric distributions, fiber bundle constructions, and
differential systems can be applied here. For many examples, however,
this theory often does not support explicit solutions for the
nonholonomic problem of rolling with no-slip and no-twist. Here we
present explicit solutions for rolling Stiefel manifolds, even for a huge class of a one parameter family of pseudo-Riemannian metrics for Stiefel manifolds. This class includes many of known examples, scattered over the literature.

We strongly believe that our work will influence future research, in
particular, when rolling motions are driven by engineering
applications. To be more specific, having solutions of the kinematic
equations of rolling at hand are helpful in deriving explicit or 
closed 
formulas for differential geometric concepts such as parallel
transport and covariant derivatives, or even
to tackle control theoretic questions. These in turn will facilitate finding solutions
for interpolation,
optimization, and path planning, or other related engineering type problems.

The paper is structured as follows.
After introducing the necessary  notations,
we recall some facts on homogeneous spaces,
with emphasis on normal naturally reductive homogeneous
spaces. The Levi-Civita connection on a normal naturally reductive homogeneous space
$G / H$ is expressed in terms of vector fields on
the Lie group $G$ which have been horizontally lifted  from $G / H$ in
Subsection~\ref{subsec:levi_civiata_covariant_derivative}.
This leads in Subsection~\ref{subsec:parallel_vector_fields} to a characterization of parallel vector fields
along curves, being important for our further investigation of rolling.

 We then come to
Section~\ref{sec:rolling_intrinsic_and_extrinsic}, 
where three different notions of rolling a pseudo-Riemannian
manifold over another one of equal dimension
are introduced.
Starting with one definition of intrinsic rolling, we continue with 
 two different definitions of extrinsic rolling, the latter being
 closely related.

Although these definitions
apply to general pseudo-Riemannian manifolds,
we turn our attention to normal naturally reductive homogeneous spaces
in Section~\ref{sec:rolling_normal_naturally_reductive_homogeneous_spaces}.
The rather simple form of rolling intrinsically pseudo-Riemannian symmetric spaces
 from~\cite{JML-2023} motivates  an Ansatz,
being an obvious generalization of this rolling. Unfortunately, this does not yield
the desired result, in general.
This discussion is summarized in
Lemma~\ref{lemma:rolling_in_general_not_of_simple_form}.
In addition, it is illustrated by the example of Stiefel manifolds
equipped with $\alpha$-metrics in Subsection~\ref{subsec:no-go-stiefel}.

Afterwards, we derive the so-called kinematic equations for rolling intrinsically normal naturally reductive homogeneous spaces.
Their solutions describe the desired rolling explicitly 
if a control curve was given a priori.

In Section~\ref{sec:rolling_stiefel_manifolds}, our findings
from Subsection~\ref{subsec:kinematic_equation_intrinsic_rolling}
are applied to Stiefel manifolds.
First, we recall some facts on Stiefel manifolds endowed with $\alpha$-metrics
from the literature. Afterwards, the intrinsic rolling of Stiefel manifolds
equipped with $\alpha$-metrics is discussed by applying results from
Subsection~\ref{subsec:kinematic_equation_intrinsic_rolling}.
For a specific choice of the parameter $\alpha$,
the $\alpha$-metric on the Stiefel manifold
$\Stiefel{n}{k}$ coincides
with the metric induced by the Euclidean metric on the embedding space $\matR{n}{k}$.
Using this fact, the extrinsic rolling of Stiefel manifolds is treated in
Subsection~\ref{subsec:stiefel_extrinsic_rolling}
by extending the intrinsic rolling from
Subsection~\ref{subsec:stiefel_intrinsic_rolling}.

In Subsection~\ref{subsec:rolling_along_special_curves},
the kinematic equations describing the rolling of Stiefel manifolds are
solved explicitly where an additional assumption
is imposed on the development curve.
More precisely, an explicit formula for the extrinsic rolling of
a tangent space of $\Stiefel{n}{k}$ over $\Stiefel{n}{k}$
 is obtained provided that the development curve is the projection
of a one-parameter subgroup in $O(n) \times O(k)$ not necessarily horizontal.

Finally, in
Subsection~\ref{subsec:comparison_with_existing_literature}
we relate our results about extrinsic rolling of Stiefel manifolds to those derived
in~\cite{hueper.kleinsteuber.leite:2008}.

\subsection{Notations and Terminology}

These are some of the notations used throughout the paper.

\begin{table}[h!]
	\begin{tabular}{lcl}
			$M, N$ && smooth manifolds \\
			$T_p M$ && tangent space at $p\in M$\\
			$d_p f \colon T_p M \to T_{f(p)} N$ && tangent
                                                               map of
                                                               $f
                                                               \colon
                                                               M \to
                                                               N$ at
                                                               $p\in M$ \\
			$N_p M$ && normal space at $p\in M$\\
			$G$ && Lie group \\
			$H$ && closed subgroup of $G$ \\
			$\liealg{g}$ && Lie algebra of $G$ \\
			$\pi \colon G \to G / H$ && canonical projection \\
			$\mathcal{H}$ && horizontal bundle of $\pi \colon G \to G / H$ \\
			$\mathcal{V}$ && vertical bundle, i.e. $\mathcal{V} = \ker(d \pi)$ \\
			$\liealg{g} = \liealg{h} \oplus \liealg{p}$ && reductive decomposition \\
			$\pr_{\liealg{p}} \colon \liealg{g} \to \liealg{p}$ && projection onto $\liealg{p}$ along $\liealg{h}$ \\
			$X\at{\liealg{p}}$ && $X\at{\liealg{p}} =
                                              \pr_{\liealg{p}}(X)$ for
                                              $X\in \liealg{g}$ \\
          $X,Y$&&smooth vector fields\\
			$\nabla_X Y $ && covariant derivative of
                                         $Y$ in direction $X$ \\
			$\nabla_{\dot{\alpha}(t)} Y$ && {covariant derivative of $Y$ along curve $\alpha$.} \\
	\end{tabular} 
\end{table}

\newpage

%%%%%%%%%%%  Section 2  %%%%%%%%%%%%

\section{Normal Naturally Reductive Homogeneous Spaces}
\label{sec:reductive}
Lowercase Latin letters for elements in a Lie group
and uppercase Latin letters  for elements in the
corresponding Lie algebra are used.
For curves in the Lie algebra it will be more
convenient to use lowercase Latin letters as well.

Assume that a Lie group $G$ acts transitively from the left on a
smooth manifold $M$ by
\begin{equation*}
	\tau\colon G\times M\to M, 
	\quad
	(g,p)\mapsto  \tau (g,p)=g.p .
\end{equation*}
Then $\tau_g\colon M\to M$, defined by
\begin{equation*}
	\tau_g(p)=\tau (g,p), 
	\quad p \in M ,
\end{equation*}
is a diffeomorphism for any $g\in G$. 

Let $ \Stab(o) \subset G$ be the isotropy subgroup of a point $o\in M$, that is,  $ \Stab(o) =\{g\in G:\ g.o=\tau(g,o)=\tau_g(o)=o\}$.
The isotropy subgroup of a point in $M$ is a closed subgroup of $G$ and any two isotropy subgroups are conjugate. To simplify notations, we may denote $ \Stab(o)$ simply by $H$.
The coset manifold $G/H$ is diffeomorphic to $M$ via $g.H \mapsto g.o$,
where $g.H \in G / H$ denotes the coset defined by $g \in G$,
and we can write $M=G/H$. The manifold $M=G/H$ is called a {\it homogeneous manifold}. Denote the corresponding Lie algebras of $G$ and $H$ by $\fg$ and $\fh$, respectively.
\\

The coset manifold is said to be {\it reductive},
see e.g.~\cite[Chap. 11, Def. 21]{ONeill:1983tv}
or~\cite[Def. 23.8]{gallier.quaintance:2020},
if there exists a
subspace $\fp\subset \fg$, such that $\fg=\fp\oplus\fh$ and
$\Ad_h(X)$ for all $X\in\fp$ and $h \in H$.
This $\Ad_H$-invariance of $\fp$ implies $[\fp,\fh]\subset\fp$.
\\

Let $\pi$ denote the projection of $G$ on the coset manifold,
i.e.
\begin{equation*}
	\pi\colon G\to G/H ,
	\quad
	g \mapsto  \pi(g) = g.H .
\end{equation*}
If $e$ is the identity element in $G$,
then the map $\pi$ and its differential
\begin{equation}
	\label{diffmap}
	d_e\pi\colon T_eG=\fg\to  T_oM
\end{equation}
have the following properties.
\begin{proposition}\label{prop1}
	\begin{enumerate}
		\item[\it{1.}]  $\pi$ is a submersion;
		\item[\it{2.}] $d_e\pi(\fh)=\{0\}\subset T_oM$;
		\item[\it{3.}] $d_e\pi\at{\liealg{p}} \colon \fp\to T_oM$ is an isomorphism.
		%\item[\it{5.}] $d_e\pi$ makes one-to-one correspondence between $Ad_H$-invariant scalar products on $\fp$ and $G$-invariant metrics on $M$.
	\end{enumerate}
\end{proposition}

Consider now $M$ endowed with a pseudo-Riemannian metric
$\langle \! \langle \cdot, \cdot \rangle \! \rangle$.
We write $\langle \! \langle \cdot, \cdot \rangle \! \rangle_p$
if we want to emphasize
the value of the metric at the point $p \in M$. 
A metric tensor $\langle \! \langle \cdot, \cdot \rangle \! \rangle$
on $M$  is said to be $G$-invariant if
\begin{equation*}
	\big\langle \! \big\langle X, Y \big\rangle \! \big\rangle_p
	=
	\big\langle \! \big\langle d_p\tau_g(X),d_p\tau_g(Y) \big\rangle \! \big\rangle_{\tau_g(p)}
\end{equation*}
for all $X, Y \in T_p M$.
%\begin{equation*}
%	h_p(X,Y)=h_{g.p}(d_p\tau_g(X),d_p\tau_g(Y)),
%	\quad
%	\forall g \in G,
%\end{equation*}
%for vector fields $X,Y$ on $M$.
In other words,
the diffeomorphism $\tau_g \colon M \to M$ is an isometry.

Next we recall the definition of a pseudo-Riemannian submersion
from~\cite[Def. 44, Chap. 7]{ONeill:1983tv}.
\begin{definition}
	Let $\big(M,\langle \! \langle \cdot, \cdot \rangle \! \rangle^M \big)$
	and let $\big(N,\langle \! \langle \cdot, \cdot \rangle \! \rangle^N \big)$
	be two pseudo-Riemannian manifolds
	and $\pi \colon N\to M$ be a submersion.
	Denote by $\mathcal{V}_n = \ker(d_n \pi)$
	the vertical space at $n \in N$.
	Then $\pi$ is called a pseudo-Riemannian submersion
	if the fibers $\pi^{-1}(p)$ are pseudo-Riemannian submanifolds of $N$
	for all $p \in M$ and the maps
	$d_n \pi \at{\mathcal{H}_n} \colon \mathcal{H}_n \to T_{\pi(n)} M$
	are isometries for all $n \in N$,
	where $\mathcal{H}_n = \mathcal{V}_n^{\perp}$.
\end{definition}

%Recall that 
%a scalar product $\langle.\,,.\rangle$ in $\fg$ is said to be $\Ad_H$-invariant if 
A scalar product $\langle\cdot, \cdot \rangle$ on $\liealg{p}$ is said
to be $\Ad_H$-invariant if 
\begin{equation*}
	\big\langle \Ad_h (X),\Ad_h (Y) \big\rangle
	=
	\langle X,Y\rangle, \, \text{for}\, \text{all}\; h\in H
	\,\text{and}\, \text{for}\, \text{all}\; X,Y\in \liealg{p}.
\end{equation*}
Next we recall~\cite[Chap. 11, Prop. 22]{ONeill:1983tv}.

\begin{proposition}\label{prop11}
	By declaring the map $d_e\pi$ an isometry, there is
	one-to-one correspondence
	between $Ad_H$-invariant scalar products on $\fp$ and
	$G$-invariant metrics on $G / H$.
\end{proposition}

\begin{definition}
	\label{definition:naturally_reductive_space}
	A coset manifold $M=G/H$ is called a naturally  reductive space, if
	\begin{enumerate}
		\item[\it{(1)}] $M=G/H$ is reductive;
		\item[\it{(2)}] $M$ carries a $G$-invariant metric;
		\item[\it{(3)}] If $\langle \cdot, \cdot \rangle$ denotes the
		$Ad_H$-invariant scalar product on $\fp$ corresponding
		to the $G$-invariant metric (described in Proposition~\ref{prop11}), 
		then it has to satisfy
		\begin{equation*}
			\big\langle [X,Y]\at{\fp},Z \big\rangle 
			=
			\big\langle X,[Y,Z]\at{\fp} \big\rangle,
			\, \text{for}\, \text{all}\; X,Y,Z\in\fp .
		\end{equation*}
	\end{enumerate}
\end{definition}
Naturally reductive homogeneous spaces are complete,
see~\cite[Chap. 11, p. 313]{ONeill:1983tv}.

Next, we introduce the notion of (pseudo-Riemannian) normal naturally reductive homogeneous space. 
This definition is a slight generalization of the
 homogeneous spaces which are considered in~\cite[Prop. 23.29]{gallier.quaintance:2020}.

\begin{definition}
	\label{definition:normal_naturally_reductive_spaces}
	{\bf{Normal Naturally Reductive Spaces}}
	Let $G$ be a Lie group equipped with a bi-invariant metric
	and denote by $\langle \cdot , \cdot \rangle$ the
	corresponding $\Ad_G$-invariant scalar product on its Lie algebra
	$\liealg{g}$. Moreover, let $H \subset G$ be a closed subgroup and
	denote its Lie algebra by $\liealg{h} \subset \liealg{g}$.
	If the orthogonal complement $\liealg{p} = \liealg{h}^{\perp}$ with respect to $\langle \cdot, \cdot \rangle$ is non-degenerated,
	we call $G / H$ equipped with the  $G$-invariant metric that turns $\pi \colon G \to G / H$ into a pseudo-Riemannian submersion
	a (pseudo-Riemannian) normal naturally reductive homogeneous space
	with reductive decomposition $\liealg{g} = \liealg{h} \oplus \liealg{p}$.
\end{definition}
By a trivial adaption of the proof of~\cite[Prop. 23.29]{gallier.quaintance:2020}, 
we show that normal naturally reductive spaces are naturally reductive.

\begin{lemma}
	\label{lemma:normal_naturally_reductive_implies_reductive}
	Let $G / H$ be normal naturally reductive. Then $G / H$ is naturally reductive.
	\begin{proof}
		Let $X \in \liealg{p} = \liealg{h}^{\perp}$.
		Then $\langle Y , X \rangle = 0$ for all $Y \in \liealg{h}$.
		The $\Ad_G$ invariance of $\langle \cdot, \cdot \rangle$ implies
		\begin{equation}
			\big\langle \Ad_h(X), \Ad_h(Y) \big\rangle = 0, 
			\quad h \in H.
		\end{equation}
		Since $\Ad_h \colon \liealg{h} \to \liealg{h}$ is an isomorphism,
		this implies $\langle \Ad_h(X), \widehat{Y} \rangle = 0$
		for $h \in H$ and all $\widehat{Y} \in \liealg{h}$ proving
		$\Ad_h(X) \in \liealg{p}$ for $h \in H$, i.e.
		$\Ad_h(\liealg{p}) \subset \liealg{p}$ for $h \in H$.
		In addition, 
		$\liealg{g} = \liealg{h} \oplus \liealg{h}^{\perp} 
		= \liealg{h} \oplus \liealg{p}$ 
		is fulfilled since
		$\liealg{h}^{\perp}$ is assumed to be non-degenerated.
		Thus $G / H$ is a reductive homogeneous space.
		
		In order to show that $G / H$ is naturally reductive, we compute
		for $X, Y, Z \in \liealg{p}$
		\begin{equation}
			\begin{split}
				0 
				&= 
				\tfrac{\D}{\D t} \langle X, Z \rangle \at{t = 0} \\
				&= 
				\tfrac{\D}{\D t} 
				\big\langle \Ad_{\exp(t Y)}(X), \Ad_{\exp(t Y)}(Z)  \big\rangle \at{t = 0} \\
				&=
				\big\langle [Y, X], Z \big\rangle + \big\langle X, [Y, Z] \big\rangle \\
				&=
				- \big\langle [X, Y], Z \big\rangle + \big\langle X, [Y, Z] \big\rangle,
			\end{split}
		\end{equation}
		where we have used the $\Ad_G$-invariance of
		$\langle \cdot, \cdot \rangle$.
		Finally, since $\liealg{p} = \liealg{h}^{\perp}$, the last identity implies 
		\begin{equation}
			\big\langle [X, Y]\at{\liealg{p}}, Z \big\rangle 
			=
			\big\langle X, [Y, Z]\at{\liealg{p}} \big\rangle,
			\quad 
			X, Y, Z \in \liealg{p},
		\end{equation}
		i.e. $G / H$ is a naturally reductive homogeneous space.
	\end{proof}
\end{lemma}

%We will use the following construction of the (pseudo)-Riemannian submersion.

%(1)\ If $M=G/H$ is a naturally reductive space, then we can extend the scalar product from $\fp$ to the entire Lie algebra $\fg$, such that the direct sum $\fp\oplus\fh$ becomes orthogonal. We then use left translations and get left-invariant distributions corresponding to $\fp$ and $\fh$, which will be orthogonal in any point on $G$. By the left translation, we also get a left-invariant metric tensor on $G$. Under this construction
%$\pi\colon G\to M$ becomes a Riemannian submersion~\cite[page 312]{ONeill:1983tv}. 

%(2)\ 
%The Lie algebra $\fg=\fp\oplus\fh$  carries an $\Ad_G$-invariant scalar product $\langle.\,,.\rangle$ making the direct sum $\fp\oplus\fh$ orthogonal. By making use of left translations we obtain left invariant horizontal and vertical distributions $\mathcal H\oplus\mathcal V$ as well as bi-invariant Riemannian metric tensor on $G$. The distributions will be orthogonal at each point. We also have that $M=G/H$ is a Riemannian manifold with a $G$-invariant metric as in Proposition~\ref{prop11}. Then the map $\pi\colon G\to M$ is a Riemannian submersion.

Let $G / H$ be a normal naturally reductive space.
Then, by definition, the map $\pi \colon G \to G / H$ is a pseudo-Riemannian 
submersion.
Obviously, the vertical bundle and horizontal bundle are given by 
\begin{equation*}
	\mathcal{V}_g 
	=
	\ker(d_g \pi) = \big( d_e L_g \big) \liealg{h}
	\quad \text{ and } \quad
	\mathcal{H}_g 
	=
	\mathcal{V}_g^{\perp} 
	=
	\big( d_e L_g \big) \liealg{p},
\end{equation*}
for $g \in G$, respectively.
From an algebraic point of view, the reductive decomposition has the following properties:
\begin{equation*}
	\mathfrak g= \fp\oplus_{\perp}\fh, \quad
	[\fh,\fh]\subset\fh,\quad
	[\fp,\fh]\subset\fp .
\end{equation*}

%From an algebraic point of view, the reductive decomposition has the following properties:
%$$
%\mathfrak g= \fp\oplus_{\perp}\fh, \quad
%[\fh,\fh]\subset\fh,\quad
%[\fp,\fh]\subset\fp,\quad
%\fh\subset [\fp,\fp]
%$$

We end this preliminary section by commenting on the regularity of curves.
%\begin{notation}
	Throughout this text, for simplicity, if not indicated otherwise,
	a curve $c \colon I \to M$ on a manifold $M$ is assumed to be smooth.
	However, we point out that many results
	can be generalized by requiring less regularity.
%\end{notation}

\subsection{Levi-Civita Connection and Covariant Derivative}
\label{subsec:levi_civiata_covariant_derivative}
We first set some notations. The  Levi-Civita connections on
$M=G / H$ and on $G$ will be denoted by
$\nabla^M$ and $\nabla^G$, respectively.
In cases when it is clear from the context, we may use
simply $\nabla$ for both.
If $Y$ is a vector field on $M=G / H$, we denote by
$\widetilde Y\in \Secinfty(T G)$ its horizontal lift to $G$.
Correspondingly, if $\alpha\colon I\to M$ is a curve in $M$ and $r\colon I\to G$ is
a lift of $\alpha$ to $G$, we write
$\nabla_{\dot\alpha (t)}Y\in T_{\alpha (t)}M$ for the covariant
derivative of $Y$ along $\alpha$, and
$\widetilde{\nabla_{\dot\alpha (t)}Y}\at{r(t)}$ for the horizontal
lift of $\nabla_{\dot\alpha (t)}Y$ to
$\mathcal{H}_{r(t)}\subset T_{r(t)}G$.

In the sequel, the lift of $\alpha$ to $G$ will be denoted by $q$ instead of $r$
if it is considered to be horizontal.
 For $g \in G$ denote by $\pr_{\mathcal H_g} \colon T_g G \to \mathcal{H}_g$ the
projection onto the horizontal bundle, explicitly given by 
\begin{equation}
	\label{equation:projection_onto_horizontal_bundle_formula}
	\pr_{\mathcal{H}_g} 
	=
	\big(d_e L_g \big) \circ \pr_{\liealg{p}} \circ \big(d_e L_g \big)^{-1}	.
\end{equation}

%In the next lemma the lift of $\alpha$ to $G$ is considered to be horizontal and will be denoted by $q$ instead of $r$.

\begin{lemma}
	\label{lemma:covariant_derivative_vector_fields}
	Let $G / H$ be a normal naturally reductive homogeneous space
	and let $X, Y$ be vector fields on $G / H$.
	Denote by $\widetilde{X}$ and $\widetilde{Y}$
	the horizontal lifts of $X$ and $Y$, respectively.
	Moreover, let $\{A_1, \ldots, A_k \mid i = 1, \ldots, k \}$ be a basis of $\liealg{p}$ 
	and denote by $\overline{A}_1, \ldots, \overline{A}_k$ the
	corresponding left invariant vector fields defined by $\overline{A}_i(g) = d_e L_g A_i$
	for $g \in G$. 
	Expanding $\widetilde{X} = \sum_{i  = 1}^{k} x_i \overline{A}_i$
	and $\widetilde{Y} = \sum_{j = 1}^k y_j \overline{A}_j$
	with smooth functions $x_i, y_j \colon G \to \field{R}$,
	we obtain for the Levi-Civita covariant derivative on $G / H$
	for $g \in G$
	\begin{equation}
		\label{equation:lemma_covariant_derivative_vector_fields_assertion_1}
		\begin{split}
			(\nabla_X^M Y )(\pi(g))
			&=
			d_{g} \pi \Big( \sum_{j = 1}^k \big(\widetilde{X}(y_j) \big)(g) \overline{A}_j(g) \\
			&\quad+ 
			\pr_{\mathcal H_g}
			\tfrac{1}{2}
			\sum_{i, j = 1}^k
			x_i(g) y_j(g) \Big[\overline{A}_i, \overline{A}_j\Big](g) \Big), 
		\end{split}
	\end{equation}
	 or, equivalently, 
	\begin{equation}
		\label{equation:lemma_covariant_derivative_vector_fields_assertion_2}
		\widetilde{\nabla_X^M Y}\at{g}
		=
		\sum_{j = 1}^k \big(\widetilde{X}(y_j) \big)(g) \overline{A}_j(g)
		+ 
		\tfrac{1}{2}
		\sum_{i, j = 1}^k
		x_i(g) y_j(g) \overline{[A_i, A_j]\at{\liealg{p}}}(g) .
	\end{equation}
	\begin{proof}
		Since the metric is bi-invariant, it follows that for left-invariant
		vector fields $V,W $ on $G$, see~\cite[page 304]{ONeill:1983tv},
		\begin{equation}\label{eq:Nabla}
			\nabla^G_VW=\tfrac{1}{2}[V,W]
		\end{equation}  
		holds.
		Since $G / H$ is a normal naturally reductive space, the map
		$\pi \colon G \to G / H$ is a pseudo-Riemannian submersion.
		Let $X$, $Y$ be
		vector fields on $M$, and $\widetilde X$, $\widetilde Y$ their horizontal
		lifts to $G$. 
		We recall that the Levi-Civita connections on $M$ and on $G$ are related by,
		see~\cite[Lemma 45, Chapter 7]{ONeill:1983tv},
		\begin{equation}
			\label{equation:Levi_civita_Lie_group_biinvariant}
			\nabla^M_XY
			=
			d_g \pi\Big(\pr_{\mathcal H_g}\nabla^G_{\widetilde X}\widetilde Y \Big) .
		\end{equation}
		Expanding  the horizontal lifts $\widetilde X$ and $\widetilde Y$  in terms of the left-invariant frame field $\{\overline A_1,\ldots, \overline A_k\}$, i.e.,  
		\begin{equation}
			\widetilde X=\sum_{i=1}^kx_i\overline A_i,
			\quad \widetilde Y=\sum_{j=1}^ky_j\overline A_j,
		\end{equation}
		we have
		\begin{equation}
			\label{equation:lemma_covariant_derivative_vector_fields_lie_group_projected}
			\begin{split}
				\nabla^G_{\widetilde X} \widetilde Y
				=
				\nabla^G_{\widetilde X}  \Big( \sum_{j=1}^k y_j\overline A_j \Big)
				=
				\sum_{j=1}^k \big(\widetilde X(y_j) \big)\overline A_j
				+
				\tfrac{1}{2}\sum_{i,j=1}^kx_iy_j\Big[\overline A_i,\overline A_j\Big].
			\end{split}
		\end{equation}
		Projecting to $\mathcal H_g$, and taking into consideration that the first
		term in the last equality belongs to $\mathcal H_g$, we obtain 
		\begin{equation}
			\label{equation:lemma_covariant_derivative_vector_fields_lie_group_projected_horizontal_bundle}
			\pr_{\mathcal H_g}\nabla^G_{\widetilde X}\widetilde Y
			=
			\sum_{j=1}^k \big(\widetilde X(y_j) \big)\overline A_j 
			+
			\pr_{\mathcal H_g}\tfrac{1}{2}\sum_{i,j=1}^kx_iy_j\Big[\overline A_i,\overline A_j\Big] .
		\end{equation}
		Identity~\eqref{equation:lemma_covariant_derivative_vector_fields_lie_group_projected_horizontal_bundle},
		together
		with~\eqref{equation:Levi_civita_Lie_group_biinvariant},
		gives~\eqref{equation:lemma_covariant_derivative_vector_fields_assertion_1}.
		Clearly, by using~\eqref{equation:projection_onto_horizontal_bundle_formula},
		one has $\pr_{\mathcal H_g} ( [\overline{A}_i, \overline{A}_j] )(g) = \overline{[A_i, A_j]\at{\liealg{p}}}(g)$.
		Hence, ~\eqref{equation:lemma_covariant_derivative_vector_fields_assertion_1}
		is equivalent to~\eqref{equation:lemma_covariant_derivative_vector_fields_assertion_2},
		as the vector field
		from~\eqref{equation:lemma_covariant_derivative_vector_fields_lie_group_projected}
		on $G$ is horizontal and $\pi$-related to $\nabla^M_X Y$
		by~\eqref{equation:lemma_covariant_derivative_vector_fields_assertion_1}.
	\end{proof}
\end{lemma}
Lemma~\ref{lemma:covariant_derivative_vector_fields}
yields an expression for the Levi-Civita covariant
derivative on $G / H$
in terms of horizontally lifted vector fields on $G$.
This expression allows for determining the covariant derivative
of vector fields along a curve in $G / H$ in terms of
horizontally lifted vector fields
along a horizontal lift of the curve, as well.
As preparation we comment on the domain of horizontal lifts.

\begin{remark}
	Let $\alpha \colon I \to G / H$ be a curve
	on a normal naturally reductive space.
	The horizontal lift $q \colon I \to G$ is indeed defined on the same interval
	as $\alpha$. 
	This can be shown by exploiting that $\mathcal{H} \subset T G$ 
	defines a principal connection which is known to be complete.
\end{remark}

\begin{lemma}
	\label{lemma:covariant_derivatives_along_curves_horizontal_lift}
	Let $M = G / H$ be a normal naturally reductive homogeneous space,
	$\alpha\colon I\to M$  a curve and $Y$ a vector field along $\alpha$.
	Let $q\colon I\to G$ be a horizontal lift of $\alpha$ and $\widetilde Y$ a
	horizontal lift of $Y$ along $q$. Then 
	\begin{eqnarray}
		\label{equation:lemma_covariant_derivatives_along_curves_horizontal_lift_1}
		\nabla^M_{\dot\alpha (t)}Y(t)
		&=&
		d_{q(t)}\pi\Big(\sum_{j=1}^k\tfrac{\D y_j(t)}{\D t} \overline A_j(t)\Big)\nonumber
		\\
		&+&
		d_{q(t)}\pi\Big(\pr_{\mathcal H_{q(t)}}\tfrac{1}{2}\sum_{i,j=1}^k x_i(t)y_j(t)\Big[\overline A_i(t),\overline A_j(t)\Big]\Big),
	\end{eqnarray}
	or equivalently
	\begin{equation}
		\label{equation:lemma_covariant_derivatives_along_curves_horizontal_lift_2}
		\begin{split}
			\widetilde{\nabla^M_{\dot\alpha (t)}Y}\at{q(t)}
			&=
			\sum_{j=1}^k\tfrac{\D y_j(t)}{\D t} \overline A_j(t)
			+
			\tfrac{1}{2}\sum_{i,j=1}^k x_i(t)y_j(t)\overline{[A_i, A_j]\at{\liealg{p}}}(t) ,
		\end{split}
	\end{equation}
	where $\{A_1,\ldots,A_k\}$ is a basis of
	$\fp$, $\overline{A}_i$ denotes the left-invariant vector field
	corresponding to $A_i$ for $i = 1, \ldots k$
	and we write $\overline{A}_i(t) = \overline{A}_i(q(t))$
	for short.
	The functions $x_i, y_j \colon I \to \field{R}$ are defined by 
	$\dot{q}(t) = \sum_{i  = 1}^k x_i(t) \overline{A}_i(t)$
	and $\widehat{Y}(t) = \sum_{j = 1}^k y_j \overline{A}_j(t)$.
	\begin{proof}
		Let $t \in I$.
		We extend the vector field $\dot{\alpha}(t)$ and $Y(t)$ locally
		to vector fields $\widehat{X}$ and $\widehat{Y}$,
		respectively,
		defined on an open neighborhood of $\alpha(t)$ in $G / H$.
		The proof of~\cite[Thm. 4.24]{Lee:2018aa} shows
			that such an extension is always possible.
		Moreover, we denote by $\widetilde{\widehat{X}}$ and $\widetilde{\widehat{Y}}$
		the horizontal lifts of $\widehat{X}$ and $\widehat{Y}$, respectively.
		These vector fields are expanded as 
		$\widetilde{\widehat{X}} = \sum_{i = 1}^k \widehat{x}_i \overline{A}_i$
		and
		$\widetilde{\widehat{Y}} = \sum_{j = 1}^k \widehat{y}_j \overline{A}_j$
		with uniquely locally defined functions $\widehat{x}_i, \widehat{y}_j$ on $G$.
		Clearly, these functions fulfill $\widehat{x}_i\big(q(t)\big) = x_i(t)$
		and $\widehat{y}_j(q(t)) = y_j(t)$ whenever both sides are defined.
		In addition, $\widetilde{\widehat{X}}\big(q(t) \big) = \dot{q}(t)$ and 
		$\widetilde{\widehat{Y}}\big (q(t) \big) = \widetilde{{Y}}(t)$ holds.
		By using Lemma~\ref{lemma:covariant_derivative_vector_fields}, 
		we compute
		\begin{equation*}
			\begin{split}
				(\nabla_X^M Y)(\pi(q(t))
				&=
				d_{q(t)} \pi \Big(
				\sum_{j=1}^k\Big(\widetilde{\widehat{X}}(\widehat{y}_j)\Big) \big (q(t) \big)\overline A_j \big (q(t) \big) \Big) \\
				&\quad +
				d_{q(t)} \pi \Big(
 				\pr_{\mathcal H_{q(t)}}\tfrac{1}{2}\sum_{i,j=1}^k \widehat{x}_i \big (q(t) \big) \widehat{y}_j \big (q(t) \big)\Big[\overline A_i,\overline A_j\Big]\big(q(t) \big)
 				\Big)
				\\
				&=
				d_{q(t)} \pi \Big(\sum_{j=1}^k\tfrac{\D y_j(t)}{\D t} \overline A_j(t) \Big) \\
				&\quad+
				d_{q(t)} \pi \Big(
				\pr_{\mathcal H_{q(t)}}\tfrac{1}{2}\sum_{i,j=1}^k x_i(t)y_j(t)\Big[\overline A_i,\overline A_j\Big](t) \Big)
			\end{split}
		\end{equation*}
		showing~\eqref{equation:lemma_covariant_derivatives_along_curves_horizontal_lift_1}.
		Clearly, this is equivalent to~\eqref{equation:lemma_covariant_derivatives_along_curves_horizontal_lift_2}
		by Lemma~\ref{lemma:covariant_derivative_vector_fields}.
	\end{proof}
\end{lemma}

\begin{remark}If $M=G/H$ is a symmetric space,  then  $[\fp,\fp ]\subset\fh$, and consequently the last summand in 
%formula~\eqref{eq:NN0}
formula~\eqref{equation:lemma_covariant_derivatives_along_curves_horizontal_lift_1}
vanishes.
So, taking into consideration that, in this case, $\nabla^G_{\dot q(t)} \widetilde Y(t)=\sum_{j=1}^k\tfrac{\D y_j(q(t))}{\D t} \overline A_j$, the
%identity $~\eqref{eq:NN0}$
identity~\eqref{equation:lemma_covariant_derivatives_along_curves_horizontal_lift_1}
reduces to 
$$
\nabla^M_{\dot\alpha (t)}Y(t)=d_{q(t)}\pi\Big(\nabla ^G_{\dot{q}(t)}\widetilde Y(t)\Big),
$$
which shows that, in case of a symmetric space, if $Y$ is a
parallel vector field along $\alpha(t)\in M$, its horizontal
lift $\widetilde Y$ is actually a parallel vector field along
the horizontal lift $q(t)\in G$ of $\alpha(t)$.
\end{remark}

As we will see below,
for non symmetric spaces the presence of the second term in~\eqref{equation:lemma_covariant_derivatives_along_curves_horizontal_lift_1} reveals that the horizontal lift $q(t)\in G$ is not a good candidate for the property of preserving parallel vector fields. In the next section we modify the "horizontal lift" in order to overcome this problem.

\subsection{Parallel Vector Fields}
\label{subsec:parallel_vector_fields}

\begin{lemma}
	\label{lemma:horizontal_lift_cov_der_arbitrary_curve}
	Let $M = G / H$ be a normal naturally
	reductive homogeneous space, $\alpha\colon I\to M$  a
	curve and  $q\colon [0,T]\to G$ a horizontal lift of $\alpha$. 
	Moreover, let $s \colon I \to H$ and define the
	curve $r \colon I \to G$ by $r(t) = q(t) s(t)$.
	Let $Z \colon I \to T M$ be a vector field along $\alpha$
	and denote by $\widetilde{Z} \colon I \to \mathcal{H}$ its horizontal
	lift along $r$.
	Then, the horizontal lift of
	$\nabla_{\dot{\alpha}(t)} Z \colon I \to T M$ along
	$r(t)$  is given by
	\begin{equation}
	\begin{split}
		\widetilde{\nabla_{\dot{\alpha}(t)} Z}\at{r(t)} 
		=&
		\sum_{j=1}^k\dot{z}_j(t) \overline{A}_j(r(t)) \\
		&+ \sum_{i,j=1}^k\tfrac{1}{2} x_i(t) z_j(t) 
		\pr_{\mathcal H_{r(t)}}\Big( \overline{[ \Ad_{s(t)^{-1}}(A_i), A_j ]}\big(r(t)\big) \Big).
	\end{split}
	\end{equation}
	Here we expanded
	$x(t) 
	= \big(d_e L_{q(t)} \big)^{-1} \dot{q}(t) 
	= \sum_{i=1}^k x_i(t) A_i \in \liealg{p}$
	and accordingly
	$z(t) 
	= \big(d_e L_{r(t)} \big)^{-1}  \widetilde{Z}(t) 
	= \sum_{i=1}^k z_i(t) A_i \in \liealg{p}$.	
	\begin{proof}
		Let $X, Z \in \Secinfty\big(T (G / H) \big)$ be vector fields with
		horizontal lifts $\widetilde{X},  \widetilde{Z} \in \Secinfty(T G)$
		and expand them by a left invariant frame
		$\overline{A}_1, \ldots, \overline{A}_k$ of the horizontal bundle of $G \to G / H$, i.e.
		$\widetilde{X} =  \sum_{i=1}^kx_i \overline{A}_i$ and
		$\widetilde{Z} = \sum_{j=1}^k z_j \overline{A}_j$.
		Then, by Lemma~\ref{lemma:covariant_derivatives_along_curves_horizontal_lift}, the Levi-Civita connection on $G / H$ can
		be expressed in terms of horizontal lifts by
		\begin{equation}
			\label{equation:lemma_horizontal_lift_cov_der_arbitrary_curve_horizontal_lift_vf}
			\widetilde{\nabla_X Z} = \sum_{j=1}^k\widetilde{X}(z_j ) \overline{A}_j + \tfrac{1}{2} \sum_{i,j=1}^kx_i z_j \overline{[A_i, A_j]\at{\liealg{p}}}.
		\end{equation} 
				Now, consider the curve $r(t) = q(t)
                                s(t)$ being a lift of $\alpha(t)$.
		A simple computation shows that 
		\begin{equation}\label{aaaa}
			\big(d_e L_{r(t)} \big)^{-1} \dot{r}(t)
			=
			\Ad_{s(t)^{-1}} \big (x(t) \big) + y(t),
		\end{equation}
		 where $y(t) := \big(d_e L_{s(t)} \big)^{-1} \dot{s}(t) \in \liealg{h}$.
		Thus, using \eqref{aaaa} and  $\pi \big(r(t) \big)=\alpha(t)$, we have
		\begin{equation}
			\label{equation:lemma:horizontal_lift_cov_der_arbitrary_curve_projection_r}
			\begin{split}
				\dot{\alpha}(t) 
				&=
				d_{r(t)} \pi \dot{r}(t) \\
				&= 
				\big( d_{r(t)} \pi \circ d_e L_{r(t)}\big) 
				\big( \Ad_{s(t)^{-1}} \big (x(t) \big) + y(t) \big) \\
				&= 
				\big( d_{r(t)} \pi \circ d_e L_{r(t)}\big) 
				\Big( \Ad_{s(t)^{-1}}\big (x(t) \big) \Big).
			\end{split}
		\end{equation}
		Here the last equality follows from the definition of
		the horizontal bundle.
		By extending $\dot{\alpha}(t)$ locally to a	vector field $X$
		on $G / H$, the horizontal
		lift $\widetilde{X}$ of $X$ satisfies 
		$\widetilde{X}\big (r(t) \big) = d_e L_{r(t)} \Big(\Ad_{s(t)^{-1}}\big (x(t) \big) \Big)$
		by~\eqref{equation:lemma:horizontal_lift_cov_der_arbitrary_curve_projection_r}.
		Moreover,  the vector field $Z$ along $\alpha$
		can be extended locally to a  vector field
		$\widehat{Z}$ on $G / H$, defined on an open neighbourhood
		of $\alpha$.
		Denote by $\widetilde{\widehat{Z}}$ the horizontal lift of
		$\widehat{Z}$. Then
		$\widetilde{\widehat{Z}}\big (r(t) \big) = \widetilde{Z}(t)$ is fulfilled.
		By~\cite[Thm. 4.24]{Lee:2018aa}, we have
		\begin{equation}
			\begin{split}
				\widetilde{\nabla_{\dot{\alpha}(t)} Z}\at{r(t)} 
				&=
				\widetilde{\nabla_{\widetilde{X}} \widetilde{\widehat{Z}}}\at{r(t)}
			\end{split}.
		\end{equation}
		The desired result follows by
	exploiting~\eqref{equation:lemma_horizontal_lift_cov_der_arbitrary_curve_horizontal_lift_vf}, similarly to what was done in the proof of Lemma~\ref{lemma:covariant_derivatives_along_curves_horizontal_lift}.
	\end{proof}
\end{lemma}
\begin{corollary}
	\label{corollary:parallel_vf}
	The vector field $Z \colon I \to T ( G/H)$ 
	along $\alpha \colon I \to G /H$
	is parallel along $\alpha$ iff its horizontal lift $\widetilde{Z}$ along
	$r(t) = q(t) s(t) \in G$, defined as in  Lemma~\ref{lemma:horizontal_lift_cov_der_arbitrary_curve} by $z(t) = (d_e L_{r(t)})^{-1}  \widetilde{Z}(t) = \sum_{i=1}^k z_i(t) A_i \in \liealg{p}$,
	satisfies
	\begin{equation}
		\dot{z}(t) 
		=
		- \tfrac{1}{2} \pr_{\liealg{p}} \Big( \big[\Ad_{s(t)^{-1}}\big (x(t) \big), z(t) \big] \Big)
	\end{equation}
	for all $t \in I$, where $x(t)
	= \big(d_e L_{q(t)} \big)^{-1} \dot{q}(t) 
	=\sum_{i=1}^k x_i(t) A_i \in \liealg{p}$.
	\begin{proof}
		Lemma~\ref{lemma:horizontal_lift_cov_der_arbitrary_curve}
                already implies the statement by applying the linear isomorphism
		$\big(d_{r(t)} \pi \circ d_e L_{r(t)}\big)^{-1}$ 
		to both sides of
		$0 = \widetilde{\nabla_{\dot{\alpha}(t)} Z}\at{r(t)}$.
	\end{proof}
\end{corollary}

 When  $s(t) = e$, for $t \in I$, Corollary~\ref{corollary:parallel_vf}
 also gives the following characterization of parallel vector fields. 
	
\begin{corollary}
	\label{corollary:parallel_vf_horizontal}
	The vector field $Z \colon I \to T ( G/H)$ 
	along $\alpha \colon I \to G /H$
	 with horizontal lift $q \colon I \to G$
	is parallel along $\alpha$ iff its horizontal lift $\widetilde{Z}$ along $q$ 
	fulfills the ODE
\begin{equation}
		\dot{z}(t) = - \tfrac{1}{2} \pr_{\liealg{p}} \big( [x(t), z(t)] \big)
	\end{equation}
	for all $t \in I$,
	where $x(t) = \big(d_e L_{q(t)} \big)^{-1} \dot{q}(t) \in \liealg{p}$
	and 
	\begin{equation}
		z(t) = \big(d_e L_{q(t)} \big)^{-1}  \circ \big(d_{q(t)} \pi\at{\mathcal{H}_{q(t)}} \big)^{-1} Z(t) \in \liealg{p}.
	\end{equation}
		
\end{corollary}

\section{Intrinsic and Extrinsic Formulation of Rolling}
\label{sec:rolling_intrinsic_and_extrinsic}

The goal of this section is to introduce the notation of rolling
a pseudo-Riemannian manifold over another one. 

In the following definitions, it is assumed that the pseudo-Riemannian
manifolds $(M,g)$ and $(\M,\widehat g)$ are of equal dimension and $g$
and $\widehat g$ have the same signature.  
%To compare the two notions of rolling, the isometries considered in the next two definitions are restricted to oriented isometries. We call an isometry $A\colon T_mM\to T_{\widehat m}\M$ oriented if it preserves the orientation of the positive definite and the negative definite subspaces of $T_mM$ and $T_{\widehat m}\M$.

\begin{definition} \label{definition:intrinsic_rolling} 
	{\bf{Intrinsic rolling.}}
	A curve $\alpha (t)$ on $M$ is said to roll on a curve
	$\widehat{\alpha}(t)$ on $\M$ intrinsically if there exists an
	 isometry $A(t): T_{\alpha (t)}M\to T_{\widehat{\alpha} (t)}\M$ satisfying the following conditions: 
	 \begin{enumerate}
	 	\item
	 	No-slip condition:
	 	$\dot{\widehat{\alpha}}(t)=A(t)\dot\alpha (t)$; 
	 	\item
	 	No-twist condition:
	 	$A(t)X(t)$ is a parallel vector field in $\M$ along $\widehat{\alpha} (t)$  iff $X(t)$ is a parallel vector field in $M$ along $\alpha (t)$.
	 \end{enumerate}
 	
	The triple  $\big (\alpha(t),\widehat{\alpha}(t),A(t) \big)$ is called a \it{rolling}
	(of $M$ over $\widehat{M}$).
	The curve $\alpha$ is called rolling curve while $\widehat{\alpha}$
	is called development curve.
\end{definition}

The next definition of extrinsic rolling is motivated by the
description of extrinsic rolling in terms of bundles,
see~\cite[Def. 2]{molina.grong.markina.leite:2012} and~\cite[Def. 3]{markina.leite:2016}.
\begin{definition}
	\label{definition:rolling_extrinsic_1}
	{\bf{Extrinsic rolling (I)}}
	Let $M$ and $\widehat{M}$ be isometrically embedded into the same pseudo-Euclidean vector space $V$.
	A quadruple $\big (\alpha(t), \widehat{\alpha}(t), A(t), C(t)
        \big)$ is called an extrinsic rolling
	(of $M$ over $\widehat{M}$),
	where $\alpha \colon I \to M$ and $\widehat{\alpha} \colon I \to \widehat{M}$
	are curves, $A(t) \colon T_{\alpha(t)} M \to T_{\widehat{\alpha}(t)} \widehat{M}$
	and $C(t) \colon N_{\alpha(t)} M \to N_{\widehat{\alpha}(t)} \widehat{M}$
	are isometries of the tangent and normal spaces, if the following conditions hold:
	 \begin{enumerate}
		\item
		No-slip condition:
		$\dot{\widehat{\alpha}}(t)=A(t)\dot\alpha (t)$; 
		\item
		No-twist condition (tangential part):
		$A(t)X(t)$ is a parallel vector field in $\M$ along $\widehat{\alpha} (t)$  if and only if $X(t)$ is a parallel vector field in $M$ along $\alpha (t)$;
		\item
		No-twist condition (normal part):
		$C(t)Z(t)$ is a normal parallel vector field in $\M$ along
		$\widehat{\alpha} (t)$  iff $Z(t)$ is a normal
		parallel vector field in $M$ along $\alpha (t)$.
	\end{enumerate}
	As in the intrinsic case, the curve $\alpha$ is called rolling curve while
	$\widehat{\alpha}$
	is called development curve.
\end{definition}

Alternatively, we define extrinsic rolling as reformulation of
a slightly generalized version of~\cite[Def. 1]{markina.leite:2016}.

\begin{definition}
	{\bf{Extrinsic rolling (II)}}
	\label{definition:extrinsic_rolling_2}
	Let $M$ and $\widehat{M}$ be isometrically embedded into
	the same pseudo-Euclidean vector space $V$.
	A curve $(\alpha, E) \colon I \to M \times E(V)$,
	where $E(V) = O(V) \ltimes V$ denotes the pseudo-Euclidean group of $V$,
	is said to be an extrinsic rolling 
	if the following conditions are satisfied:
	\begin{enumerate}
		\item
		\label{item:def_extrinsic_2_claim_1}
		$\widehat{\alpha}(t) := E(t) \alpha(t) \in \widehat{M}$;
		\item
		\label{item:def_extrinsic_2_claim_2}
		$d_{\alpha(t)} E(t) (T_{\alpha(t)} M) = T_{\widehat{\alpha}(t)} \widehat{M}$;
		\item
		\label{item:def_extrinsic_2_claim_3}
		No-slip condition:
		$\dot{\widehat{\alpha}}(t) = d_{\alpha(t)} E(t) \dot{\alpha}(t)$
		\item
		\label{item:def_extrinsic_2_claim_4}
		No-twist condition (tangential part): $d_{\alpha(t)} E(t) X(t)$ is parallel along $\widehat{\alpha}$ iff $X$ is parallel along $\alpha$; 
		\item
		\label{item:def_extrinsic_2_claim_5}
		No-twist condition (normal part): 
		$d_{\alpha(t)} E(t) Z(t)$ is normal parallel along $\widehat{\alpha}$ iff
		$Z$ is normal parallel along $\alpha$. 
	\end{enumerate}
	The curve $\alpha$ is called rolling curve and the $\widehat{\alpha}$ 
	is the development curve.
\end{definition}

\begin{remark}
	The discussion in~\cite[Sec. 3]{JML-2023}
	reveals that a rolling in the sense of Definition~\ref{definition:extrinsic_rolling_2}
	is closely related to the classical definition of rolling in~\cite[Ap. B, Def. 1.1]{sharpe:1997}.
	Indeed, the conditions
        Definition~\ref{definition:extrinsic_rolling_2} and 
	Claims~\ref{item:def_extrinsic_2_claim_1}-\ref{item:def_extrinsic_2_claim_5}
	are equivalent to the conditions from~\cite[Def. 1.1]{sharpe:1997}.
	Thus the essential difference between Definition~\ref{definition:extrinsic_rolling_2}
	and~\cite[Def. 1.1]{sharpe:1997}
	is that the rolling curve is already part of the Definition.
	This is motivated by~\cite[Ex. 1]{molina.grong.markina.leite:2012}.
\end{remark}

Motivated by~\cite[Prop. 3]{JML-2023}, we 
relate
the two different notions of extrinsic rolling from Definition~\ref{definition:rolling_extrinsic_1}
and Definition~\ref{definition:extrinsic_rolling_2}.

\begin{proposition}
	\label{proposition:relation_extrinsic_rolling}
	Let $\big(\alpha(t), \widehat{\alpha}(t), A(t), C(t) \big)$ be an extrinsic rolling in the sense of Definition~\ref{definition:rolling_extrinsic_1}.
	Then the curve
	$g(t) = \big (\alpha(t), \big (R(t) , s(t) \big) \big) \in M \times E(V)$,
	where
	\begin{equation}
		\label{equation:proposition_relation_extrinsic_rolling_curve_in_M_E_V}
		\begin{split}
			R(t) \at{T_{\alpha(t)} M} &= A(t), \\
			R(t) \at{N_{\alpha(t)} M} &= C(t), \\
			s(t) &= \widehat{\alpha}(t) - R(t) \alpha(t),
		\end{split}
	\end{equation}
	is an extrinsic rolling in the sense of Definition~\ref{definition:extrinsic_rolling_2}.
	
	Conversely, given an extrinsic rolling $\big (\alpha(t), (R(t), s(t)) \big)$
	in the sense of Definition~\ref{definition:extrinsic_rolling_2},
	then $\big (\alpha(t), \widehat{\alpha}(t), A(t), C(t) \big)$
	defines an extrinsic rolling in the sense of Definition~\ref{definition:rolling_extrinsic_1},
	where 
	\begin{equation}
		\label{equation:proposition_relation_extrinsic_rolling_curve_quad-tuple}
		\begin{split}
			A(t) &= R(t) \at{T_{\alpha(t)} M}, \\
			C(t) &= R(t) \at{N_{\alpha(t)} M}, \\
			\widehat{\alpha}(t) &= s(t) + R(t) \alpha(t) .
		\end{split}
	\end{equation}
	\begin{proof}
		Since this proposition follows analogously to~\cite[Prop. 3]{JML-2023},
		we only sketch the proof.
		Let $\big (\alpha(t), \widehat{\alpha}(t), A(t), C(t) \big)$ be an extrinsic rolling
		in the sense of Definition~\ref{definition:extrinsic_rolling_2}
		and define 
		$I \ni t \mapsto \big (\alpha(t), (R(t), s(t)) \big) \in M \times E(V)$
		by~\eqref{equation:proposition_relation_extrinsic_rolling_curve_in_M_E_V}.
		We obtain
		\begin{equation}
			\label{equation:proposition_relation_extrinsic_rolling_rolling_curves}
			\begin{split}
				E(t) \alpha(t)
				&=
				R(t) \alpha(t) + s(t) \\
				&=
				R(t) \alpha(t) + \big( \widehat{\alpha}(t) - R(t) \alpha(t) \big) \\
				&=
				\widehat{\alpha}(t) \in \widehat{M}
			\end{split}
		\end{equation}
		showing Definition~\ref{definition:extrinsic_rolling_2},
		Claim~\ref{item:def_extrinsic_2_claim_1}.
		Let $\gamma \colon (- \epsilon, \epsilon) \to M$ be a curve
		with $\gamma(0) = \alpha(t)$ and $\dot{\gamma}(0) = Z \in V$.
		Then
		\begin{equation}
			\label{equation:proposition_relation_extrinsic_rolling_tangent_map_E}
			d_{\alpha(t)} E(t) Z = \tfrac{\D}{\D \tau} \big (R(t) \gamma(\tau) + s(t) \big) \at{\tau = 0} = R(t) Z
		\end{equation}
		holds.
		Using~\eqref{equation:proposition_relation_extrinsic_rolling_tangent_map_E}
		it is straightforward to verify that
		Definition~\eqref{definition:extrinsic_rolling_2} and 
		Claims~\ref{item:def_extrinsic_2_claim_2}-\ref{item:def_extrinsic_2_claim_5}
		are fulfilled.
		
		Conversely, assume that $I \ni t \mapsto M \times E(V)$ is a rolling
		in the sense of Definition~\ref{definition:extrinsic_rolling_2}.
		We now show that the quadruple
		$\big (\alpha(t), \widehat{\alpha}(t), A(t), C(t) \big)$, given 
		by~\eqref{equation:proposition_relation_extrinsic_rolling_curve_quad-tuple},
		is an extrinsic rolling in the sense of
		Definition~\ref{definition:rolling_extrinsic_1}.
		To this end, we note that
		$\widehat{\alpha}(t) 
		= 
		s(t) + R(t) \alpha(t)
		=
		E(t) \alpha(t)$
		holds by
		Definition~\ref{definition:extrinsic_rolling_2},
		Claim~\ref{item:def_extrinsic_2_claim_1}.
		Hence, by Definition~\ref{definition:extrinsic_rolling_2},
		Claim~\ref{item:def_extrinsic_2_claim_2}, the map
		\begin{equation}
			A(t)
			=
			R(t)\at{T_{\alpha(t)} M} 
			= \big( d_{\alpha(t)} E(t) \big) \at{T_{\alpha(t)} M} 
			\colon T_{\alpha(t)} M \to T_{\widehat{\alpha}(t)} M
		\end{equation}
		is indeed a well-defined isometry.
		Obviously, this implies that
		$C(t) = R(t)\at{N_{\alpha(t) }M} 
		= \big( d_{\alpha(t)} E(t) \big) \at{N_{\alpha(t)} M} 
		\colon N_{\alpha(t)} M \to N_{\widehat{\alpha}(t)} M$
		is a well-defined isometry, as well.
		Using Definition~\ref{definition:extrinsic_rolling_2}, Claim~\ref{item:def_extrinsic_2_claim_3}-\ref{item:def_extrinsic_2_claim_5},
		it is straightforward to show that $\big (\alpha(t),  \widehat{\alpha}(t), A(t), C(t) \big)$ 
		is indeed a rolling in the sense of Definition~\ref{definition:rolling_extrinsic_1}.
	\end{proof}
\end{proposition}

Below, in Section~\ref{sec:rolling_stiefel_manifolds},
we use Proposition~\ref{proposition:relation_extrinsic_rolling}
to relate the rolling of Stiefel manifolds constructed in this paper
to rolling maps of Stiefel manifolds known from the literature.

%-----------------
%
%In Section~\ref{sec:rolling_stiefel_manifolds} below, we clarify the relation of
%Definition~\ref{definition:rolling_extrinsic_1} and
%Definition~\ref{definition:extrinsic_rolling_2} on
%the example of the Stiefel manifold equipped with the Euclidean metric.
%
%
%-----------------

\section{Rolling Normal Naturally Reductive Homogeneous Spaces Intrinsically}
\label{sec:rolling_normal_naturally_reductive_homogeneous_spaces}

We first formulate an Ansatz for the rolling of normal naturally reductive homogeneous spaces, which is a generalization
of the rolling of pseudo-Riemannian symmetric spaces. It turns out,
however, that such an assumption does not work in general.

\subsection{No-Go Lemma}
\label{subsec:no_go_theorem}

Assume that $G / H$ is a pseudo-Riemannian symmetric space. 
Then, by~\cite[Sec. 4.2]{JML-2023},
a rolling of $\liealg{p}$ over $G / H$ along a given rolling curve
can be viewed as a triple
$\big (\alpha(t), \widehat{\alpha}(t), A(t) \big)$,
where 
\begin{equation}
  \begin{split}
    A(t) \colon T_{\alpha(t)} \liealg{p} \cong \liealg{p} &\to
    T_{\widehat{\alpha}(t)} (G / H) ,\\A(t) &= d_{q(t)} \pi \circ d_{e}
    L_{q(t)} .
  \end{split}
\end{equation}
Here $q \colon I \to G$ is defined by the initial value problem
\begin{equation}
	\dot{q}(t) = d_e L_{q(t)} \dot{\alpha}(t), \quad q(0) = e,
\end{equation}
whose solution is
the horizontal lift of the development
curve $\widehat{\alpha}(t) = \pi(q(t))$
through $q(0) = e$.

{ Note that  in~\cite{JML-2023},  $G / H$ is always rolled over $\liealg{p}$, 
while in our work we consider $\liealg{p}$ rolling over $G / H$. This choice 
is more convenient for us, since there is no need to invert $q(t)$, as in~\cite[Eq. 26]{JML-2023}.}

Motivated by this rather simple form of the intrinsic rolling
for symmetric spaces,
we make the following Ansatz for the
rolling of $\liealg{p}$ over $G  /H$, { where $q(t)$ will
be replaced by another lift of $\hat\alpha$, $r(t):=q(t)s(t)$,
$s(t)$ being a correction term, still to be specified, see below.\\}

\noindent\fbox{%
    \parbox{\textwidth}{%
 \begin{quotation}

{\bf Ansatz:}

Given a rolling curve
$\alpha \colon I \to \liealg{p}$, let $u \colon I \ni t \mapsto u(t) = 
\dot{\alpha }(t) \in \liealg{p}$, and define
the development curve 
$\widehat{\alpha} \colon I \to G / H$ by $\widehat{\alpha} (t)= \pi\big(q(t)\big)$,
with $q \colon I \to G$ being the horizontal curve
defined by
the initial value problem
\begin{equation}
	\label{equation:horizontal_lift_rolling_definition}
	\dot{q}(t) = d_e L_{q(t)} \big(\Ad_{s(t)} (u(t)) \big), 
	\quad
q(0) = e .
\end{equation}
Here $s \colon I \to H$ is a smooth curve that still needs to be specified.
The definition of $q$ in~\eqref{equation:horizontal_lift_rolling_definition}
is chosen such that the no-slip condition is satisfied, as will become clear in the computation \eqref{eq:22} below. 
As a candidate for the isometry
$A(t) \colon T_{\alpha (t)} \liealg{p} \cong \liealg{p} \to T_{\widehat{\alpha}(t)} (G / H)$,
we define
\begin{equation}
	\label{equation:ansatz_isometry}
	A(t)(Z) = \big( d_{r(t)} \pi \circ d_e L_{r(t)} \big) (Z),
	\quad Z \in T_{\alpha (t)} \liealg{p} \cong \liealg{p}
\end{equation}
where $r \colon I \ni t \mapsto q(t) s(t) \in G$
for some $s \colon I \to H$.
\end{quotation}}}
\begin{remark}
	If $G / H$ is a symmetric space,
	this yields a rolling of $\liealg{p}$ over $G / H$ for $s(t) = e$, see~\cite{JML-2023}.
\end{remark}

The more general situation, where $G / H$ is a naturally reductive
homogeneous space is considered in the following.
Our Ansatz satisfies the no-slip condition due to 
\begin{equation}\label{eq:22}
	\begin{split}
		A(t) \dot{\alpha}(t)
		&=
		d_{r(t)} \pi \circ d_e L_{r(t)} \, u(t) \\
		&=
		d_e (\pi \circ L_{q(t)} \circ  L_{s(t)}) \, u(t) \\
		&=
		d_e (\tau_{q(t)} \circ \pi \circ L_{s(t)} ) \, u(t) \\
		&= 
		d_e ( \tau_{q(t)} \circ \tau_{s(t)} \circ \pi ) \, u(t) \\
		&=
		d_{\pi(q(t))} \tau_{q(t)} \circ d_e \pi \circ \Ad_{s(t)}\, u(t) \\
		&=
		d_{q(t)} \pi \circ d_e L_{q(t)} \Ad_{s(t)}\, u(t) \\
		&= 
		d_{q(t)} \pi \,  \dot{q}(t) \\
		&= \dot{\widehat{\alpha}}(t)	,	
	\end{split} 
\end{equation}
where
$\tau \colon G \times G  / H \ni (g, g^\prime H) \mapsto (g g^\prime) H \in G/H$ denotes the $G$-action on $G / H$ from the left which
fulfills $\tau_g \circ \pi = \pi \circ L_g$, for $g \in G$.
Moreover, we exploited that the isotropy representation of $G / H$ and
the representation $\Ad \colon H \to GL(\liealg{p})$ are equivalent,
to be more precise,
$d_{\pi(e)} \tau_h \circ d_e \pi = d_e \pi \circ \Ad_h$, for $h \in H$,
see e.g.~\cite[Sec. 23.4, page 692]{gallier.quaintance:2020}.

Next we try to specify the curve $s \colon I \to H$ by
imposing the no-twist condition.
To this end, let $Z \colon I \ni t \mapsto \big(\alpha(t), Z_2(t)\big)  \in \liealg{p} \times \liealg{p} \cong T \liealg{p}$ be
a parallel vector field along $\alpha$. 
By identifying $Z$ with its second component $Z_2$,
$Z$ can be expressed by
$Z(t) = z$ for some $z \in \liealg{p}$.
We need to determine $s \colon I \to H$
such
that the vector field
$t \mapsto A(t) Z(t) = \big( d_{r(t)} \pi \circ d_e L_{r(t)} \big) z $
along $\widehat{\alpha}$ is parallel. Note that using~\eqref{equation:horizontal_lift_rolling_definition}
the curve $x(t) = \big(d_e L_{q(t)} \big)^{-1} \dot{q}(t)$ from 
Corollary~\ref{corollary:parallel_vf} corresponds
to $x(t) = \Ad_{s(t)}(u(t))$. Moreover, also due to
\begin{equation}\label{eq:1}
	\big(d_e L_{r(t)}\big)^{-1} \circ
	\big(d_{r(t)} \pi \at{\mathcal{H}_{r(t)}} \big)^{-1} A(t) (z) = z = const,
	\quad t \in I,        
\end{equation}
the condition $A(t) Z(t)$ being parallel tells us that 
\begin{equation}\label{eq:2}
	\begin{split}
		0 &=  - \tfrac{1}{2} \pr_{\liealg{p}} \big( \big[\Ad_{s(t)^{-1}}(\Ad_{s(t)}(u(t)), z \big]  \big) \\
		&= 
		- \tfrac{1}{2} \pr_{\liealg{p}} \big( [u(t), z] \big) \\
		&=
		- \tfrac{1}{2} \pr_{\liealg{p}} \big( [\dot{\alpha}(t), z] \big).
	\end{split}
\end{equation}
Assuming that for a given $0 \neq \dot\alpha(t) \in \liealg{p}$
there is a $z \in \liealg{p}$ such that
$0 \neq [\dot\alpha(t), z] \in \liealg{p}$ holds, 
\eqref{eq:2} cannot be satisfied independent of the choice of
$s \colon I \to H$.
We summarize the above discussion in the following lemma.
\begin{lemma}
	\label{lemma:rolling_in_general_not_of_simple_form} {\bf (No-Go)}
	Let $\alpha \colon I \to \liealg{p}$ be a curve
	so that
	$0 \neq \pr_{\liealg{p}}\big([\dot\alpha(t), z] \big)$
	holds for
	some $z \in \liealg{p}$ and some $t \in I$.
	%such
	%that there exists a $t \in I$ with
	%$0 \neq \pr_{\liealg{p}}\big([\dot\alpha(t), z] \big)$ for
	%some $z  \in \liealg{p}$.
	Then $\big(\alpha(t), \widehat{\alpha}(t), A(t)\big)$, as
	defined in the Ansatz at the beginning of this section, does not define a rolling of
	$\liealg{p}$ over $G / H$ no matter how $s \colon I \to H$ is chosen.
	To be more precise, the no-twist condition will never be fulfilled.
\end{lemma}

\subsection{Example: Stiefel Manifolds}
\label{subsec:no-go-stiefel}

We now specialize the above discussion
to the Stiefel manifold $\Stiefel{n}{k}$ (for the definition and more details see Section~\ref{subsec:stiefel_manifolds_alpha_metrics}), equipped with the $\alpha$-metrics
introduced in~\cite{HML-2021}. 
These metrics will be recalled in
Subsection~\ref{subsec:stiefel_manifolds_alpha_metrics}, below.
However, we think that it is convenient to apply
Lemma~\ref{lemma:rolling_in_general_not_of_simple_form}
to a non-trivial example here.
According to~\cite[Eq. (37)]{HML-2021}, for 
$E = \begin{bsmallmatrix}
	I_k \\
	0
\end{bsmallmatrix}$
and $\alpha \neq - 1$,
the projection
$\pr_{\liealg{p}} \colon \liealg{so}(n) \times \liealg{so}(k) \to \liealg{p}$
is given by 
\begin{equation}
	\label{equation:example_nogo_stiefel_projection}
	\pr_{\liealg{p}} \Big( \begin{bsmallmatrix}
		A && - B^{\top} \\
		B  && C
	\end{bsmallmatrix}, \Psi\Big)
	=
	\Big(
	\begin{bsmallmatrix}
		\tfrac{A -   \Psi}{\alpha + 1} && - B^{\top} \\
		B && 0		
	\end{bsmallmatrix}, 
	\tfrac{\alpha(\Psi -A)}{\alpha + 1}
	\Big) .
\end{equation}
We first assume $1 \leq k \leq n - 1$.
Setting $\Psi = A$,  we get
elements of the form
$\Big( \begin{bsmallmatrix}
	0 && - B^{\top} \\
	B && 0
\end{bsmallmatrix}, 0 \Big)
\in \liealg{p}$, where $B \in \matR{(n - k)}{k}$.
Using~\eqref{equation:example_nogo_stiefel_projection},
we can write
\begin{equation}
	\label{equation:example_nogo_stiefel_projection_computation_projection_lie_bracket}
	\begin{split}
		&\pr_{\liealg{p}}\Big[ \Big(\begin{bsmallmatrix}
			0 && - B_1^{\top} \\
			B_1 && 0
		\end{bsmallmatrix}, 0 \Big]\Big),
		\Big( \begin{bsmallmatrix}
			0 && - B_2^{\top} \\
			B_2 && 0
		\end{bsmallmatrix}, 0 \Big)\Big] \\
		&=
		\pr_{\liealg{p}}\Big(\begin{bsmallmatrix}
			- B_1^{\top} B_2 + B_2^{\top} B_1 && 0 \\
			0	&& - B_1 B_2^{\top} + B_2 B_1^{\top}
		\end{bsmallmatrix}, 0\Big) \\
		&=
		\Big(
		\begin{bsmallmatrix}
			\tfrac{	- B_1^{\top} B_2 + B_2^{\top} B_1}{\alpha + 1} 	&& 0 \\
			0 && 0
		\end{bsmallmatrix},
		\tfrac{\alpha}{\alpha + 1} (B_1^{\top} B_2 - B_2^{\top} B_1)
		\Big) .
	\end{split}
\end{equation}

Obviously, for $k = 1$, i.e. $B_1, B_2 \in \matR{(n - 1)}{1}$,
one has $B_2^{\top} B_1 = B_1^{\top} B_2$ implying
that~\eqref{equation:example_nogo_stiefel_projection_computation_projection_lie_bracket} is vanishing for $k = 1$.
Thus for $\Stiefel{n}{1}\cong S^{n-1}$, the
		Ansatz actually yields a rolling.

Next assume $k > 1$. Then, there are
$B_1, B_2 \in \matR{(n - k)}{k}$ such that
$B_2^{\top} B_1 - B_1^{\top} B_2 \neq 0$ holds.
Indeed, choosing $B_2 = E_{1 2}$ given by
$(E_{12})_{ij} = \delta_{1i} \delta_{2j}$,
where $\delta_{1i}$ and $\delta_{2j}$ are Kronecker deltas, and
$B_1 \in \matR{(n - k)}{k}$ with $(B_1)_{1 2} \neq 0$,
we obtain
\begin{equation}
	\begin{split}
		\big(B_2^{\top} B_1 - B_1^{\top} B_2\big)_{2 2}
		&=
		\sum_{\ell = 1}^{n - k}
		\big((E_{12})_{k 2} (B_1)_{k 2} - (B_1)_{k 2} (E_{12})_{k 2}  \big) \\
		&=
		\sum_{\ell = 1}^{n - k}
		\big( \delta_{1 k} \delta_{2 2} (B_1)_{k 2} - (B_1)_{k 2} \delta_{2 k} \delta_{1 2}  \big) \\
		&= 
		(B_1)_{12} 
		\neq 0 .
	\end{split}
\end{equation}
Thus~\eqref{equation:example_nogo_stiefel_projection_computation_projection_lie_bracket}
does not vanish identically for $1 < k < n$.
It remains to consider the case $k = n$. 
This yields $\Stiefel{n}{n} \cong (O(n) \times O(n)) / O(n)$ and
for $(A, \Psi) \in \liealg{so}(n) \times \liealg{so}(n)$
the projection~\eqref{equation:example_nogo_stiefel_projection} reduces
to
\begin{equation}
	\pr_{\liealg{p}} \big( A, \Psi\big)
	=
	\Big(
	\tfrac{A -   \Psi}{\alpha + 1} ,
	\tfrac{\alpha(\Psi -A)}{\alpha + 1}
	\Big) .
\end{equation}
Parameterize $\liealg{p}$ by 
\begin{equation}
	\liealg{p}
	=
	\Big\{ (\tfrac{A}{\alpha + 1}, - \tfrac{\alpha A}{\alpha + 1} ) \mid A \in \liealg{so}(n)\Big\} .
\end{equation}
Consequently, we obtain for $A_1, A_2 \in \liealg{so}(k)$
\begin{equation}
	\begin{split}
		\pr_{\liealg{p}}\Big[ (\tfrac{A_1}{\alpha + 1}, 
		- \tfrac{\alpha A_1}{\alpha + 1} ), (\tfrac{A_2}{\alpha + 1}, - \tfrac{\alpha A_2}{\alpha + 1} ) \Big]
		&=
		\pr_{\liealg{p}} \Big( \tfrac{[A_1, A_2 ]}{(\alpha + 1)^2}, 
		\tfrac{\alpha^2 [A_1, A_2]}{(\alpha + 1)^2} \Big) \\
		&=
		\Big( \tfrac{[A_1, A_2]- \alpha^2 [A_1, A_2]}{(\alpha +1)^3},
		\tfrac{\alpha ( \alpha^2 [A_1, A_2] -  [A_1, A_2]}{
			(\alpha + 1)^3} \Big) .
	\end{split}
\end{equation}
Clearly, this equation vanishes for $k = n = 1$ and all
$\alpha \in \field{R} \setminus \{-1\}$.
Moreover, it vanishes for $k = n > 1$ and  all $A_1, A_2 \in \liealg{so}(n)$
iff $\alpha = 1$ holds. (Note that $\alpha = -1$ is excluded by the definition of the $\alpha$-metrics in~\cite[Def. 3.1]{HML-2021})
We summarize this computations in the next corollary.

\begin{corollary}
	Let $1 < k < n$ and let $\alpha \in \field{R} \setminus \{- 1, 0\}$.
	Then the Ansatz from Subsection~\ref{subsec:no_go_theorem}
	does not yield an intrinsic rolling, with respect to any $\alpha$-metric,
	of a tangent space of the Stiefel manifold over the
	Stiefel manifold $\Stiefel{n}{k}$.
	However, for the case $k = n > 1$ the Ansatz yields only a rolling for
	$\alpha = 1$.
\end{corollary}

%\newpage

\subsection{Kinematic Equations for Intrinsic Rolling}
\label{subsec:kinematic_equation_intrinsic_rolling}

We aim to find the triple $\big(\alpha(t),\widehat\alpha (t), A(t)\big)$
satisfying Definition~\ref{definition:intrinsic_rolling}
for a rolling of $\liealg{p}$ over the normal naturally reductive 
homogeneous space $G / H$.

More precisely, our goal is to find a system of ODEs,
the so-called kinematic equations, which, for a prescribed
rolling curve $\alpha \colon I \to \liealg{p}$, determines the development
curve $\widehat{\alpha} \colon I \to G / H $
as well as the curve of isometries
$A(t) \colon T_{\alpha(t)} \liealg{p} \cong \liealg{p} 
\to T_{\widehat{\alpha}(t)} (G / H)$.

The new  terminology in the next definition is motivated by
the theory of control, since the kinematic equations can be
written as a control system whose control function is
precisely  $\dot{\alpha}(t)$.

\begin{definition}
	Given a rolling curve $\alpha \colon I \to \liealg{p}$, we call the curve
	$u \colon I \to \liealg{p}$,
	defined by $u(t) = \dot{\alpha}(t)$, the associated control curve.
\end{definition}
Note that a prescribed control curve $u \colon I \to \liealg{p}$ determines uniquely the rolling curve $\alpha \colon I \to \liealg{p}$ up to the initial condition $\alpha(0) = \alpha_0 \in \liealg{p}$.

In order to derive the kinematic equations we start with the following remark.
\begin{remark}
	\label{remark:intrinsic_rolling_stiefel_remark_isometries}
	Let $V$ and $W$ be 
	finite dimensional pseudo-Euclidean vector spaces whose scalar
	products have the same signature and let $\phi \colon V \to W$
	be an isometry.
	Then, the set of isometries between $V$ and $W$ is given by
	$\{ \phi \circ S \colon V \to W \mid S \in O(V) \}$.
	Indeed, for $S \in O(V)$, $\phi \circ S$ is 	a composition
        of isometries, so it is an isometry, as well. Conversely, given an isometry
	$\psi \colon V \to W$, define the isometry
	$S=\phi^{-1} \circ \psi \colon V \to V$ which is an element of $O(V)$, and clearly $\psi=\phi \circ S$.
\end{remark}
In view of Remark~\ref{remark:intrinsic_rolling_stiefel_remark_isometries},
a possible candidate for the curve of isometries $A(t) \colon T_{\alpha(t)} \liealg{p} \cong \liealg{p} 
\to T_{\widehat{\alpha}(t)} (G / H) $ that is required
for an intrinsic rolling is of the form 
\begin{equation}
	A(t)  = \big( d_{q(t)} \pi \big) \circ \big( d_e L_{q(t)} \big) \circ S(t),
\end{equation}
where $q \colon I \to G$ is the horizontal lift
of the devolopement curve $\widehat{\alpha} \colon I \to G / H$
through $q(0) = e$ and $S \colon I \to O(\liealg{p})$ is a curve
in the orthogonal group of $\liealg{p}$ through $S(0) = \id_{\liealg{p}}$.

In the next theorem, we reproduce from~\cite{schlarb:2023} the kinematic equations for the rolling of $\liealg{p}$ over $G / H$. This statement holds for general  
 normal 
naturally reductive homogeneous spaces, and the proof is provided to keep the paper as  
 self-contained as possible.

\begin{theorem}
	\label{theorem:intrinsic_rolling}
	Let $G / H$ be a normal naturally reductive homogeneous space,  $\alpha \colon I \to \liealg{p}$  a given curve,
	and  $u \colon I \to \liealg{p}$ defined by
	$u(t) = \dot{\alpha}(t)$
	 the associated control curve.
	Moreover, let $S \colon I \to O(\liealg{p})$ and 
	$q \colon I \to G$ be determined by the initial value problem
	\begin{equation}
		\label{kinematicsFL}
		\begin{split}
			\dot{S}(t) 
			&=
			- \tfrac{1}{2} \pr_{\liealg{p}} \circ \ad_{S(t) u(t)} \circ S(t) ,
			\quad S(0) = \id_{\liealg{p}}, \\
			\dot{q}(t) 
			&= 
			\big( (d_e L_{q(t)}) \circ S(t)\big) u(t),
			\quad
			q(0) = e.
		\end{split} 
	\end{equation}
	Then, the triple $\big(\alpha(t), \widehat{\alpha}(t), A(t)\big)$, where
	\begin{equation}
		\widehat{\alpha} \colon I \to G /  H,
		\quad
		t \mapsto \widehat{\alpha}(t) = (\pi \circ q)(t)
	\end{equation}
	and 
	 {\begin{equation}
		t \mapsto A(t)=(d_{q(t)} \pi) \circ (d_e L_{q(t)}) \circ S(t) \colon T_{\alpha(t)} \liealg{p} \cong \liealg{p}
		\to T_{\widehat{\alpha}(t)} (G / H),
	\end{equation}
	is an intrinsic rolling of $\liealg{p}$ over $G/H$.}
	\begin{proof}
		We show that $(\alpha(t), \widehat{\alpha}(t), A(t))$ satisfies
		the conditions of Definition~\ref{definition:intrinsic_rolling}.
		The solution $S$
		of the first equation in~\eqref{kinematicsFL} is
		indeed a curve in $O(\liealg{p})$
		since
		$- \tfrac{1}{2} \pr_{\liealg{p}} \circ \ad_{S u} \colon \liealg{p} \to \liealg{p}$
		is skew adjoint for all $S \in O(\liealg{p})$ and $u \in \liealg{p}$
		with respect to the scalar-product on $\liealg{p}$ defined by means
		of the bi-invariant metric on $G$.
		In fact, by exploiting that $G / H$ is naturally reductive, using Definition~\ref{definition:naturally_reductive_space}, we obtain
		for $X, Y \in \liealg{p}$.
		\begin{equation}
			\begin{split}
				\big\langle - \tfrac{1}{2} \pr_{\liealg{p}} \circ \ad_{S u} (X), Y  \big\rangle
				&=
				\big\langle - \tfrac{1}{2}\pr_{\liealg{p}} \big([S u, X] \big), Y \big\rangle \\
				&= 
				\big\langle X , \tfrac{1}{2} \pr_{\liealg{p}} \circ \ad_{S u}(Y) \big\rangle,
			\end{split}
		\end{equation}
		showing that
		$- \tfrac{1}{2} \pr_{\liealg{p}} \circ \ad_{S u} \in \liealg{so}(\liealg{p})$.
		Thus $S(t) \in O(\liealg{p})$  since it is the integral curve of the
		time-variant vector field  $- \tfrac{1}{2} \pr_{\liealg{p}} \circ \ad_{S u(t)} \circ S$ on $O(\liealg{p})$.
		
		Next, we set $\widehat{\alpha}(t) = (\pi \circ q)(t)$.
		Obviously, the ODE for $q$ in~\eqref{kinematicsFL}
		implies that $q \colon I \to G$ is the horizontal lift of $\widehat{\alpha}$ through $q(0) = e$.
		Moreover, the map
		$A(t) \colon T_{\alpha(t)} \liealg{p} \cong \liealg{p}
		\to T_{\widehat{\alpha}(t)} (G / H)$ 
		is well-defined and an isometry since it is a composition of isometries.
		
		We now show the no-slip condition.
		Indeed, by the chain-rule
		\begin{equation}
			\begin{split}
				\dot{\widehat{\alpha}}(t)
				&= 
				\tfrac{\D}{\D t}
				(\pi \circ q(t)) \\
				&=
				(d_{q(t)} \pi ) \dot{q}(t) \\
				&= 
				\big(d_{q(t)} \pi \big) \big( d_e L_{q(t)} \circ S(t) \big) u(t) \\
				&=
				A(t) \dot{\alpha}(t) .
			\end{split}
		\end{equation}
		
		It remains to show the no-twist condition.
		Let $Z \colon I \to \liealg{p}$ be a parallel vector field along $\alpha \colon I \to \liealg{p}$, i.e. $Z$ can be viewed as a constant function
		$Z(t) = Z_0$ for all $t \in I$ and some $Z_0 \in \liealg{p}$.
		We prove that the vector field $\widehat{Z}(t) = A(t) Z_0$ is parallel along the curve $\widehat{\alpha}$, by exploiting the result in Corollary~\ref{corollary:parallel_vf_horizontal}.
		The curve $z \colon I \to \liealg{p}$ defined by
		\begin{equation}
			z(t) 
			=
			(d_e L_{q(t)} )^{-1}  \circ (d_{q(t)} \pi )^{-1} A(t) Z_0 
			=
			S(t) Z_0
		\end{equation}
		fulfills
		\begin{equation}
			\begin{split}
				\dot{z}(t)
				&= 
				\dot{S}(t) Z_0 \\
				&=
				- \tfrac{1}{2} \circ \pr_{\liealg{p}} \circ \ad_{S(t) u(t)} \circ S(t)(Z_0) \\
				&= 
				- \tfrac{1}{2} \big[S(t) u(t), S(t)( Z_0) \big]\at[\Big]{\liealg{p}} \\
				&= 
				- \tfrac{1}{2} \big[S(t) u(t), z(t) \big]\at[\Big]{\liealg{p}} .
			\end{split}
		\end{equation}
		Thus $Z(t) = A(t)Z_0$ is parallel along $\widehat{\alpha}(t) = (\pi \circ q)(t)$ 
		by Corollary~\ref{corollary:parallel_vf_horizontal},
		due to $\big(d_e L_{q(t)} \big)^{-1} \dot{q}(t) = S(t) u(t)$.
		
		Conversely, assume that $A(t) Z(t)$ is parallel along
		$\widehat{\alpha}$ for some vector field $Z(t)$ along $\alpha$.
		We define the parallel frame $A_i(t) = A(t) A_i$,
		where $\{A_1,\cdots , A_k\}$  
		forms a basis of $\liealg{p}$, and expand $A(t) Z(t)$ in this basis to obtain
		$A(t) Z(t) = \sum_{i = 1}^k z_i A_i(t)$,
		where the coefficients $z_i \in \field{R}$ are constant, since $A(t) Z(t)$ is assumed to be parallel, see~\cite[Chap. 4, p.109]{Lee:2018aa}.
				By the linearity of $A(t)$, we obtain
		\begin{equation}
			A(t) Z(t)
			=
			\sum_{i = 1}^k z_i A_i(t)
			=
			 A(t) \Bigg( \sum_{i = 1}^k z_i A_i \Bigg)
			= A(t) Z_0,
		\end{equation}
		for $Z_0 = \sum_{i = 1}^k z_i A_i \in \liealg{p}$,
		i.e. $Z(t) = Z_0$
		is constant.
		Thus $Z(t)$ is a parallel vector field along $\alpha$
		as desired.
	\end{proof}
\end{theorem}

\begin{remark}
	\label{remark:completeness_intrinsic_rolling_S}
	It is not clear whether the curve $S \colon I \to O(\liealg{p})$
	from Theorem~\ref{theorem:intrinsic_rolling}
	is defined on the same interval $I$ as the control
	curve $u \colon I \to \liealg{p}$
	due to the nonlinearity
	of~\eqref{kinematicsFL}.
	We cannot rule out that $S$ is defined only on a proper
	subintervall $I^{\prime} \subsetneq I$ with $0 \in I^{\prime}$.
	By abuse of notation, we write $S \colon I \to O(\liealg{p})$, nevertheless, 
	even if $S$ was defined on a proper subinterval.
	However, we are not aware of an example.
\end{remark}
If $G / H$ is a Riemannian normal naturally reductive space,
i.e. if the metric is positive definite,
and the control
defined on $\field{R}$ is bounded, following~\cite{schlarb:2023},
we can prove that $S$ is defined on the
whole interval $\field{R}$.
This is the next lemma.

\begin{lemma}
	\label{lemma:completeness_intrinsic_rolling_S}
	Let $u \colon \field{R} \to \liealg{p}$ be bounded
	and let $G / H$ be a Riemannian normal naturally reductive
	homogeneous space.
	Then the vector field given by
	\begin{equation}
		X(t, S) = \big(1, - \tfrac{1}{2} \pr_{\liealg{p}} \circ \ad_{S(t)u(t)} \circ S(t)\big)
	\end{equation}
	on $\field{R} \times O(\liealg{p})$ is complete.
	\begin{proof}
			We will show that this vector field is bounded in a complete Riemannian metric on $\field{R} \times O(\liealg{p})$.
		Completeness then 
		follows by~\cite[Prop. 23.9]{michor:2008}.
		To this end, we view $O(\liealg{p})$ as subset of $\End(\liealg{p})$.
		Since $G / H$ is Riemannian,
		the corresponding scalar product on
		$\liealg{p}$ denoted by $\langle \cdot, \cdot \rangle$ is positive definite,
		i.e. an inner product.
		The norm on $\liealg{p}$ induced by this inner product is
		denoted $\Vert \cdot \Vert$.
		We denote an extension of $\langle \cdot , \cdot \rangle$ to an inner product on $\liealg{g}$ by $\langle \cdot , \cdot \rangle$, too.
		The corresponding norm is denoted by $\Vert \cdot \Vert$, as well.
		We now endow $\End(\liealg{p})$ with the Frobenius scalar product
		given by $\langle S, T \rangle_F = \tr(S^{\top} T)$,
		where $S^{\top}$ is the adjoint of $S$ with respect
		to $\langle \cdot, \cdot \rangle$.
		Then $\langle \cdot, \cdot \rangle_F$ induces
		a bi-invariant and hence a complete metric on $O(\liealg{p})$.
		Moreover, the norm $\Vert \cdot \Vert_F$ defined by the Frobenius
		scalar product is equivalent to the operator norm $\Vert \cdot \Vert_2$.
		In particular, there is a $C > 0$ such that
		$\Vert S \Vert_F \leq C \Vert S \Vert_2$ holds
		for all $S \in \End({\liealg{p}})$.
		In addition, on the $\field{R}$-component,
		define the metric to be the Euclidean metric.
		In other words, the Riemannian metric on $\field{R} \times O(\liealg{p})$ is given by 
		\begin{equation}
			\langle (v, V), (w, W) \rangle_{(s, S)} 
			=
			v w + \tr(V^{\top} W) 
		\end{equation}
		for all $(s, S) \in \field{R} \times O(\liealg{p})$ and
		$(v, V), (w, W)  \in T_{(s, S)} (\field{R} \times O(\liealg{p}))$.
		Moreover, $\ad \colon \liealg{p} \times \liealg{p} \to \liealg{g}$ is
		bounded since $\liealg{p}$ is finite dimensional.
		Hence there exists a $C^{\prime} \geq 0$
		with $\Vert \ad(X, Y) \Vert \leq C^{\prime} \Vert X \Vert \Vert Y \Vert$.
		Consequently, for fixed $X \in \liealg{p}$,
		the operator norm of $\ad_X \colon \liealg{p} \to \liealg{g}$
		can be estimated by $\Vert \ad(X, \cdot) \Vert_2 \leq C^{\prime} \Vert X \Vert$.
		By this notation, we compute
		\begin{equation}
			\begin{split}
				\Vert X(t, S) \Vert^2_{\field{R} \times O(\liealg{p})}
				&=
				1 +  \Vert \tfrac{1}{2} \pr_{\liealg{p}} \circ \ad_{S u(t)} \circ S \Vert_F^2 \\
				&\leq
				1 + \tfrac{C^2}{4} \Vert \pr_{\liealg{p}} \circ \ad_{S u(t)} \circ S \Vert_2^2 \\
				&\leq
				1 + \tfrac{C^2}{4} \Vert \pr_{\liealg{p}} \Vert_2^2 \Vert \ad_{S u(t)} \Vert_2^2 \Vert S \Vert_2^2 \\
				&\leq
				1 + \tfrac{(C C^{\prime})^2}{4}   \Vert S \Vert_2^2 \Vert u(t) \Vert^2 \\
				&\leq
				1+ \tfrac{(C C^{\prime})^2}{4}  \Vert u \Vert_{\infty}^2 < \infty,
			\end{split}
		\end{equation}
		where $\Vert u \Vert_{\infty}$ denotes the supremum norm of $u$ and
		we exploited $\Vert S \Vert_2 = 1$ due to
		$S \in O(\liealg{p})$ and $\Vert \pr_{\liealg{p}} \Vert_2 \leq 1$
		showing that $X$ is bounded in a complete Riemannian metric.
	\end{proof}
\end{lemma}

\section{Rolling Stiefel Manifolds}
\label{sec:rolling_stiefel_manifolds}

A first attempt to generalize the rolling for pseudo-Riemannian
symmetric spaces, as discussed  in Section~\ref{sec:rolling_normal_naturally_reductive_homogeneous_spaces},
does not work for Stiefel manifolds by Subsection~\ref{subsec:no-go-stiefel}.
However,  rolling maps for Stiefel manifolds have already appeared
in~\cite{hueper.kleinsteuber.leite:2008}, and more recently also in~\cite[Sec. 5]{JML-2023}.

In this section, we reformulate the  most  recent results
in~\cite{schlarb:2023}, without using fiber-bundle techniques,
to describe the intrinsic rolling of Stiefel manifolds equipped with the
so-called $\alpha$-metrics defined in~\cite{HML-2021}.
 Although, up to now, we have used the Greek letter $\alpha$ for rolling curves, in the first part of this section we will use the same letter $\alpha$ for the real parameter that defines a family of metrics  on Stiefel manifolds. This will not create difficulties, since it will be clear from the context. 
In order to reach the above mentioned objective, we specialize Theorem~\ref{theorem:intrinsic_rolling}
to Stiefel manifolds.
Eventually, this rolling is extended to an extrinsic rolling for
the Euclidean metric.
Finally, we show that our findings coincide with the rolling results
from~\cite{hueper.kleinsteuber.leite:2008}.

\subsection{Stiefel Manifolds Equipped with $\alpha$-Metrics as Normal Naturally Reductive Homogeneous Spaces}
\label{subsec:stiefel_manifolds_alpha_metrics}

The Stiefel manifold $\Stiefel{n}{k}$
can be viewed as the embedded submanifold
\begin{equation}
	\Stiefel{n}{k} = \{X \in \matR{n}{k} \mid X^{\top} X  = I_k \}, \quad 1 \leq k \leq n 
\end{equation}
of $\matR{n}{k}$.
In the sequel, we recall the so-called $\alpha$-metrics on $\Stiefel{n}{k}$
introduced in~\cite{HML-2021} and show that $\Stiefel{n}{k}$ equipped with an $\alpha$-metric can be viewed as a normal naturally reductive homogeneous space.
The $\big(O(n) \times O(k)\big)$-left action
\begin{equation}
	\Phi \colon \big(O(n) \times O(k)\big) \times \matR{n}{k }\to \matR{n}{k},
	\quad
	\big((R, \theta), X\big) \to R X \theta^{\top}
\end{equation}
by linear isomorphisms restricts to a transitive action
\begin{equation}
		(O(n) \times O(k)) \times \Stiefel{n}{k }\to \Stiefel{n}{k},
	\quad
	((R, \theta), X) \to R X \theta^{\top}
\end{equation}
on $\Stiefel{n}{k}$, also denoted by $\Phi$, which coincides with the action from~\cite[Eq. (13)]{HML-2021}.
Next, let $X \in \Stiefel{n}{k}$ be fixed, and denote by
$H = \Stab(X) \subset O(n) \times O(k)$ the isotropy
subgroup of $X$ under the action $\Phi$.
Moreover, we write $G = O(n) \times O(k)$. 
Then, the Stiefel manifold $\Stiefel{n}{k}$ is diffeomorphic to
the homogeneous space $G / H$.
Moreover, the map
\begin{equation}
	\label{equation:embedding_stiefel_manifold_G_equivariant}
	\iota_X \colon G / H \ni 
	(R, \theta) \cdot H \mapsto R X \theta^{\top} 
	\in \Stiefel{n}{k} \subset \matR{n}{k}
\end{equation}
is a $G$-equivariant embedding, where $(R, \theta) \cdot H$ denotes
the coset in $G / H$ represented by $(R, \theta) \in G$.

In order to construct the $\alpha$-metrics, 
the map 
\begin{equation}
		\langle \cdot, \cdot \rangle^{\alpha}_{\liealg{so}(n) \times \liealg{so}(k)} 
		\colon 
		\liealg{so}(n) \times \liealg{so}(k)
		\to \field{R}
\end{equation}
is defined on $\liealg{so}(n) \times \liealg{so}(k)$, for $\alpha \in \field{R} \setminus \{0\}$,
by \begin{equation}
	\label{equation:alpha_metrics_lie_algebra_def}
	\big\langle (\Omega_1, \Psi_1), (\Omega_2, \Psi_2) \big\rangle^{\alpha}_{\liealg{so}(n) \times \liealg{so}(k)} 
	= 
	- \tr(\Omega_1 \Omega_2) - \tfrac{1}{\alpha} \tr(\Psi_1 \Psi_2),
\end{equation}
see~\cite[Eq. (21)]{HML-2021}.

Obviously, $\langle \cdot, \cdot \rangle^{\alpha}_{\liealg{so}(n) \times \liealg{so}(k)} $
yields a symmetric bilinear form 
on $\liealg{g} =  \liealg{so}(n) \times \liealg{so}(k)$
which is
$\Ad_G$-invariant.
Moreover, by~\cite[Prop. 2]{HML-2021}, the subspace $\liealg{h} \subset \liealg{g}$
being the Lie algebra of $H = \Stab(X)$ for $X \in \Stiefel{n}{k}$
is non-degenerated for all $\alpha \in \field{R} \setminus \{-1, 0\}$.

After this preparation, we are in the position to
reformulate~\cite[Def. 3.3]{HML-2021}.
\begin{definition}
	Let $\alpha \in \field{R} \setminus \{-1, 0\}$. 
	The $\alpha$-metric on $\Stiefel{n}{k} \cong G / H$ is defined
	as the $G = O(n) \times O(k)$-invariant metric on $G / H$ 
	that turns the canonical projection $\pi \colon G \to G / H$ into a pseudo-Riemannian submersion,
	where $G$ is equipped with the bi-invariant metric defined by 
	means of the scalar product
	from~\eqref{equation:alpha_metrics_lie_algebra_def}.
\end{definition}
This definition turns $G / H$ into a normal
naturally reductive homogeneous space.

\begin{lemma}
	\label{lemma:stiefel_naturally_reductive}
	Let $\alpha \in \field{R} \setminus \{-1, 0\}$. 
	Then $G / H \cong \Stiefel{n}{k}$ equipped with
	an $\alpha$-metric is a normal naturally reductive space.
	In particular, it is a naturally reductive homogeneous space.
	\begin{proof}
		Obviously, $\Stiefel{n}{k}  \cong G / H$
		is a normal naturally reductive homogeneous space
		since the metric on $G$ is bi-invariant and 
		$\liealg{h} \subset \liealg{g}$ is a non-degenerated subspace.
		Hence it is naturally reductive by
		Lemma~\ref{lemma:normal_naturally_reductive_implies_reductive}.
	\end{proof}
\end{lemma}
By requiring that $\iota_X \colon G / H \to \Stiefel{n}{k}$
from~\eqref{equation:embedding_stiefel_manifold_G_equivariant}
is an isometry, the $\alpha$-metric
on $\Stiefel{n}{k}$ for $\alpha \in \field{R} \setminus\{-1, 0\}$, viewed
as an embedded submanifold
of $\matR{n}{k}$, is given by
\begin{equation}
	\langle V, W \rangle_X^{(\alpha)} = 2 \tr(V^{\top} W)
	 + \tfrac{2 \alpha +1}{\alpha + 1} \tr(V^{\top} X X^{\top} W),
\end{equation}
where $X \in \Stiefel{n}{k}$ and $V, W \in T_X \Stiefel{n}{k}$
by~\cite[Cor. 2]{HML-2021}.
In addition, if $\Stiefel{n}{k}$ is equipped with an $\alpha$-metric,
and $O(n) \times O(k)$ is equipped with the corresponding bi-invariant metric
defined by the scalar product
from~\eqref{equation:alpha_metrics_lie_algebra_def}, the map
\begin{equation}
	\label{equation:alpha_metrics_pseudo_Riemannian_submersion}
	\Phi_X = \iota_X \circ \pi \colon O(n) \times O(k) \to \Stiefel{n}{k},
	\quad
	(R, \theta) \mapsto R X \theta^{\top}
\end{equation}
is a pseudo-Riemannian submersion, where $X \in \Stiefel{n}{k}$ is
arbitrary but fixed.

For considering the intrinsic rolling of $\Stiefel{n}{k} \cong G / H$,
we need a formula for the orthogonal projection
$\pr_{\liealg{p}} \colon \liealg{so}(n) \times \liealg{so}(k) \to \liealg{p}$
with respect to the metric defined in~\eqref{equation:alpha_metrics_lie_algebra_def}, where $\liealg{p} = \liealg{h}^{\perp}$,
 $\liealg{h}$ being the Lie algebra of $H = \Stab(X) \subset G$
for a fixed $X \in \Stiefel{n}{k}$.
This is the next lemma which is taken from~\cite[Lem. 3.2]{HML-2021}.

\begin{lemma}
	\label{lemma:stiefel_orthogonal_projection_on_p}
	Let $\alpha \in \field{R} \setminus \{-1, 0\}$.
	The orthogonal projection
	\begin{equation}
		\pr_{\liealg{p}}  
		\colon \liealg{so}(n) \times \liealg{so}(k)
		\to \liealg{p},
		\quad
		(\Omega, \eta) \mapsto (\Omega^{\perp_X}, \eta^{\perp_X})
	\end{equation}
	is given by
	\begin{equation}
		\begin{split}
			\Omega^{\perp_X} 
			&=
			X X^{\top} \Omega 
			+ \Omega X X^{\top} 
			- \tfrac{2 \alpha + 1}{\alpha +1} X X^{\top} \Omega X X^{\top} 
			- \tfrac{1}{\alpha + 1} X \eta X^{\top}
			\\
			\eta^{\perp_X} &= \tfrac{\alpha}{\alpha + 1} \big( \eta - X^{\top} \Omega X  \big) .
		\end{split}
	\end{equation}
	\begin{proof}
		This is just a reformulation of \cite[Lem. 3.2]{HML-2021}.
	\end{proof}
\end{lemma}
Since $\pi \colon G \to G / H$ is a pseudo-Riemannian submersion
whose horizontal bundle is defined
point-wise by
$\mathcal{H}_g = (d_{(I_n, I_k)} L_g)(\liealg{p}) \subset T_g G$
and $\iota_X \colon G / H \to \Stiefel{n}{k}$ is an isometry,
the map
\begin{equation}
	d_{(I_n, I_k)} (\iota_X \circ \pi)\at{\liealg{p}} \colon \liealg{p} \to T_X \Stiefel{n}{k}, 
	\quad
	(\Omega, \eta) \mapsto \Omega X - X \eta
\end{equation}
as well as its inverse are linear isometries.
For the discussion of rolling Stiefel manifolds,
we need an explicit formula for
\begin{equation}
	\big(d_{(I_n, I_k)} (\iota_X \circ \pi)\at{\liealg{p}}\big)^{-1} \colon T_X \Stiefel{n}{k} \to \liealg{p} .
\end{equation}
Such a formula is given in the following lemma
which is a trivial reformulation
of~\cite[Prop. 3]{HML-2021}.

\begin{lemma}
	\label{lemma:stiefel_identification_tangent_space_p}
	Let $\alpha \in \field{R} \setminus \{-1, 0\}$ and $X \in \Stiefel{n}{k}$.
	The map 
	\begin{equation}
		\big(d_{(I_n, I_k)} (\iota_X \circ \pi)\at{\liealg{p}}\big)^{-1} \colon T_X \Stiefel{n}{k} \to \liealg{p},
		\quad
		\quad
		V \mapsto \big(\Omega(V)^{\perp_X}, \eta(V)^{\perp_X} \big)
	\end{equation}
	is given by 
	\begin{equation}
		\begin{split}
			\Omega(V)^{\perp_X} 
			&=
			V X^{\top} 
			- X V^{\top} 
			+ \tfrac{2 \alpha + 1}{\alpha + 1} X V^{\top} X X^{\top}
			\\
			\eta(V)^{\perp_X}
			&= 
			-\tfrac{\alpha}{\alpha + 1} X^{\top} V .
		\end{split}
	\end{equation}
	\begin{proof}
		This is a consequence of \cite[Prop. 3]{HML-2021}.
	\end{proof}
\end{lemma}

Finally, we specialize the previous two lemmas for $\alpha = - \tfrac{1}{2}$.
For this choice, the $\alpha$-metric coincides with the Euclidean metric,
scaled by the factor $2$, see~\cite[Sec. 4.2]{HML-2021}.
Therefore this special case will be important for discussing the extrinsic rolling
of Stiefel manifolds equipped with the Euclidean metric.

\begin{corollary}
	\label{corollary:stiefel_orthogonal_projection_identification_euclidean_metric}
	Let $\alpha = - \tfrac{1}{2}$.
	Using the notation of Lemma~\ref{lemma:stiefel_identification_tangent_space_p},
	the following assertions are fulfilled:
	\begin{enumerate}
		\item
		The projection 
		$\pr_{\liealg{p}} 
		\colon \liealg{so}(n) \times \liealg{so}(k)
		\to
		\liealg{p}$
		is given by
		\begin{equation}
			\begin{split}
				\Omega^{\perp_X} 
				&=
				X X^{\top} \Omega 
				+ \Omega X X^{\top} 
				- 2 X \eta X^{\top},
				\\
				\eta^{\perp_X} &= -
				\big( \eta - X^{\top} \Omega X  \big) .
			\end{split}
		\end{equation}
		\item
		The map
		$\big(d_{(I_n, I_k)} (\iota_X \circ \pi)\at{\liealg{p}}\big)^{-1}
		\colon T_X \Stiefel{n}{k}
		\to
		\liealg{p}$
		is given by
		\begin{equation}
			V \mapsto \big(\Omega(V)^{\perp_X}, \eta(V)^{\perp_X} \big)=\big( V X^{\top} 
			- X V^{\top},   X^{\top} V   \big) .
	\end{equation}
	\end{enumerate}
	\begin{proof}
		This is a consequence of Lemma~\ref{lemma:stiefel_orthogonal_projection_on_p} and
		Lemma~\ref{lemma:stiefel_identification_tangent_space_p}.
	\end{proof}
\end{corollary}

\subsection{Intrinsic Rolling}
\label{subsec:stiefel_intrinsic_rolling}

In this section, using ideas from~\cite{schlarb:2023},
we apply Theorem~\ref{theorem:intrinsic_rolling}
to $\Stiefel{n}{k}$ equipped with an $\alpha$-metric.
More precisely, we use the isometry
\begin{equation}
	\label{equation:intrinsic_rolling_equivariant_diffeo}
	\iota_X \colon G / H \to \Stiefel{n}{k}
\end{equation}
to identify $\Stiefel{n}{k} \cong G / H$ as 
a normal naturally reductive homogeneous space
as well as the linear isometry
\begin{equation}
	\label{equation:intrinsic:rolling_isometry_tangent_spaces}
	\big(d_{(I_n, I_k)} (\iota_X \circ \pi)\at{\liealg{p}}\big)^{-1} \colon T_X \Stiefel{n}{k} \to \liealg{p}
\end{equation}
identifying $T_X \Stiefel{n}{k} \cong \liealg{p}$ as 
vector spaces equipped with the scalar product
from Subsection~\ref{subsec:stiefel_manifolds_alpha_metrics}.

Throughout this section, if not indicated otherwise, we always
assume that the maps 
from~\eqref{equation:intrinsic_rolling_equivariant_diffeo}
and~\eqref{equation:intrinsic:rolling_isometry_tangent_spaces}
are used to identify $G / H \cong \Stiefel{n}{k}$
and $\liealg{p} \cong T_X \Stiefel{n}{k}$, respectively.

These identifications allow the construction of an intrinsic rolling of 
$T_X \Stiefel{n}{k}$ over $\Stiefel{n}{k}$, 
where both manifolds are considered as
embedded into $\matR{n}{k}$.
We state the next definition in order to make this notion more precise.

{ Although, in the first part of this section, we have used the Greek letter $\alpha$ for the real parameter that defines a family of metrics  on Stiefel, the same letter will be used later for rolling curves. This will not create difficulties, since it will be clear from the context.}

\begin{definition}
	\label{definition:rolling_intrinsic_stiefel}
	Consider the Stiefel manifold $\Stiefel{n}{k} \subset \matR{n}{k}$,
	equipped with an $\alpha$-metric,
	as a submanifold of $\matR{n}{k}$.
	Moreover, let $X \in \Stiefel{n}{k}$ be fixed.
	Consider the triple $(\beta(t), \widehat{\beta(t)}, B(t))$, 
	where $\beta \colon I \to T_X \Stiefel{n}{k} \subset \matR{n}{k}$ 
	and
	$\widehat{\beta} \colon I \to \Stiefel{n}{k} \subset \matR{n}{k}$
	are curves and
	$B(t)
	\colon T_{\beta(t)} (T_X \Stiefel{n}{k}) \cong T_X \Stiefel{n}{k} 
	\to T_{\widehat{\beta}(t)} \Stiefel{n}{k}$
	is a linear isometry.
	This triple is called an intrinsic rolling of $T_X \Stiefel{n}{k}$
	over $\Stiefel{n}{k}$, with both manifolds embedded into $\matR{n}{k}$,
	if the following conditions hold:
	\begin{enumerate}
		\item
		no-slip condition:
		$\widehat{\beta}(t) = B(t) \dot{\beta}(t)$;
		\item
		no-twist condition:
		$B(t) Z(t)$ is a parallel vector field along $\widehat{\beta}(t)$ 
		iff $Z(t)$ is
		a parallel vector field along $\beta(t)$.
	\end{enumerate}
	The curve $\beta$ is called rolling curve and $\widehat{\beta}$ is called development curve.
\end{definition}
The next lemma uses Theorem~\ref{theorem:intrinsic_rolling}
to obtain a rolling of $T_X \Stiefel{n}{k}$ over $\Stiefel{n}{k}$ 
in the sense of Definition~\ref{definition:rolling_intrinsic_stiefel}.
 
 \begin{lemma}
 	\label{lemma:stiefel_intrinsic_rolling}
 	Let $\beta \colon I \to T_X \Stiefel{n}{k} \subset \matR{n}{k}$
 	be a curve and define the curve 
 	$\alpha \colon I \to \liealg{p}$ 
 	by $\alpha(t) = \big(d_{(I_n, I_k)} (\iota_X \circ \pi)\at{\liealg{p}}\big)^{-1} \big( \beta(t) \big)$
 	for $t \in I$.
 	Let $(\alpha(t), \widehat{\alpha}(t), A(t))$ be the triple
 	obtained in Theorem~\ref{theorem:intrinsic_rolling}
 	for the rolling  along $\alpha$ of $T_X \Stiefel{n}{k}$  
 	(identified with $\liealg{p}$), over $G / H$
 	(identified with $\Stiefel{n}{k}$).
 	Moreover, define the curve
 	\begin{equation}
 		\widehat{\beta} \colon I \to \Stiefel{n}{k},
 		\quad
 		t \mapsto \widehat\beta(t)  = \iota_X(\widehat{\alpha}(t) )
 	\end{equation}
 	and the isometry
 	$B(t) \colon T_{\beta(t)} (T_X \Stiefel{n}{k}) \cong T_X \Stiefel{n}{k} \to T_{\widehat\beta (t)} \Stiefel{n}{k}$
 	by
 	\begin{equation}
 		\begin{split}
 			B(t) = \big( d_{\widehat{\alpha}(t)} \iota_X \big) \circ A(t) \circ \big( d_{(I_n, I_k)} (\iota_X \circ \pi)^{-1} \big)
 		\end{split}.
 	\end{equation}
 	Then, the triple $(\beta(t), \widehat{\beta}(t), B(t))$ defines an intrinsic rolling of $T_X \Stiefel{n}{k}$ over $\Stiefel{n}{k}$
 	in the sense of Definition~\ref{definition:rolling_intrinsic_stiefel}.
 	\begin{proof}
 		The proof follows by applying Theorem~\ref{theorem:intrinsic_rolling}
 		since $G / H$ can be isometrically and
                $G$-equivariantly identified with $\Stiefel{n}{k}$ via
                $\iota_X \colon G / H \to \Stiefel{n}{k}$. Moreover, parallel vector fields are mapped to parallel vector fields by isometries.
 		
 		In more detail, the no-slip condition holds as
 		\begin{equation}
 			\begin{split}
 				\dot{\widehat{\beta}}(t)
 				&= 
 				\tfrac{\D}{\D t}
 				(\iota_X \circ \widehat{\alpha})(t) \\
 				&=
 				(d_{\widehat{\alpha}(t)} \iota_X ) \dot{\widehat{\alpha}}(t) \\
 				&=
 				(d_{\widehat{\alpha}(t)} \iota_X ) \big( A(t) \dot{\alpha}(t) \big) \\
 				&=
 				(d_{\widehat{\alpha}(t)} \iota_X ) \circ  A(t)
 				\circ \big(d_{(I_n, I_k)} (\iota_X \circ \pi)\at{\liealg{p}}\big)^{-1} \big( \dot{\beta}(t) \big) \\
 				&= 
 				B(t) \dot{\beta}(t).
 			\end{split}
 		\end{equation}
 		Next, we consider a parallel vector field
 		$V \colon I \to T (T_X \Stiefel{n}{k})$ along $\beta$,
 		i.e. $V$ can be viewed as the constant map
 		$V(t) = V_0$ for $t \in I$ and some
 		$V_0 \in T_X \Stiefel{n}{k}$.
 		Clearly,
 		$Z(t) = \big( d_{(I_n, I_k)} (\iota_X \circ \pi)^{-1} \big)(V(t)) = Z_0$
 		 is constant, with
 		$Z_0 = \big( d_{(I_n, I_k)} (\iota_X \circ \pi)^{-1} \big)(V_0)$,
 		i.e. $Z(t)$ is a parallel vector field along the curve
 		$\widehat{\alpha}$.
 		Thus, by Theorem~\ref{theorem:intrinsic_rolling}, the vector field $A(t) Z(t)$ is parallel along $\widehat{\alpha}$.
 		Since $\iota_X \colon G / H \to \Stiefel{n}{k}$ is an isometry, 
 		this parallel vector field is mapped to the parallel vector field
 		$d_{\widehat{\alpha}(t) } \iota_X \big(A(t) Z(t) \big)$
 		along the curve $\widehat{\beta}(t) = \iota_X(\widehat\alpha(t))$.
 		
 		Conversely, assuming that $d_{\widehat{\alpha}(t)} \iota_X (A(t) Z(t))$
 		is parallel along $\widehat{\beta}$, one shows by exploiting
 		Theorem~\ref{theorem:intrinsic_rolling} that $Z(t)$ is parallel along
 		$\widehat{\alpha}$ since $\iota_X^{-1} \colon \Stiefel{n}{k} \to G / H$
 		is an isometry.
 		Hence $V(t) = T_{(I_n, I_k)} (\iota_X \circ \pi)(Z(t))$ is parallel along $\beta$.
  	\end{proof}
 \end{lemma}
As a corollary, we reformulate the kinematic equations for the
intrinsic rolling of Stiefel manifolds
in the sense of Definition~\ref{definition:rolling_intrinsic_stiefel}.
\begin{corollary}
	\label{corollary:stiefel_kinematic_equations}
	Let $\beta \colon I \to \Stiefel{n}{k}$ be a curve and let
	$u \colon I \to \liealg{p}$ be the associated control curve, so that
	$u(t) =
	\big(d_{(I_n, I_k)} (\iota_X \circ \pi)\at{\liealg{p}}\big)^{-1}
	\big( \dot{\beta}(t) \big)$
	for $t \in I$.
	Consider the curves
	$S \colon I \to O(\liealg{p})$ as well as
	$q \colon I \ni t \mapsto q(t) = (R(t), \theta(t)) \in O(n) \times O(k)$
	defined by the initial value problems
	\begin{equation}
		\begin{split}
			\dot{S}(t) 
			&=
			- \tfrac{1}{2} \pr_{\liealg{p}} \circ \ad_{S(t) u(t)} \circ S(t) ,
			\quad S(0) = \id_{\liealg{p}} \\
			\dot{q}(t) 
			&= 
			\big(d_{(I_n, I_k)}L_{q(t)} \big) S(t) u(t),
			\quad
			q(0) = (I_n, I_k) .
		\end{split}
	\end{equation}
	Then the triple $(\beta(t), \widehat{\beta}(t), B(t))$
	defines an intrinsic rolling of $T_X \Stiefel{n}{k}$ over $\Stiefel{n}{k}$, 
	where
	\begin{equation}
		\widehat{\beta} \colon I \to \Stiefel{n}{k},
		\quad
		t \mapsto (\iota_X \circ \pi)(q(t)) = R(t) X \theta(t)^{\top}
	\end{equation}
	and 
	\begin{equation}
		B(t)
		=
		d_{q(t)} (\iota_X \circ \pi) \circ  d_e L_{q(t)} \circ S(t)  \circ  d_{(I_n, I_k)} (\iota_X \circ \pi)^{-1}.
	\end{equation}
	\begin{proof}
		This is a consequence of Lemma~\ref{lemma:stiefel_intrinsic_rolling}
		combined with Theorem~\ref{theorem:intrinsic_rolling}.
	\end{proof}
\end{corollary}

\subsection{Extrinsic Rolling}
\label{subsec:stiefel_extrinsic_rolling}

We now consider $\Stiefel{n}{k}$ embedded into $\matR{n}{k}$,
equipped with the metric induced by the Frobenius scalar product
scaled by the factor of two, i.e. 
the metric on $\Stiefel{n}{k}$ is given by 
\begin{equation}
	\langle V, W \rangle_X 
	=
	2 \tr\big(V^{\top} W \big), 
	\quad X \in \Stiefel{n}{k}, \ V, W \in T_X \Stiefel{n}{k} .
\end{equation}
This metric corresponds to the $\alpha$-metric, when $\alpha = - \tfrac{1}{2}$.
In the sequel, we will refer to this metric as the Euclidean metric.

We now construct a quadruple
$(\beta(t), \widehat{\beta}(t), B(t), C(t))$
which satisfies Definition~\ref{definition:rolling_extrinsic_1},
borrowing ideas from~\cite{schlarb:2023a}.

To this end, we first recall that a vector field
$\widehat{Z} \colon I \to N \Stiefel{n}{k}$ along a curve
$\widehat{\beta} \colon I \to  \Stiefel{n}{k}$ is normal parallel
if 
\begin{equation}
	\nabla_{\dot{\widehat{\beta}}(t)} \widehat{Z}(t) 
	=
	P_{\widehat{\beta}(t)}^{\perp}\big(\tfrac{\D}{\D t} \widehat{Z}(t)\big) 
	= 0,
	\quad t \in I
\end{equation}
holds, where $P_X^{\perp} \colon \matR{n}{k} \to N_X \Stiefel{n}{k}$
denotes the orthogonal projection onto the normal space 
$N_X \Stiefel{n}{k} = (T_X \Stiefel{n}{k})^{\perp}$
of $\Stiefel{n}{k}$ at the point $X$
with respect to the Euclidean metric.
This projection is given by
\begin{equation}
	P_X^{\perp}(V) = \tfrac{1}{2} X (X^{\top} V + V^{\top} X),
	\quad V \in \matR{n}{k}
\end{equation}
see e.g.~\cite{absil.mahony.sepulchre:2008}.

In order to determine the curve $\exIso \colon I \to O(N_X \Stiefel{n}{k})$,
we derive an ODE that is satisfied by a curve associated to
a normal vector field iff the vector field is parallel.
To this end, we first recall that
$\Phi_X = \iota_X \circ \pi \colon O(n) \times O(k) \to \Stiefel{n}{k}$
from~\eqref{equation:alpha_metrics_pseudo_Riemannian_submersion}
is a pseudo-Riemannian submersion.
Hence it makes sense to consider the horizontal lift of a
curve $\widehat{\beta} \colon I \to \Stiefel{n}{k}$.
In addtion,  for fixed 
$(\xi_1, \xi_2) \in \liealg{so}(n) \times \liealg{so}(k)$, we define the linear map
\begin{equation}
	\label{definiton:f_xi_xi_extrinsic_rolling_stiefel}
	f_{(\xi_1, \xi_2)} \colon \matR{n}{k} \to \matR{n}{k},
	\quad
	V \mapsto \xi_1 V - V \xi_2 .
\end{equation}
By this notation, we obtain the next lemma.
\begin{lemma}
	\label{lemma:normal_parallel_vector_fields_stiefel_euclidean_metric}
	Let $X \in \Stiefel{n}{k}$ be fixed,
	 $\widehat{\beta} \colon I \to \Stiefel{n}{k}$
	a curve and  $\widehat{Z} \colon I \to N \Stiefel{n}{k}$ a normal
	vector field along $\widehat\beta$.
	Moreover, let
	$q \colon I \ni t \mapsto q(t) = (R(t), \theta(t)) \in O(n) \times O(k)$
	be a horizontal
	lift of $\widehat{\beta}$.
	Then $\widehat{Z}$ is parallel along $\widehat{\beta}$ iff 
	the curve 
	\begin{equation}
		z^{\perp} \colon I \to N_X \Stiefel{n}{k},
		\quad
		t \mapsto z^{\perp}(t) 
		=
		\Phi_{q(t)^{-1}}\big(\widehat{Z}(t) \big)
		=
		R(t)^{\top} \widehat{Z}(t) \theta(t)
	\end{equation}
	satisfies the ODE 
	\begin{equation}
		\dot{z}^{\perp}(t)
		=
		-
		\big(P_X^{\perp} \circ f_{(\xi_1(t), \xi_2(t))} \big)\big(z^{\perp}(t) \big),
		\quad t \in I ,
	\end{equation}
	where $(\xi_1(t), \xi_2(t)) = \big(R(t)^{\top} \dot{R}(t), \theta(t)^{\top} \dot{\theta}(t) \big) \in \liealg{so}(n) \times \liealg{so}(k)$.
	\begin{proof}
		Let $(R, \theta) \in O(n) \times O(k)$ and $X \in \Stiefel{n}{k}$.
		Then 
		\begin{equation}
			\label{equation:orthogonal_projection_normal_bundle_invariance}
			P^{\perp}_{\Phi_{(R, \theta)}(X)} (V)
			=
			\Phi_{(R, \theta)}
			\circ P_X^{\perp} \circ 
			\Phi_{(R^{\top}, \theta^{\top})}(V)
		\end{equation}
		holds for $V \in \matR{n}{k}$ by the $\Phi$-invariance of the Euclidean metric.
		Since $q(t) = (R(t), \theta(t))$ is a horizontal lift of $\widehat{\beta}$, 
		i.e.
		$\widehat{\beta}(t) 
		= (\iota_X \circ \pi)(q(t)) = R(t) X \theta(t)^{\top}$,
		\eqref{equation:orthogonal_projection_normal_bundle_invariance} implies
		\begin{equation}
			\label{equation:parallel_normal_vector_field_1}
			\begin{split}
				P^{\perp}_{\widehat{\beta}(t)}(V)
%				=
%				\Phi_{(R(t), \theta(t))}
%				\circ P_X^{\perp} \circ 
%				\Phi_{(R(t)^{\top}, \theta(t)^{\top})}(V)
				= 
				\Phi_{(R(t), \theta(t))}
				\circ P_X^{\perp}
				\big( R(t)^{\top} V \theta(t) \big) .
			\end{split}
		\end{equation}
		Moreover, the condition
		$P_{\widehat{\beta}(t)}^{\perp}\big(\tfrac{\D}{\D t} \widehat{Z}(t)\big) = 0$
		is equivalent to
		\begin{equation}
			\label{equation:parallel_normal_vector_field_2}
			P_X^{\perp} \big( R(t)^{\top} \big(\tfrac{\D}{\D t} \widehat{Z}(t) \big) \theta(t) \big) = 0
		\end{equation}
		by~\eqref{equation:parallel_normal_vector_field_1}
		since $\Phi_{(R(t), \theta(t))} \colon \matR{n}{k} \to \matR{n}{k}$
		is a linear isomorphism.
		Obviously, by the definition of $z^{\perp}$, we have
		\begin{equation}
			\label{equation:parallel_normal_vector_field_z_Z_relation}
			\widehat{Z}(t) = R(t) z^{\perp}(t) \theta(t)^{\top}.
		\end{equation}
		Plugging~\eqref{equation:parallel_normal_vector_field_z_Z_relation} into~\eqref{equation:parallel_normal_vector_field_2} yields
		\begin{equation}
			\label{equation:parallel_normal_vector_field_3}
			\begin{split}
				0
				&=
				P_X^{\perp} \big( R(t)^{\top} \big(\tfrac{\D}{\D t} \big(  R(t) z^{\perp}(t) \theta(t)^{\top}  \big) \big) \theta(t) \big)  \\
				&=
				P_X^{\perp} \big( R(t)^{\top} \big( 
				 \dot{R}(t) z^{\perp}(t) \theta(t)^{\top}
				 + R(t) \dot{z^{\perp}}(t) \theta(t)^{\top}
				 + R(t) z^{\perp}(t) \dot{\theta}(t)^{\top}  \big) \theta(t) \big) \\
				 &=
				 P_X^{\perp}  \big( R(t)^{\top} \dot{R}(t) z^{\perp}(t)
				 + \dot{z^{\perp}}(t) 
				 + z^{\perp}(t) \dot{\theta}(t)^{\top} \theta(t) \big) .
			\end{split}
		\end{equation}
		Using 
		$(\xi_1(t), \xi_2(t)) 
		= \big(R(t)^{\top} \dot{R}(t), \theta(t)^{\top} \dot{\theta}(t) \big)$
		and 
		$\theta(t)^{\top} \dot{\theta}(t) 
		= - \dot{\theta}(t)^{\top} \theta(t)$
		as well as $P_X^{\perp}\big(\dot{z^{\perp}}(t) \big) = \dot{z^{\perp}}(t)$
		due to $z^{\perp}(t) \in N_X \Stiefel{n}{k}$,
		we can equivalently
		rewrite~\eqref{equation:parallel_normal_vector_field_3}
		by 
		\begin{equation}
			\begin{split}
				0
				&= 
				P_X^{\perp}\big( \xi_1(t) z^{\perp}(t) 
				+ \dot{z^{\perp}}(t) 
				- z^{\perp}(t) \xi_2(t) \big) \\
				&=
				\dot{z^{\perp}}(t) + \big( P_X^{\perp} \circ f_{(\xi_1(t), \xi_2(t))} \big) \big(z^{\perp}(t) \big)  .
			\end{split}
		\end{equation}
		This yields the desired result.	
	\end{proof}
\end{lemma}
After this preparation, we are in the position to determine
the extrinsic rolling of $T_X \Stiefel{n}{k}$ over
$\Stiefel{n}{k}$ with respect to the Euclidean metric
in the sense of Definition~\ref{definition:rolling_extrinsic_1}.

\begin{theorem}
	\label{theorem:stiefel_rolling_extrinsic_euclidean_metric}
	Let $X \in \Stiefel{n}{k}$ be fixed and let { $\beta \colon I \to T_X\Stiefel{n}{k}$} be a curve.
	Moreover, let $(\beta(t), \widehat{\beta}(t), B(t))$ denote the
	intrinsic rolling of $T_X \Stiefel{n}{k}$ over $\Stiefel{n}{k}$
	from Lemma~\ref{lemma:stiefel_intrinsic_rolling}
	for $\alpha = - \tfrac{1}{2}$.
	Furthermore, let 
	$q \colon I \ni t \mapsto q(t) = (R(t), \theta(t)) \in O(n) \times O(k)$
	be the horizontal lift of $\widehat{\beta} \colon I \to \Stiefel{n}{k}$
	through $q(0) = (I_n, I_k)$
	and define
	$(\xi_1, \xi_2) \colon I \to \liealg{so}(n) \times \liealg{so}(k)$ by
	\begin{equation}
		(\xi_1(t), \xi_2(t)) 
		=
		 \big(d_{(I_n, I_k)} L_{q(t)} \big)^{-1} \dot{q}(t) 
		=
		 \big(R(t)^{\top} \dot{R}(t), \theta(t)^{\top} \dot{\theta}(t) \big)
	\end{equation}
	for $t \in I$.
	Let $\exIso \colon I \to O(N_X \Stiefel{n}{k})$ be the solution of
	the initial value problem
	\begin{equation}
		\label{equation:theorem_stiefel_rolling_extrinsic_euclidean_metric_ODE_for_P}
		\dot{\exIso}(t) 
		= 
		- P_X^{\perp} \circ f_{(\xi_1(t), \xi_2(t))} \circ \exIso(t),
		\quad
		\exIso(0) = \id_{N_X \Stiefel{n}{k}}.
	\end{equation}
	Then the quadruple $(\beta(t), \widehat{\beta}(t), B(t), C(t))$, 
	with
	\begin{equation}
		C(t) \colon N_{\beta(t)}(T_X \Stiefel{n}{k}) 
		\cong N_X \Stiefel{n}{k} \to N_{\widehat{\beta}(t)} \Stiefel{n}{k}
	\end{equation}
	defined by
	\begin{equation}
		C(t) = \Phi_{(R(t), \theta(t))} \circ \exIso(t),
	\end{equation}
	is an extrinsic rolling of $T_X \Stiefel{n}{k}$ over $\Stiefel{n}{k}$
	with respect to the Euclidean metric.
	\begin{proof}
		We only need to show the normal no-twist condition
		since the tangential no-twist condition and the no-slip
		condition are fulfilled
		by Lemma~\ref{lemma:stiefel_intrinsic_rolling}.
		We start with proving that $\exIso(t) \in O(N_X \Stiefel{n}{k})$
		for $t \in I$.
		For that, we compute
		\begin{equation}
			\begin{split}
				\big\langle ( - P_X^{\perp} \circ f_{(\xi_1(t), \xi_2(t)) } )(Y), Z \big\rangle _X 
				&= 
				- \big\langle f_{(\xi_1(t), \xi_2(t)) } (Y), Z \big\rangle_X \\
				&=
				- 2 \tr\big( \big(\xi_1(t) Y - Y \xi_2(t) \big)^{\top} Z \big) \\
				&=
				2 \tr\big( Y^{\top} \xi_1(t) Z - \xi_2(t) Y^{\top} Z \big) \\
				&=
				2 \tr\big( Y^{\top} \big(\xi_1(t) Z - Z \xi_2(t) \big) \big) \\
				&= 
				\big\langle Y, (P_X^{\perp} \circ f_{(\xi_1(t), \xi_2(t))} ) (Z) \big\rangle
			\end{split}
		\end{equation}
		for $Y, Z \in N_X \Stiefel{n}{k}$ by exploiting $(\xi_1(t), \xi_2(t))  \in \liealg{so}(n) \times \liealg{so}(k)$.
		Thus, $ - P_X^{\perp} \circ f_{(\xi_1(t), \xi_2(t)) } 
		\colon N_X \Stiefel{n}{k} \to \Stiefel{n}{k}$
		is skew-adjoint with respect to the Euclidean metric,
		implying that
		$- P_X^{\perp} \circ f_{(\xi_1(t), \xi_2(t)) } ) \circ \exIso$,
		for $\exIso \in O(N_X \Stiefel{n}{k})$, can be viewed as 
		a time-variant vector field on $O(N_X \Stiefel{n}{k})$.
		
		Next we note that
		$C(t) \colon N_{\beta(t)}(T_X \Stiefel{n}{k}) 
		\cong N_X \Stiefel{n}{k} \to N_{\widehat{\beta}(t)} \Stiefel{n}{k}$
		is an isometry (as composition of isometries).
		Now, let $Z^{\perp} \colon I \to N(T_X \Stiefel{n}{k})$ be
		a normal parallel vector field along
		$\beta \colon I \to T_X \Stiefel{n}{k}$.
		Then $Z^{\perp}$ can be viewed as the constant curve
		$Z^{\perp}(t) = Z_0^{\perp}$, for $t \in I$ and some
		$Z_0^{\perp} \in N_X \Stiefel{n}{k}$.
		Obviously, $\widehat{Z}^{\perp} \colon I \to N \Stiefel{n}{k}$ 
		given by
		\begin{equation}
			\widehat{Z}^{\perp}(t) 
			=
			C(t) Z^{\perp}(t)
			=
			 \big(\Phi_{(R(t), \theta(t))} \circ \exIso(t) \big) (Z_0^{\perp}),
			\quad
			t \in I,
		\end{equation}
		is a normal vector field along the curve $\widehat{\beta}$.
		It remains to show that $\widehat{Z}^{\perp}$ is parallel
		along $\widehat{\beta}$.
		To this end, we exploit
		Lemma~\ref{lemma:normal_parallel_vector_fields_stiefel_euclidean_metric}.
		We consider the curve 
		$z^{\perp} \colon I \to N_X \Stiefel{n}{k}$ given by
		\begin{equation}
			z^{\perp}(t) 
			=
			\Phi_{(R(t)^{\top}, \theta(t)^{\top})} \big(\widehat{Z}^{\perp}(t) \big)
			= 
			\exIso(t) (Z_0^{\perp}) 
		\end{equation}
		and obtain
		\begin{equation}
			\begin{split}
				\dot{z^{\perp}}(t) 
				&=
				\dot{\exIso}(t)(Z_0^{\perp}) \\
				&=
				- \big(P_X^{\perp} \circ f_{(\xi_1(t), \xi_2(t))} \circ \exIso(t) \big)(Z_0^{\perp}) \\
				&= 
				- \big(P_X^{\perp} \circ f_{(\xi_1(t), \xi_2(t))} \big) \big(z^{\perp}(t) \big)
			\end{split}
		\end{equation}
		due
		to~\eqref{equation:theorem_stiefel_rolling_extrinsic_euclidean_metric_ODE_for_P}.
		Thus $\widehat{Z}^{\perp}$ is parallel along $\widehat{\beta}$ by
		Lemma~\ref{lemma:normal_parallel_vector_fields_stiefel_euclidean_metric}.
		
		Conversely, assume that $\widehat{Z}^{\perp} \colon I \to N\Stiefel{n}{k}$ given by
		$Z^{\perp}(t) = C(t) Z(t)^{\perp}$ for some
		$Z^{\perp} \colon I \to N_X \Stiefel{n}{k}$ is normal parallel
		along $\widehat{\beta}$.
		We define the normal parallel frame along $\widehat{\beta}$ by 
		$A^{\perp}_i(t) = C(t) A_i$, where the vectors $A_i^{\perp} \in N_X \Stiefel{n}{k}$
		for $i \in \{1, \ldots, \ell_n\}$
		with
		$\ell_n = \dim(N_X \Stiefel{n}{k})$
		form a basis.
		Then, analogously to~\cite[Chap. 4, p. 106]{Lee:2018aa},
		one shows that $\widehat{Z}^{\perp}$ is normal parallel along
		$\widehat{\beta}$ iff the coefficient functions
		$z^i \colon I \to \field{R}$ defined by
		$\widehat{Z}^{\perp}(t)
		= 
		\sum_{i = 1}^{\ell_n} z_i(t) A_i^{\perp}(t)$
		are constant.
		Since $\widehat{Z}^{\perp}$ is assumed to be normal parallel,
		there exists a uniquely determined $z_i \in \field{R}$ such that $Z^{\perp}(t) 
		=
		\sum_{i = 1}^{\ell_n} z_i A_i^{\perp}(t)$
		is fulfilled.
		Hence, by the linearity of
		$C(t) \colon N_{\beta(t)}(T_X \Stiefel{n}{k}) \cong N_X \Stiefel{n}{k}
		\to N_{\widehat{\beta}(t)} \Stiefel{n}{k}$, we obtain
		\begin{equation}
			\widehat{Z}^{\perp}(t) 
			=
	 		\sum_{i = 1}^{\ell_n}
			z_i A_i^{\perp}(t) 
			=
	 		\sum_{i = 1}^{\ell_n}
			z_i C(t) A_i^{\perp}
			=
			\sum_{i = 1}^{\ell_n}
			C(t) (z^i A_i^{\perp})
			=
			C(t) Z^{\perp},
		\end{equation}
		where
		$Z^{\perp}
		=
		\sum_{i = 1}^{\ell_n} z_i A_i^{\perp}$
		is viewed as a normal
		vector field along $\beta$ which is clearly normal parallel.
		This yields the desired result.
	\end{proof}
 \end{theorem}
As a corollary of Theorem~\ref{theorem:stiefel_rolling_extrinsic_euclidean_metric}
we obtain the kinematic equations for the extrinsic rolling of
$T_X \Stiefel{n}{k}$ over $\Stiefel{n}{k}$ with respect
to the Euclidean metric.

\begin{corollary}
	\label{corollary:extrinsic_rolling_stiefel_kinematic_equations}
	Let $X \in \Stiefel{n}{k}$ be fixed and let { $\beta \colon I \to T_X\Stiefel{n}{k}$} be a prescribed rolling curve
	with associated control curve 
	\begin{equation}
		u \colon I \ni t \mapsto \big(d_{(I_n, I_k)} (\iota_X \circ \pi) \at{\liealg{p}} \big)^{-1}(\dot{\beta}(t)) \in \liealg{p}
	\end{equation}
	viewed as a curve in $\liealg{p}$,
	where 
	\begin{equation}
		\big(d_{(I_n, I_k)} (\iota_X \circ \pi) \at{\liealg{p}} \big)^{-1} \colon T_X \Stiefel{n}{k} \to \liealg{p}
	\end{equation}
	is given by
	Corollary~\ref{corollary:stiefel_orthogonal_projection_identification_euclidean_metric}.
	Moreover, let the curves $S \colon I \to O(\liealg{p})$ and
	$q \colon I \to O(n) \times O(k)$, as well as
	$\exIso \colon I \to O(N_X \Stiefel{n}{k})$,
	be defined by the initial value problem
	\begin{equation}
		\label{equation:kinematic_equations_stiefel_extrinsic}
		\begin{split}
			\dot{S}(t) 
			&=
			- \tfrac{1}{2} \pr_{\liealg{p}} \circ \ad_{S(t) u(t)} \circ S(t) ,
			\quad S(0) = \id_{\liealg{p}} , \\
			\dot{q}(t) 
			&= 
			\big(d_{(I_n, I_k)}L_{q(t)} \big) S(t) u(t),
			\quad
			q(0) = (I_n, I_k) , \\
			\dot{\exIso}(t) 
			&= 
			- P_X^{\perp} \circ f_{(\xi_1(t), \xi_2(t))} \circ \exIso(t),\quad
			\exIso(0) = \id_{N_X \Stiefel{n}{k}} , \\
		\end{split}
	\end{equation}
	where $f_{(\xi_1, \xi_2)} \colon \matR{n}{k} \to \matR{n}{k}$
	is given by~\eqref{definiton:f_xi_xi_extrinsic_rolling_stiefel}
	and $\pr_{\liealg{p}} \colon \liealg{so}(n) \times \liealg{so}(k) \to \liealg{p}$ is determined in Corollary~\ref{corollary:stiefel_orthogonal_projection_identification_euclidean_metric}.
	Then,  $(\beta(t), \widehat{\beta}(t), B(t), C(t))$
	defines an extrinsic rolling of $T_X \Stiefel{n}{k}$ over $\Stiefel{n}{k}$
	with respect to the Euclidean metric,
	where
	\begin{equation}
		\widehat{\beta} \colon I \to \Stiefel{n}{k},
		\quad
		t \mapsto (\iota_X \circ \pi)(q(t)) = R(t) X \theta(t)^{\top},
	\end{equation}
	\begin{equation}
		B(t)
		=
		d_{(q(t))} (\iota_X \circ \pi) \circ \big( (d_e L_{q(t)}) \circ S(t) \big) \circ d_{(I_n, I_k)} (\iota_X \circ \pi)^{-1},
	\end{equation}
	and 
	\begin{equation}
		C(t) = \Phi_{(R(t), \theta(t))} \circ \exIso(t) .
	\end{equation}
	We call~\eqref{equation:kinematic_equations_stiefel_extrinsic}
	the kinematic equations for the extrinsic rolling of $T_X \Stiefel{n}{k}$
	over $\Stiefel{n}{k}$ with respect to the Euclidean metric.
\end{corollary}

\subsection{Rolling Along Special Curves}
\label{subsec:rolling_along_special_curves}

In this subsection
we consider a rolling of $T_X \Stiefel{n}{k}$
over $\Stiefel{n}{k}$ such that its development curve
$\widehat{\beta} \colon I \to \Stiefel{n}{k}$ is the projection
of a \emph{not} necessarily horizontal one-parameter subgroup,
i.e. a curve
\begin{equation}
	\label{equation:rolling_along_special_curves_development_curve}
	\widehat{\beta} \colon I \to \Stiefel{n}{k},
	\quad
	t \mapsto (\iota_X \circ \pi)(\exp(t \xi))
	= 
	\e^{t \xi_1} X \e^{-t \xi_2}
\end{equation}
for some
$(\xi_1, \xi_2) \in \liealg{so}(n) \times \liealg{so}(k)$,
where $X \in \Stiefel{n}{k}$ is fixed.
For this special case, which includes the curves considered in~\cite{kakowski.machado.leite:2017},
we determine an extrinsic rolling $(\beta(t), \widehat{\beta}(t), B(t), C(t))$ 
explicitly. 
To this end, we proceed as in~\cite{schlarb:2023},
where the intrinsic rolling of general reductive spaces 
along such a curve are determined explicitly.
However, for the following discussion, we will restrict  to the study of Stiefel manifolds
equipped with the Euclidean metric, as it
allows  simplifying some arguments.

Before we continue, we fix some notations.
Let $\xi = (\xi_1, \xi_2) \in 
\liealg{so}(n) \times \liealg{so}(k)$.
Let
$\xi_{\liealg{h}} = (\xi_{1, \liealg{h}}, \xi_{2, \liealg{h}})$
and $\xi_{\liealg{p}} = (\xi_{1, \liealg{p}}, \xi_{2, \liealg{p}})$
denote the projections of $\xi$ onto $\liealg{h}$ and onto $\liealg{p}$,
respectively.
Here, the reductive decomposition is always understood
to be taken with respect to the $\alpha$-metric,
where $\alpha = - \tfrac{1}{2}$.
In particular, the subspaces $\liealg{h}$ and $\liealg{p}$ 
	of $\liealg{so}(n) \times \liealg{so}(k)$ are
	orthogonal with respect to the
	scalar product
	$\langle \cdot, \cdot \rangle^{\alpha}_{\liealg{so}(n) \times \liealg{so}(k)}$
	defined in~\eqref{equation:alpha_metrics_lie_algebra_def}.

We first consider the horizontal lift of a curve given
by~\eqref{equation:rolling_along_special_curves_development_curve}.

\begin{lemma}
	\label{lemma:horizontal_lift_special_curve}
	Let $X \in \Stiefel{n}{k}$ and $\xi = (\xi_1, \xi_2) \in 
	\liealg{so}(n) \times \liealg{so}(k)$.
	The horizontal lift of
	\begin{equation}
		\widehat{\beta} \colon I \to \Stiefel{n}{k},
		\quad
		t \mapsto \widehat{\beta}(t) = \big(\iota_X \circ \pi \big)\big(\exp(t \xi) \big)
		= \e^{t \xi_1} X \e^{- t \xi_2}
	\end{equation}
	through $q(0) = (I_n, I_k)$ is given by
	\begin{equation}\label{horiz-curve-F}
		\begin{split}
			q \colon I 
			&\to O(n) \times O(k), \\
			t 
			&\mapsto \exp(t \xi) \exp(- t \xi_{\liealg{h}})
			=
			\big(\e^{t \xi_1} e^{- t \xi_{1, \liealg{h}}}, \e^{t \xi_2} \e^{- t \xi_{2, \liealg{h}}} \big) .
		\end{split}
	\end{equation}
	Moreover, it is the solution of the initial value problem
	\begin{equation}
		\label{equation_lemma_horizontal_lift_special_curve_IVP}
		\dot{q}(t) 
		=
		\big(d_{(I_n, I_k)} L_{q(t)} \big) \Ad_{\exp( t \xi_{\liealg{h}})}(\xi_{\liealg{p}}),
		\quad
		q(0) = (I_n, I_k) .
	\end{equation}
	\begin{proof}
		Obviously, $q(0) = (I_n, I_k)$ holds and
		$\widehat{\beta}(t) 
		=
		\big(\iota_X \circ \pi \big)(\exp(t \xi)) 
		=
		\big(\iota_X \circ \pi \big) \big(\exp(t \xi) \exp(-t \xi_{\liealg{h}}) \big)$ 
		is fulfilled
		since $t \mapsto \exp(- t \xi_{\liealg{h}})$ 
		is a curve in $H \subset O(n) \times O(k)$.
		
		We claim that $q$ is horizontal.
		Indeed, 
		by using the well-known properties
		of matrix exponential
		\begin{equation}
			\tfrac{\D}{\D t} \exp(t \xi) = \exp(t \xi) \xi
			\quad \text{ and } \quad
			\tfrac{\D}{\D t} \exp(t \xi_{ \liealg{h}}) = \xi_{ \liealg{h}} \exp(t \xi_{ \liealg{h}}),
	\end{equation}
	we compute
		\begin{equation}
			\dot{q}(t) 
			= 
			\exp(t \xi) \xi \exp( - t \xi_{ \liealg{h}}) 
			- \exp( t \xi) \xi_{\liealg{h}} \exp( - t \xi_{ \liealg{h}})
		\end{equation}
		yielding
		\begin{equation}
			\label{equation_lemma:horizontal_lift_special_curve_computation}
			\begin{split}
				(d_{(I_n, I_k)} L_{q(t)} )^{-1} \dot{q}(t)
				&=
				\exp( t \xi_{\liealg{h}}) \exp(- t \xi) \dot{q}(t)\\
				&=
				\exp( t \xi_{ \liealg{h}}) \exp(- t \xi)
				\big(\exp(t \xi) \xi \exp( - t \xi_{ \liealg{h}}) \\
				&\quad
				- \exp( t \xi) \xi_{\liealg{h}} \exp(t \xi_{\liealg{h}})\big) \\
				&= 
				\exp(t \xi_{\liealg{h}}) \xi \exp(- t \xi_{\liealg{h}}) \\
				&\quad
				- \exp(t \xi_{\liealg{h}}) \xi_{ \liealg{h}} \exp( -t \xi_{ \liealg{h}}) \\
				&= 
				\exp(t \xi_{ \liealg{h}}) \xi_{\liealg{p}}\exp( -t  \xi_{\liealg{h}} ) \\
				&= 
				\Ad_{\exp( t \xi_{\liealg{h}})} (\xi_{\liealg{p}}).
			\end{split}
		\end{equation}
		Here we exploited
		that $O(n) \times O(k)$ can be viewed as a matrix Lie group.
		Hence $q \colon I \to O(n) \times O(k)$ is horizontal
		due to $\Ad_{\exp(t \xi_{\liealg{h}}}(\xi_{\liealg{p}}) \in \liealg{p}$ 
		since 
		$\liealg{so}(n) \times \liealg{so}(k) = \liealg{h} \oplus \liealg{p}$
		is a reductive decomposition.
		In addition,~\eqref{equation_lemma:horizontal_lift_special_curve_computation}
		implies that $q$ is the solution of~\eqref{equation_lemma_horizontal_lift_special_curve_IVP},
		as desired.
	\end{proof}
\end{lemma}

Next we determine the intrinsic rolling
$(\alpha(t), \widehat{\alpha}(t), A(t))$ of $T_X \Stiefel{n}{k} \cong \liealg{p}$ over $\Stiefel{n}{k} \cong (O(n) \times O(k)) / H$ 
viewed as a normal naturally reductive homogeneous space,
where $\widehat{\alpha}(t) = \pi(\exp(t \xi))$ for some
$\xi \in \liealg{so}(n) \times \liealg{so}(k)$.

To this end, we recall the kinematic equations from Theorem~\ref{theorem:intrinsic_rolling}.
They are given by 
\begin{equation}
	\label{equation:kinematic_equation_intrinsic_rolling_stiefel_special_curves}
	\begin{split}
		\dot{S}(t) 
		&=
		- \tfrac{1}{2} \pr_{\liealg{p}} \circ \ad_{S(t) u(t)} \circ S(t) ,
		\quad S(0) = \id_{\liealg{p}} \\
		\dot{q}(t) 
		&= 
		\big(d_{(I_n, I_k)} L_{q(t)} \big) S(t) u(t),
		\quad
		q(0) = (I_n, I_k) ,
	\end{split} 
\end{equation}
where 
\begin{equation}
	\label{equation:rolling_along_special_curves_S_u_Ad}
	S(t) u(t) =\Ad_{\exp(t \xi_{ \liealg{h}})} (\xi_{\liealg{p}})
\end{equation}
for $t \in I$ by the definition of
$\widehat{\alpha}(t) = \pi(\exp(t \xi)) = (\iota_X)^{-1}(\widehat{\beta}(t))$
and Lemma~\ref{lemma:horizontal_lift_special_curve}.
Thus, the ODE for $S \colon I \to O(\liealg{p})$ in~\eqref{equation:kinematic_equation_intrinsic_rolling_stiefel_special_curves}
becomes 
\begin{equation}
	\dot{S}(t)
	%	=
	%	- \tfrac{1}{2} \pr_{\liealg{p}} \circ \ad_{S(t) u(t)} \circ S(t)
	=
	- \tfrac{1}{2} \pr_{\liealg{p}} \circ \ad_{\Ad_{\exp(t \xi_{ \liealg{h}})} (\xi_{\liealg{p}})} \circ S(t),
	\quad S(0) = \id_{\liealg{p}} .
\end{equation}
In order to determine the intrinsic rolling explicitly,
we need to solve this equation.

As a preparation, we state a lemma on time-variant linear
ODEs which is inspired
by~\cite[p. 48]{brockett:1970}.

\begin{lemma}
	\label{lemma:ode_time_variant_linear_special_solution}
	Let $V$ be a finite dimensional  real vector space and let
	$A, B \in \End(V)$ be linear maps on $V$.
	Consider the curve $S \colon I \to GL(V)$ defined by the initial value problem
	\begin{equation}
		\dot{S}(t) = \exp(t A) \circ B \circ \exp(-t A) \circ S(t),
		\quad
		S(0) = S_0 \in GL(V) .
	\end{equation}
	Then $S$ is given by
	\begin{equation}
		S(t) = \exp(t A) \circ \exp(t (B-A))\circ S_0 .
	\end{equation}
	\begin{proof}
		Define $\widetilde{S} \colon I \to GL(V)$ 
		by $\widetilde{S}(t) = \exp(-t A) \circ S(t)$.
		Then
		\begin{equation}
			\begin{split}
				\dot{\widetilde{S}}(t)
				&= 
				- A \circ \exp(-t A) \circ S(t) 
				+ \exp(-t A) \circ \dot{S}(t) \\
				&=
				- A \circ \widetilde{S}(t) + \exp(-t A) \circ \exp(t A) \circ B \circ \exp(-t A) \circ S(t) \\
				&=
				(B - A) \circ \widetilde{S}(t) 
			\end{split}
		\end{equation}
		for $t \in I$ implying $\widetilde{S}(t) = \exp(t (B - A)) \circ S_0$.
		Consequently, by the  definition of
                $\widetilde{S}$, we obtain
		\begin{equation}
			S(t) 
			=
			\exp(t A) \circ \widetilde{S}(t) 
			=
			\exp(t A) \circ \exp(t (B - A)) \circ S_0,
			\quad t \in I.
		\end{equation}
	\end{proof}
\end{lemma}

\begin{lemma}
	\label{lemma:ode_rolling_along_special_curves_intrinsic_solution}
	Let $\xi \in \liealg{so}(n) \times \liealg{so}(k)$.
	The solution of the initial value problem 
	\begin{equation}
		\label{equation:lemma_ode_rolling_along_special_curves_intrinsic_solution_1}
		\dot{S}(t)
		=
		- \tfrac{1}{2} \pr_{\liealg{p}} \circ \ad_{\Ad_{\exp(t \xi_{ \liealg{h}})} (\xi_{\liealg{p}})} \circ S(t),
		\quad S(0) = \id_{\liealg{p}} 
	\end{equation}
	is given by 
	\begin{equation}
		S \colon I \to O(\liealg{p}),
		\quad
		t \mapsto
		\Ad_{\exp(t \xi_{\liealg{h}} )} \circ 
		\exp\big(- t \big( \ad_{\xi_{ \liealg{h}}} + \tfrac{1}{2} (\pr_{\liealg{p}} \circ \ad_{\xi_{\liealg{p}}}) \big) \big) .
	\end{equation}
	\begin{proof}
		Rewrite~\eqref{equation:lemma_ode_rolling_along_special_curves_intrinsic_solution_1}
		such that Lemma~\ref{lemma:ode_time_variant_linear_special_solution}
		can be applied.
		We compute
		\begin{equation}
			\begin{split}
				\dot{S}(t)
				&=
				- \tfrac{1}{2} \pr_{\liealg{p}} \circ \ad_{\Ad_{\exp(t \xi_{ \liealg{h}})} (\xi_{\liealg{p}})} \circ S(t) \\
				&=
				- \tfrac{1}{2} 
				\Ad_{\exp(t \xi_{ \liealg{h}})} \circ \pr_{\liealg{p}} \circ \ad_{\xi_{\liealg{p}}} \circ \Ad_{\exp( -t \xi_{\liealg{h}})} \circ S(t)
				\\
				&=
				- \tfrac{1}{2} 
				\exp(t \ad_{\xi_{ \liealg{h}}} )\circ \pr_{\liealg{p}}
				\circ \ad_{\xi_{\liealg{p}}} \circ \exp(-t \ad_{\xi_{ \liealg{h}}} ) \circ S(t),
			\end{split}
		\end{equation}
		where, in the first equality, we used 
		$\Ad_{\exp(t \xi_{ \liealg{h}})} \colon \liealg{g} \to
                \liealg{g}$ being a Lie algebra morphism, and moreover, that
		$\Ad_h \circ \pr_{\liealg{p}} 
		= \pr_{\liealg{p}} \circ \Ad_h$
		holds due to
		$\Ad_h(\liealg{p}) \subset \liealg{p}$ as well as
		$\Ad_h(\liealg{h}) \subset \liealg{h}$
		for $h \in H$.
		For the second equality, 		$\Ad_{\exp(t \xi_{ \liealg{h}})} = \exp(\ad_{t \xi_{ \liealg{h}}})$
		is used.
		Hence, we can apply Lemma~\ref{lemma:ode_time_variant_linear_special_solution},
		where we set $A = \ad_{\xi_{ \liealg{h}}}$ and $B = - \tfrac{1}{2}  \pr_{\liealg{p}} \circ \ad_{\xi_{\liealg{p}}}$.
		This yields
		\begin{equation}
			\begin{split}
				S(t) 
				&=
				\exp(t \ad_{\xi_{ \liealg{h}}}) \circ 
				\exp\big( t \big( - \tfrac{1}{2}  \pr_{\liealg{p}} \circ \ad_{\xi_{\liealg{p}}} 
					- \ad_{\xi_{ \liealg{h}}} \big) \big) \circ \id_{\liealg{p}} \\
				&=
				\Ad_{\exp(t \xi_{\liealg{h}})} 
				\circ
				\exp\big({- t \big( \ad_{\xi_{ \liealg{h}}} + \tfrac{1}{2}  \pr_{\liealg{p}} \circ \ad_{\xi_{\liealg{p}}} \big)} \big)
			\end{split}
		\end{equation}
		as desired.
	\end{proof}
\end{lemma}

We proceed with determining
the intrinsic rolling $(\alpha(t), \widehat{\alpha}(t), A(t))$.
Recall that the control curve
$u \colon I \to \liealg{p}$ is defined by $u(t) = \dot{\alpha}(t)$.
Hence~\eqref{equation:rolling_along_special_curves_S_u_Ad} yields
\begin{equation}
	\dot{\alpha}(t) 
	= 
	S(t)^{-1} \Ad_{\exp(t \xi_{ \liealg{h}})}(\xi_{\liealg{p}})
	=
	\exp\big( t \big( \ad_{\xi_{ \liealg{h}}} + \tfrac{1}{2} (\pr_{\liealg{p}} \circ \ad_{\xi_{\liealg{p}}}) \big) \big) (\xi_{\liealg{p}}),
\end{equation}
where we used the formula for $S \colon I \to O(\liealg{p})$ from
Lemma~\ref{lemma:ode_rolling_along_special_curves_intrinsic_solution}.
Therefore,
\begin{equation}
	\alpha(t)
	= 
	\int_{0}^{t}
	\exp\big( s \big( \ad_{\xi_{ \liealg{h}}} + \tfrac{1}{2} (\pr_{\liealg{p}} \circ \ad_{\xi_{\liealg{p}}}) \big) \big) (\xi_{ \liealg{p}}) \D s 
\end{equation}
is the rolling curve $\alpha \colon I \to \liealg{p}$.

We summarize our findings for the intrinsic rolling of
$T_X \Stiefel{n}{k}$ over $\Stiefel{n}{k}$ in the next proposition.

\begin{proposition}
	\label{proposition:rolling_along_special_curves_stiefel_intrinsic}
	Let $\xi = (\xi_1, \xi_2) \in \liealg{so}(n) \times \liealg{so}(k)$ and $X \in \Stiefel{n}{k}$.
	Then the triple $(\alpha(t), \widehat{\alpha}(t), A(t))$
	with
	\begin{equation}
		\begin{split}
			\alpha(t) 
			&= 
			\int_{0}^{t}
			\exp\big( s \big( \ad_{\xi_{ \liealg{h}}} + \tfrac{1}{2} (\pr_{\liealg{p}} \circ \ad_{\xi_{\liealg{p}}}) \big) \big) (\xi_{ \liealg{p}}) \D s ,\\
			\widehat{\alpha}(t) 
			&=
			\pi(\exp(t \xi)) ,\\
			A(t) 
			&=
			(d_{q(t)} \pi ) \circ d_{(I_n, I_k)} L_{q(t)} \circ S(t),
		\end{split}
	\end{equation}
	for $t \in I$,
	where $q \colon I \ni t \mapsto \exp(t \xi) \exp( -t \xi_{ \liealg{h}}) \in O(n) \times O(k)$
	and 
	\begin{equation}
		S \colon I \to O(\liealg{p}),
		\quad
		t \mapsto 	\Ad_{\exp(t \xi_{\liealg{h}} )} \circ 
		\exp\big(- t \big( \ad_{\xi_{ \liealg{h}}} + \tfrac{1}{2} (\pr_{\liealg{p}} \circ \ad_{\xi_{\liealg{p}}}) \big) \big) ,
	\end{equation}
	is an intrinsic rolling of $T_X \Stiefel{n}{k} \cong \liealg{p}$
	over $\Stiefel{n}{k} \cong (O(n) \times O(k)) / H$
	viewed as normal naturally reductive homogeneous space.
\end{proposition}

\begin{remark}
	Obviously, proceeding analogously to the proof of Proposition~\ref{proposition:rolling_along_special_curves_stiefel_intrinsic},
	one derives an explicit expression for the
	intrinsic rolling $(\alpha(t), \widehat{\alpha}(t), A(t))$
	of $T_X \Stiefel{n}{k} \cong \liealg{p}$ over $\Stiefel{n}{k}$,
	where $\widehat{\alpha}(t) = \pi(\exp(t \xi))$
	for $\xi \in \liealg{so}(n) \times \liealg{so}(k)$
	for any $\alpha$-metric, where $\alpha \in \field{R} \setminus \{-1, 0\}$.
	Indeed, an explicit expression for  
	rolling of general reductive homogeneous spaces $G / H$
	whose development curve is given by
	$t \mapsto \pi(\exp(t \xi))$ for $\xi \in \liealg{g}$ is known,
	see~\cite{schlarb:2023}.
\end{remark}

 From now on, whenever convenient, we may interchangeably  use two different notations,
$\e^A$ and $\exp(A)$,
for the exponential of a matrix.

To determine an extrinsic rolling $(\beta(t), \widehat{\beta}(t), B(t), C(t))$ of $T_X \Stiefel{n}{k}$ over $\Stiefel{n}{k}$,
with respect to the Euclidean metric
whose development curve is given by 
$\widehat{\beta} \colon I \ni t \mapsto \big(\iota_X \circ \pi \big) \big(\exp(t \xi) \big)\in \Stiefel{n}{k}$, we recall
from
Corollary~\ref{corollary:extrinsic_rolling_stiefel_kinematic_equations}
that the normal part $C(t)$ is given by
\begin{equation}
	C(t) = \Phi_{(R(t), \theta(t))} \circ \exIso(t), \quad t \in I.
\end{equation}
Here $\exIso \colon I \to O(N_X \Stiefel{n}{k})$ is the solution of the initial value problem
\begin{equation}\label{equation:ivp}
	\dot{\exIso}(t) = - P_X^{\perp} \circ f_{(\xi_1(t), \xi_2(t))} \circ \exIso(t), 
	\quad \exIso(0) = \id_{N_X \Stiefel{n}{k}},
\end{equation}
and the horizontal lift $q \colon I \to O(n) \times O(k)$
of $\widehat{\beta}$ and $S(t) u(t)$ are, as in the intrinsic case, given by 
~\eqref{horiz-curve-F} and ~\eqref{equation:rolling_along_special_curves_S_u_Ad}, respectively. That is,
\begin{equation*}
	q(t) = \exp(t \xi) \exp(- t \xi_{ \liealg{h}}) \in O(n) \times O(k),
\end{equation*}
\begin{equation}
	\label{equation:extrinsic_rolling_special_curves_normal_part_xi_1_xi_2}
	\begin{split}
		(\xi_1(t), \xi_2(t)) 
		=
		S(t) u(t) 
		=&
		\Ad_{\exp(t \xi_{ \liealg{h}})}(\xi_{ \liealg{p}})\\
		=&
		\big(
		\e^{t \xi_{1, \liealg{h}}} \xi_{1, \liealg{p}} \e^{- t \xi_{1, \liealg{h}}},
		\e^{t \xi_{2, \liealg{h}}} \xi_{2, \liealg{p}} \e^{- t \xi_{2, \liealg{h}}} 
		\big).
	\end{split}
\end{equation}
In order to determine the normal part of the extrinsic rolling explicitly, we  need to solve  \eqref{equation:ivp}.

\begin{lemma}
	\label{lemma:ode_rolling_along_special_curves_extrinsic_solution}
	Let $X \in \Stiefel{n}{k}$ and 
	$\xi = (\xi_1, \xi_2) \in \liealg{so}(n) \times \liealg{so}(k)$.
	Then the initial value problem
	\begin{equation}
		\label{equation:lemma_ode_rolling_along_special_curves_extrinsic_solution_statement}
		\begin{split}
			\dot{\exIso}(t) 
			&= 
			- P_X^{\perp} \circ
			f_{ (\e^{t \xi_{1, \liealg{h}}} \xi_{1, \liealg{p}} \e^{- t \xi_{1, \liealg{h}}},
				\e^{t \xi_{2, \liealg{h}}} \xi_{2, \liealg{p}} \e^{- t \xi_{2, \liealg{h}}} ) } 
			\circ \exIso(t), \\
			\exIso(0) 
			&=
			\id_{N_X \Stiefel{n}{k}} 
		\end{split}
	\end{equation}
	has the unique solution $\exIso \colon I \to N_X \Stiefel{n}{k}$ given by
	\begin{equation}
		\label{equation:lemma_ode_rolling_along_special_curves_extrinsic_solution_statement_sol}
		\exIso(t)
		= 
		\Phi_{\exp(t \xi_{ \liealg{h}})} \circ 
		\exp\big( - t (P_X^{\perp} \circ f_{(\xi_1, \xi_2)} )\big) .
	\end{equation}
	\begin{proof}
		By direct computation, we verify that $T$
		from~\eqref{equation:lemma_ode_rolling_along_special_curves_extrinsic_solution_statement_sol} is indeed a solution.
		We first calculate two alternative formulas for $\tfrac{\D}{\D t} \Phi_{\exp(t \xi_{ \liealg{h}})}(V)$,
		with $V \in \matR{n}{k}$, as follows:
		\begin{equation}
			\label{equation:lemma_ode_rolling_along_special_curves_extrinsic_solution_comp_1}
			\begin{split}
				\tfrac{\D}{\D t} \Phi_{\exp(t \xi_{ \liealg{h}})}(V)
				&= 
				\tfrac{\D}{\D t}
				\big(\e^{t \xi_{1, \liealg{h}}} V \e^{- t \xi_{2, \liealg{h}}} \big)  \\
				&=
				\e^{t \xi_{1, \liealg{h}}} \xi_{ 1,\liealg{h}} V \e^{- t \xi_{2, \liealg{h}}}
				- \e^{t \xi_{1, \liealg{h}}} V \xi_{2, \liealg{h}} \e^{- t \xi_{2, \liealg{h}}} \\
				&=
				\e^{t \xi_{1, \liealg{h}}} \xi_{ 1,\liealg{h}} \e^{- t \xi_{1, \liealg{h}}} \e^{t \xi_{1, \liealg{h}}}  V \e^{- t \xi_{2, \liealg{h}}}
				- \e^{t \xi_{1, \liealg{h}}} V \e^{- t \xi_{2, \liealg{h}}} \e^{t \xi_{2, \liealg{h}}} \xi_{ 2,\liealg{h}} \e^{- t \xi_{2, \liealg{h}}} \\
				&=
				f_{ (\e^{t \xi_{1, \liealg{h}}} \xi_{1, \liealg{h}} \e^{- t \xi_{1, \liealg{h}}},
					\e^{t \xi_{2, \liealg{h}}} \xi_{2, \liealg{h}} \e^{- t \xi_{2, \liealg{h}}} ) } 
				\circ \Phi_{\exp(t \xi_{ \liealg{h}})}(V),
			\end{split}
		\end{equation}
		and also		\begin{equation}
			\label{equation:lemma_ode_rolling_along_special_curves_extrinsic_solution_comp_2}
			\begin{split}
				\Phi_{\exp(t \xi_{ \liealg{h}})} \circ f_{(\xi_{1, \liealg{p}}, \xi_{2, \liealg{p}})} (V)
				&=
				\e^{t \xi_{1, \liealg{h}}} \big( \xi_{1, \liealg{p}} V - V \xi_{2, \liealg{p}} \big) \e^{- t \xi_{2, \liealg{h}}} \\
				&=
				(\e^{t \xi_{1, \liealg{h}}} \xi_{1, \liealg{p}} \e^{ - t \xi_{1, \liealg{h}}}) (\e^{t \xi_{1, \liealg{h}}} V  \e^{- t \xi_{2, \liealg{h}}})  \\
				&\quad
				- 
				(\e^{t \xi_{1, \liealg{h}}} V  \e^{- t \xi_{2, \liealg{h}}} ) ( \e^{ t \xi_{2, \liealg{h}}} \xi_{2, \liealg{p}} e^{- t \xi_{2, \liealg{h}}} ) \\
				&=
				f_{ (\e^{t \xi_{1, \liealg{h}}} \xi_{1, \liealg{h}} \e^{- t \xi_{1, \liealg{h}}},
					\e^{t \xi_{2, \liealg{h}}} \xi_{2, \liealg{h}} \e^{- t \xi_{2, \liealg{h}}} ) }
				\big( \e^{t \xi_{1, \liealg{h}}} V  \e^{- t \xi_{2, \liealg{h}}}\big) \\
				&=
				f_{ (\e^{t \xi_{1, \liealg{p}}} \xi_{1, \liealg{h}} \e^{- t \xi_{2, \liealg{h}}},
					\e^{t \xi_{2, \liealg{h}}} \xi_{2, \liealg{h}} \e^{- t \xi_{2, \liealg{h}}} ) }
				\circ \Phi_{\exp(t \xi_{ \liealg{h}})} (V) \\
				&=
				\tfrac{\D}{\D t} \Phi_{\exp(t \xi_{ \liealg{h}})}(V).
			\end{split}
		\end{equation}
		By
		exploiting~\eqref{equation:lemma_ode_rolling_along_special_curves_extrinsic_solution_comp_1}
		and~\eqref{equation:lemma_ode_rolling_along_special_curves_extrinsic_solution_comp_2},
		we can write:
		\begin{equation}
			\begin{split}
				\dot{\exIso}(t)
				&= 
				\tfrac{\D}{\D t} \big( \Phi_{\exp(t \xi_{ \liealg{h}})} \circ 
				\exp\big( - t (P_X^{\perp} \circ f_{(\xi_1, \xi_2)} )\big) \big) \\				&\stackrel{\eqref{equation:lemma_ode_rolling_along_special_curves_extrinsic_solution_comp_1}}{=}
				f_{ (\e^{t \xi_{1, \liealg{h}}} \xi_{1, \liealg{h}} \e^{- t \xi_{1, \liealg{h}}},
					\e^{t \xi_{2, \liealg{h}}} \xi_{2, \liealg{h}} \e^{- t \xi_{2, \liealg{h}}} ) } 
				\circ \Phi_{\exp(t \xi_{ \liealg{h}})} \circ 	\exp\big( - t (P_X^{\perp} \circ f_{(\xi_1, \xi_2)} )\big) \big) \\
				&\quad
				-
				\Phi_{\exp(t \xi_{ \liealg{h}})} \circ (P_X^{\perp} \circ f_{(\xi_1, \xi_2)} ) \circ \exp\big( - t (P_X^{\perp} \circ f_{(\xi_1, \xi_2)} )\big) \\
				&=
				P_X^{\perp} \circ \Phi_{\exp(t \xi_{ \liealg{h}})} \circ f_{(\xi_{1, \liealg{h}}, \xi_{2, \liealg{h}})} \circ \exp\big( - t (P_X^{\perp} \circ f_{(\xi_1, \xi_2)} )\big) \big) \\
				&\quad
				-
				(P_X^{\perp} \circ \Phi_{\exp(t \xi_{ \liealg{h}})} \circ f_{(\xi_1, \xi_2)} ) \circ \exp\big( - t (P_X^{\perp} \circ f_{(\xi_1, \xi_2)} )\big) \\
				&= 
				P_X^{\perp} \circ \Phi_{\exp(t \xi_{ \liealg{h}})} \circ
				\big(f_{(\xi_{1, \liealg{h}}, \xi_{2, \liealg{h}})} -  f_{(\xi_1, \xi_2)}  \big)
				\circ \exp\big( - t (P_X^{\perp} \circ f_{(\xi_1, \xi_2)} )\big) \\
				&=
				- P_X^{\perp} \circ \Phi_{\exp(t \xi_{ \liealg{h}})} \circ f_{(\xi_{1, \liealg{p}}, \xi_{2, \liealg{p}})} \circ  \exp\big( - t (P_X^{\perp} \circ f_{(\xi_1, \xi_2)} )\big) \\				 &\stackrel{\eqref{equation:lemma_ode_rolling_along_special_curves_extrinsic_solution_comp_2}}{=}
				- P_X^{\perp} \circ
				f_{ (\e^{t \xi_{1, \liealg{h}}} \xi_{1, \liealg{p}} \e^{- t \xi_{1, \liealg{h}}},
					\e^{t \xi_{2, \liealg{h}}} \xi_{2, \liealg{p}} \e^{- t \xi_{2, \liealg{h}}} ) } 
				\circ \Phi_{\exp(t \xi_{ \liealg{h}})} \circ \exp\big( - t (P_X^{\perp} \circ f_{(\xi_1, \xi_2)} )\big) \\
				&=
				- P_X^{\perp} \circ
				f_{ (\e^{t \xi_{1, \liealg{h}}} \xi_{1, \liealg{p}} \e^{- t \xi_{1, \liealg{h}}},
					\e^{t \xi_{2, \liealg{h}}} \xi_{2, \liealg{p}} \e^{- t \xi_{2, \liealg{h}}} ) } 
				\circ \exIso(t),
			\end{split}
		\end{equation}
		where we used 
		$P_X^{\perp} \circ f_{(\xi_{1, \liealg{h}}, \xi_{2, \liealg{h}})} 
		=
		f_{(\xi_{1, \liealg{h}}, \xi_{2, \liealg{h}})} \circ P_X^{\perp}$.
		Together with the obvious observation that the initial
		condition $\exIso(0) 
		=
		\id_{N_X \Stiefel{n}{k}}$ is satisfied, this gives the desired result.
	\end{proof}
\end{lemma}

Now we are in the position to give an
explicit expression for the extrinsic rolling of $T_X \Stiefel{n}{k}$
over $\Stiefel{n}{k}$
with respect to the Euclidean metric
whose development curve is of the desired form.

\begin{proposition}
	\label{proposition:stiefel_rolling_along_special_curves_extrinsic}
	Let $\xi = (\xi_1, \xi_2) \in \liealg{so}(n) \times \liealg{so}(k)$
	and $X \in \Stiefel{n}{k}$.
	Then the quadruple $(\beta(t), \widehat{\beta}(t), B(t), C(t))$
	is an extrinsic rolling of $T_X \Stiefel{n}{k}$ over $\Stiefel{n}{k}$
	with respect to the Euclidean metric, where
	\begin{equation}
		\begin{split}
			\beta(t) 
			&=
			(d_{(I_n, I_k)} \iota_X \circ \pi)(\alpha(t)) \\
			\widehat{\beta}(t) 
			&=
			(\iota_X \circ \pi)(\exp(t \xi)) 
			= \e^{t \xi_1} X \e^{- t \xi_2} \\
			B(t) 
			&=
			(d_{(I_n, I_k)} \iota_X) \circ A(t) \circ (d_{(I_n, I_k)} \iota_X \circ \pi)^{-1} \\
			C(t)
			&=
			\Phi_{q(t)} \circ \exIso(t) ,
		\end{split}
	\end{equation}
	for $t \in I$ and
	\begin{equation}
		\begin{split}
			\alpha(t) 
			&= 
			\int_{0}^{t}
			\exp\big( s \big( \ad_{\xi_{ \liealg{h}}} + \tfrac{1}{2} (\pr_{\liealg{p}} \circ \ad_{\xi_{\liealg{p}}}) \big) \big) (\xi_{ \liealg{p}}) \D s \\
			q(t)
			&=
			\exp(t \xi) \exp( -t \xi_{ \liealg{h}})
			=
			\big( \e^{t \xi_1} \e^{- t \xi_{1, \liealg{h}}}, \e^{t \xi_2} \e^{- t \xi_{2, \liealg{h}}} \big) \\
			S(t)
			&=
			\Ad_{\exp(t \xi_{\liealg{h}} )} \circ 
			\exp\big(- t \big( \ad_{\xi_{ \liealg{h}}} + \tfrac{1}{2} (\pr_{\liealg{p}} \circ \ad_{\xi_{\liealg{p}}}) \big) \big) \\
			A(t) 
			&=
			(d_{q(t)} \pi ) \circ (d_{(I_n, I_k)} L_{q(t)}) \circ S(t) \\
			\exIso(t) 
			&=
			\Phi_{\exp(t \xi_{ \liealg{h}})} \circ 
			\exp\big( - t (P_X^{\perp} \circ f_{(\xi_1, \xi_2)} )\big) .
		\end{split}
	\end{equation}
	\begin{proof}
		This is a consequence of the above discussion.
		Essentially,
		the assertion follows by combining
		Proposition~\ref{proposition:rolling_along_special_curves_stiefel_intrinsic},
		Lemma~\ref{lemma:ode_rolling_along_special_curves_extrinsic_solution}
		and Theorem~\ref{theorem:stiefel_rolling_extrinsic_euclidean_metric}.	
	\end{proof}
\end{proposition}

Proposition~\ref{proposition:stiefel_rolling_along_special_curves_extrinsic}
implies an explicit expression for the rolling along geodesics.
In fact, by exploiting that geodesics 
on naturally reductive homogeneous spaces 
are projections of horizontal one-parameter groups,
we obtain the next corollary.

\begin{corollary}
	Let $\xi = (\xi_1, \xi_2) \in \liealg{p}$ and $X \in \Stiefel{n}{k}$.
	Then the quadruple $(\beta(t), \widehat{\beta}(t), B(t), C(t))$
	is an extrinsic rolling of $T_X \Stiefel{n}{k}$ over $\Stiefel{n}{k}$
	with respect to the Euclidean metric, where
	\begin{equation}
		\begin{split}
			\beta(t) 
			&=
			(d_{(I_n, I_k)} \iota_X \circ \pi)(t \xi_1, t \xi_2)
			= t (\xi_1 X - X \xi_2)\\
			\widehat{\beta}(t) 
			&=
			(\iota_X \circ \pi)(\exp(t \xi)) 
			= \e^{t \xi_1} X \e^{- t \xi_2} \\
			B(t) 
			&=
			(d_{(I_n, I_k)} \iota_X \circ \pi)
				 \circ (d_{(I_n, I_k)} L_{(\e^{t \xi_1} , \e^{t \xi_2})} ) \\
				 &\quad
				 \circ \exp\big(-\tfrac{1}{2}  t  (\pr_{\liealg{p}} 
				 \circ \ad_{\xi_{\liealg{p}}}) \big)
				 \circ (d_{(I_n, I_k)} \iota_X \circ \pi)^{-1} \\
			C(t)
			&=
			\Phi_{(\e^{t \xi_1} , \e^{t \xi_2})} \circ \exp\big( - t (P_X^{\perp} \circ f_{(\xi_1, \xi_2)} )\big)
		\end{split}
	\end{equation}
	for $t \in I$
	is a rolling, whose
	development curve is a geodesic.
	\begin{proof}
		Clearly, $\xi \in \liealg{p}$ implies $\xi_{\liealg{h}} = 0$.
		Thus the assertion follows by
		Proposition~\ref{proposition:stiefel_rolling_along_special_curves_extrinsic}.
	\end{proof}
\end{corollary}

\subsection{Comparison with Existing Literature}
\label{subsec:comparison_with_existing_literature}

In this final section, we relate our results with the known
rolling of Stiefel manifolds from~\cite{hueper.kleinsteuber.leite:2008}.

We discuss how the rolling of $T_X \Stiefel{n}{k}$ over $\Stiefel{n}{k}$ is related to the rolling obtained in~\cite{hueper.kleinsteuber.leite:2008}.
As in~\cite{hueper.kleinsteuber.leite:2008},
we specify $X = E = \begin{bsmallmatrix}
	I_k \\
	0
\end{bsmallmatrix}$.
It is well-known that 
\begin{equation}
	T_{E} \Stiefel{n}{k} 
	=
	\Big\{ \begin{bsmallmatrix}
		\Omega \\
		B
	\end{bsmallmatrix} 
	\mid \Omega \in \liealg{so}(k), B \in \matR{(n - k)}{k} \Big\}
	= 
	\liealg{so}(n) E
\end{equation}
holds.
We now recall the rolling map from~\cite{hueper.kleinsteuber.leite:2008},
where we make trivial modifications concerning the terminology
and notations
in order to adapt it to our notation.

Let $\alpha \colon I \to \Stiefel{n}{k}$ be a rolling curve with $\alpha(0) = E$.
Then there exists
a curve $U \colon I \to SO(n)$ such that $\alpha(t) = U(t) E$.
Denote
\begin{equation}
	G 
	=
	\{ W \in SO(n k) \mid W = V \tensor U, V \in SO(k), U \in SO(n)  \}
	\subset
	SO(n k)
\end{equation}
and 
\begin{equation}
	\mathcal{U}(t) = \{ Q(t) \in G \mid Q(t) \vect(E) = (V(t) \tensor U(t)) \vect(E) = \vect(\alpha(t)) \} .
\end{equation}
The rotational part $R(t) \in SO(nk)$ describing the rolling of $T_E \Stiefel{n}{k}$ over $\Stiefel{n}{k}$ is obtained in~\cite{hueper.kleinsteuber.leite:2008} by the following Ansatz:
\begin{equation}
	R(t) = Q(t) \widetilde{S}(t) ,
\end{equation}
where $Q(t) \in \mathcal{U}(t)$ and $\widetilde{S}(t)$ is a curve in the isotropy subgroup of $E$ under the $SO(n k)$-action on $\field{R}^{n k} = \vect(\matR{n}{k})$, i.e.
\begin{equation}
	\begin{split}
		\widetilde{S}(t) 
		&\in \{  R \in SO(n k) \mid R \vect(E) = \vect(E)\} \\
		&\quad
		= \Stab(\vect(E)) \cong SO(n k - 1) ,
	\end{split}
\end{equation}
where the isomorphism in the above equation is obtained by choosing an orthogonal transformation 
$P_0 \in O(n k)$ such that $P_0 E  \in \Span\{e_{n k}\}$ holds,
as well as 
\begin{equation}
	\label{equation:stiefel_extrinsic_classical_P_0}
	\begin{split}
		P_0 (T_E \Stiefel{n}{k}) 
		&= 
		\Span\{e_1, \ldots e_{\ell_t}\} , \\
%		\, , \,
		P_0 (N_E  \Stiefel{n}{k}) 
		&= 
		\Span\{e_{\ell_t + 1}, \ldots e_{n k}\},
	\end{split}
\end{equation}
where $\ell_t = \dim(\Stiefel{n}{k})$ and $\ell_n = \dim(N_X \Stiefel{n}{k})$,
yielding
\begin{equation}
	P_0 \big( \Stab (E)\big) P_0^{\top} 
	=
	\Big\{ \begin{bsmallmatrix}
		S && 0 \\
		0 && 1
	\end{bsmallmatrix} \mid S \in SO(n k - 1)\Big\} .
\end{equation}
Note that $\widetilde{S}$ in this text corresponds to $S$
in~\cite{hueper.kleinsteuber.leite:2008}.
By this notation, it is shown in~\cite{hueper.kleinsteuber.leite:2008}
that $\widetilde{S}(t)$ needs to fulfill
\begin{equation}
	\widetilde{S}(t) \in 
	\begin{bsmallmatrix}
		O(\ell_t) && 0 && 0 \\
		0 && O(\ell_n-1) && 0 \\
		&& 0 && 1
	\end{bsmallmatrix} \cap SO(n k),
\end{equation}
where $\ell_t = \dim(T_E \Stiefel{n}{k})$ and 
	$\ell_n = \dim(N_E \Stiefel{n}{k})$.

The orthogonal projection of a matrix $A \in \matR{nk}{nk}$ onto
a matrix with the structure given in the above equation
is denoted by $A_{\mathrm{bl-diag}}$. 
Using this notation, we recall~\cite[Lem. 3.2]{hueper.kleinsteuber.leite:2008}.

\begin{lemma}
	\label{lemma:stiefel_rolling_classical}
	Let $h = (R, s)$ be a rolling map for the
	Stiefel manifold $\Stiefel{n}{k}$. 
	If $Q(t) \in \mathcal{U}(t)$ and $R(t) = Q(t) \widetilde{S}(t)$
	with $\widetilde{S}(t) \in \Stab(E)$, then $\widetilde{S}(t)$ obeys the ODE
	\begin{equation}
		\label{equation:lemma_stiefel_rolling_classical_S}
		\begin{split}
			\dot{\widetilde{S}}(t)
			&=
			P_0^{\top} \big( P_0 \dot{Q}(t)^{\top} Q(t) P_0^{\top} \big)_{\mathrm{bl-diag}} P_0 \widetilde{S}(t) \\
			&=
			- P_0^{\top} \big( P_0 Q(t)^{\top} \dot{Q}(t) P_0^{\top} \big)_{\mathrm{bl-diag}} P_0 \widetilde{S}(t),
		\end{split}
	\end{equation}
	where $s \colon I \to \field{R}^{n k}$ fulfills the ODE
	\begin{equation}
		\label{equation:lemma_stiefel_rolling_classical_s}
		\dot{s}(t) 
		=
		- \widetilde{S}(t) \dot{Q}(t)^{\top} Q(t) \vect(E) 
		= \widetilde{S}(t) Q(t)^{\top} \dot{Q}(t) \vect(E)
	\end{equation}
	by~\cite[Eq. (44)]{hueper.kleinsteuber.leite:2008}.
\end{lemma}
Note that the second equations
in~\eqref{equation:lemma_stiefel_rolling_classical_S} and~\eqref{equation:lemma_stiefel_rolling_classical_S}
of Lemma~\ref{lemma:stiefel_rolling_classical}
are correct by $Q^\top \dot{Q} \in \liealg{so}(n k )$
because of $Q \colon I \to SO(n k)$.

The goal of the remainder part of this subsection is to show that the
extrinsic rolling of the Stiefel manifold obtained in Subsection~\ref{subsec:stiefel_extrinsic_rolling}
fulfills Lemma~\ref{lemma:stiefel_rolling_classical}.
To this end, we recall, that the extrinsic rolling 
$(\beta(t), \widehat{\beta}(t), B(t), C(t))$ from 
Subsection~\ref{subsec:stiefel_extrinsic_rolling}
is constructed by using the kinematic equations
\begin{equation}
	\label{equation:stiefel_extrinsic_rolling_kinematic_equation_comparison}
	\begin{split}
		u(t) 
		&=
		(d_{(I_n, I_k)} (\iota_E \circ \pi) \at{p})^{-1}(\dot{\beta}(t)) , \\
		\dot{S}(t) 
		&=
		- \tfrac{1}{2} \pr_{\liealg{p}} \circ \ad_{S(t) u(t)} \circ S(t), 
		\quad S(0) = S_0 = \id_{\liealg{p}}\in O(\liealg{p}), \\
		\dot{q}(t) 
		&= (d_{(I_n, I_k)} L_{q(t)} ) \circ S(t) u(t), 
		\quad q(0) = (I_n, I_k) \in O(n) \times O(k), \\
		\dot{\exIso}(t) 
		&=
		-  P^{\perp}_{E} \circ 
		f_{(\xi_1(t), \xi_2(t))} \circ \exIso(t), 
		\quad \exIso(0) = \id_{N_E \Stiefel{n}{k}} \in O(N_E \Stiefel{n}{k}) ,
	\end{split}
\end{equation}
according to
Corollary~\ref{corollary:extrinsic_rolling_stiefel_kinematic_equations}
for $X = E$.
The development curve reads
\begin{equation}
	\widehat{\beta}(t) = R(t) E \theta(t)^{\top} .
\end{equation}
Hence $q(t) = (R(t), \theta(t)) \in \mathcal{U}(t)$ is
fulfilled by the definition of $\mathcal{U}(t)$
after identifying $q(t)$ with $Q(t) = \theta(t) \tensor R(t)$
by the map
\begin{equation}
	O(n) \times O(k) \ni (R, \theta) \mapsto \theta \tensor R \in O(k) \tensor O(n)
\end{equation}
which is an isomorphism of Lie groups onto its images.
Using this identification, we obtain
that
\begin{equation}
	\big(d_e L_{q(t)} \big)^{-1} \dot{q}(t) = S(t) u(t) =  (\xi_1(t), \xi_2(t))
\end{equation}
corresponds to
\begin{equation}
	\label{equation:comparison_stiefel_xi_1_xi_2_kronecker}
	\begin{split}
		Q^{\top} \dot{Q} 
		&=
		\xi_2(t) \tensor I_n + I_k \tensor \xi_1(t)
	\end{split}
\end{equation}
by using properties of the Kronecker product, see e.g.~\cite[Sec. 7.1]{bernstein:2009}.

It remains to relate the curves $S(t)$ and $\exIso(t)$ 
from~\eqref{equation:stiefel_extrinsic_rolling_kinematic_equation_comparison}
to the curve
$\widetilde{S}(t)$ considered in Lemma~\ref{lemma:stiefel_rolling_classical}.

We first consider the normal part. 
We show that $E \in N_E \Stiefel{n}{k}$ is invariant
under $\exIso \colon I \to O(N_E \Stiefel{n}{k})$,
where $\exIso$ is defined by the kinematic equation.
We obtain by the definition of $f_{(\xi_1, \xi_2)}$ for $X = E$
\begin{equation}
	f_{(\xi_1, \xi_2)} \big( E \big)
	= 
	\big( \xi_1 E - E \xi_2 \big) 
\end{equation}
implying
$f_{(\xi_1, \xi_2)}\big(\Span\{E\} \big) \subset T_E \Stiefel{n}{k}$
by the linearity of $f_{(\xi_1, \xi_2)} \colon \matR{n}{k} \to \matR{n}{k}$.
Next, we consider the curve $I \ni t \mapsto E(t) = \exIso(t)( E) $,
where $\exIso \colon I \to O(N_X \Stiefel{n}{k})$ is given by
the kinematic equation.
We may view $E(t)$ as a solution of the initial value problem
\begin{equation}
	\dot{E}(t) 
	=
	- \big( P_E^{\perp} \circ f_{(\xi_1(t), \xi_2(t))} \big) ( E(t) ),
	\quad
	E(0) = E .
\end{equation}
The unique solution of this ODE is given by $E(t) =  E$ for $t \in I$,
since $E(0) = E$ is clearly fulfilled and
\begin{equation}
	- \big( P_E^{\perp} \circ f_{(\xi_1(t), \xi_2(t))} \big) (E)
	=
	0
	= \dot{E}(t)
\end{equation}
holds due to
$f_{(\xi_1, \xi_2)}\big(\Span\{E\} \big) \subset T_E \Stiefel{n}{k}$.
In other words, since $\exIso(0) = \id_{N_E \Stiefel{n}{k}}$, one has 
\begin{equation}
	\label{equation:stiefel_rolling_classical_E_invariant}
	\exIso(t) E = E, \quad t \in I .
\end{equation}

Clearly, by choosing $P_0 \in O(n k)$ such that~\eqref{equation:stiefel_extrinsic_classical_P_0} holds,
one obtains for $v \in \field{R}^{n k}$
\begin{equation}
	\begin{split}
		P_0 \vect( P_E^{\perp} (\vect^{-1}(v))) 
		&=
		\begin{bsmallmatrix}
			0_{\ell_t} && 0 \\
			0 && I_{\ell_{n}}
		\end{bsmallmatrix}
		P_0 v
		\in \field{R}^{n k} ,\\
		P_0 \vect(P_E (\vect^{-1}(v)))
		&=
		\begin{bsmallmatrix}
			I_{\ell_t} && 0 \\
			0 && 0_{\ell_{n}}
		\end{bsmallmatrix} P_0 v \in \field{R}^{n k} ,
	\end{split} 
\end{equation}
which implies, for $v \in \field{R}^{n k}$,
\begin{equation}
	\label{equation:stiefel_rolling_classical_projections_matrix_rep}
	\begin{split}
	\vect \circ P_E \circ \vect^{-1} (v)
	&= 
	P_0^{\top} 	\begin{bsmallmatrix}
		I_{\ell_t} && 0 \\
		0 && 0_{\ell_{n}}
	\end{bsmallmatrix} P_0 v,
	\\
%	&\\
	\vect \circ P_E^{\perp} \circ \vect^{-1} (v)
	&=
	P_0^{\top}
	\begin{bsmallmatrix}
		0_{\ell_t} && 0 \\
		0 && I_{\ell_{n}}
	\end{bsmallmatrix}
	P_0 v.
	\end{split} 
\end{equation}

We now identify the curve $S \colon I \to O(\liealg{p})$ with the curve
$\widehat{S} \colon I \to O(T_E \Stiefel{n}{k})$ via 
\begin{equation}
	 \label{equation:comparison_tangential_part_rectangular_matrices}
	\widehat{S}(t) 
	=
	\big( d_{(I_n, I_k)} (\iota_E \circ \pi) \big)
	 \circ S(t) \circ 
	 \big(d_{(I_n, I_k)} (\iota_E \circ \pi)\at{\liealg{p}} \big)^{-1} .
\end{equation}
Afterwards, we find a matrix representation for $\widehat{S}$,
roughly speaking,
via considering $\mathcal{S} = \vect \circ \widehat{S} \circ \vect^{-1}$.

We start with
computing~\eqref{equation:comparison_tangential_part_rectangular_matrices} explicitly.
The ODE~\eqref{equation:stiefel_extrinsic_rolling_kinematic_equation_comparison}
for $S(t) \in O(\liealg{p})$ can be equivalently rewritten as

\begin{equation}
	\dot{S}(t) \circ S(t)^{-1}
	=
	- \tfrac{1}{2} \pr_{\liealg{p}} \circ \ad_{(\xi_1(t), \xi_2(t))} 
\end{equation}
and therefore
\begin{equation}
	\label{equation:comparison_hat_S_ode}
	\begin{split}
	&\dot{\widehat{S}}(t) \circ \widehat{S}(t)^{-1} \\
	&=
	\big(d_{(I_n, I_k)} (\iota_E \circ \pi)\big)
	\circ 
	\big( - \tfrac{1}{2} \pr_{\liealg{p}} \circ \ad_{(\xi_1(t), \xi_2(t))} \big)
	\circ 
	(d_{(I_n, I_k)} (\iota_E \circ \pi)\at{\liealg{p}})^{-1} \\
	&=
	\big(d_{(I_n, I_k)} (\iota_E \circ \pi)\big)
	\circ 
	\big( - \tfrac{1}{2}  \ad_{(\xi_1(t), \xi_2(t))} \big)
	\circ 
	(d_{(I_n, I_k)} (\iota_E \circ \pi)\at{\liealg{p}})^{-1} ,
	\end{split}
\end{equation} 
where, for the last equality, we use the fact that $\liealg{h}$ belongs
to the kernel of $d_{(I_n, I_k)} (\iota_E \circ \pi)$.

We now compute the right-hand side of the above equation.
To this end, we write
\begin{equation}
	V = \begin{bsmallmatrix}
		\Omega \\
		C
	\end{bsmallmatrix}
	\in T_E \Stiefel{n}{k}
	\quad \text{ and }  \quad
	(\xi_1, \xi_2) 
	=
	\Big(
	\begin{bsmallmatrix}
		2 \Psi && - B^{\top} \\
		B && 0
	\end{bsmallmatrix},
	\Psi \Big)
	\in \liealg{p}.
\end{equation}
Taking into account that $\Omega^T=-\Omega$, $\Psi^T=-\Psi$, and 
\begin{equation}
	d_{(I_n, I_k)} (\iota_E \circ \pi)\Big(
	\begin{bsmallmatrix}
		2 \Psi && - B^{\top} \\
		B && 0
	\end{bsmallmatrix},
	\Psi \Big)= \begin{bsmallmatrix}
		2 \Psi && - B^{\top} \\
		B && 0
	\end{bsmallmatrix}\begin{bsmallmatrix}
		I_k \\
		 0
	\end{bsmallmatrix}-
	\begin{bsmallmatrix}
		I_k \\
		 0
	\end{bsmallmatrix}\Psi =\begin{bsmallmatrix}
		\Psi \\
		 B
	\end{bsmallmatrix},
\end{equation}
we can write
\begin{equation}
	\label{equation:comparison_tangental_part_identification}
	\begin{split}
		&\big(d_{(I_n, I_k)} (\iota_E \circ \pi) \big) \circ \dot{S}(t) \circ S(t)^{-1} \circ (d_{(I_n, I_k)} (\iota_E \circ \pi)\at{\liealg{p}})^{-1}(V) \\
		&=
		- \tfrac{1}{2} d_{(I_n, I_k)} (\iota_X \circ \pi) 
		\Big( [\xi_1, ( E V^{\top} - V E^{\top} ], [\xi_2, E^{\top} V ] \Big) \\
		&=
		- \tfrac{1}{2} d_{(I_n, I_k)} (\iota_X \circ \pi)
		\Big( \begin{bsmallmatrix}
			2 \Psi && - B^{\top} \\
			B && 0
		\end{bsmallmatrix}
		\begin{bsmallmatrix}
			2 \Omega && - C^{\top} \\
			C && 0
		\end{bsmallmatrix}
		-\begin{bsmallmatrix}
		2 \Omega && - C^{\top} \\
		C && 0
		\end{bsmallmatrix} 
		\begin{bsmallmatrix}
			2 \Psi && - B^{\top} \\
			B && 0
		\end{bsmallmatrix},
		\Psi \Omega - \Omega \Psi \Big) \\
		&=
		- \tfrac{1}{2} d_{(I_n, I_k)} (\iota_X \circ \pi)
		\Big(
		\begin{bsmallmatrix}
			4 \Psi \Omega - B^{\top} C && - 2 \Psi C^{\top} \\
			2 B \Omega && - B C^{\top}
		\end{bsmallmatrix}
		- 
		\begin{bsmallmatrix}
			4 \Omega \Psi - C^{\top} B && - 2 \Omega B^{\top}\\
			2 C  \Psi && - C B^{\top}
		\end{bsmallmatrix}, 
		\Psi \Omega - \Omega \Psi
		\Big) \\
		&=
		- \tfrac{1}{2} d_{(I_n, I_k)} (\iota_X \circ \pi)
		\Big(
		\begin{bsmallmatrix}
			4 \Psi \Omega - B^{\top} C - 4 \Omega \Psi + C^{\top} B && - 2 \Psi C^{\top} + 2 \Omega B^{\top} \\
			2 B \Omega - 2 C \Psi && - B C^{\top} + C B^{\top}
		\end{bsmallmatrix},
		\Psi \Omega - \Omega \Psi\Big) \\
		&=
		- 
		\tfrac{1}{2} \Big(
		\begin{bsmallmatrix}
			4 \Psi \Omega - B^{\top} C - 4 \Omega \Psi + C^{\top} B && - 2 \Psi C^{\top} + 2 \Omega B^{\top} \\
			2 B \Omega - 2 C \Psi && - B C^{\top} + C B^{\top}
		\end{bsmallmatrix}
		\begin{bsmallmatrix}
			I_k \\
			0
		\end{bsmallmatrix}
		- 
		\begin{bsmallmatrix}
			I_k \\
			0
		\end{bsmallmatrix}
		(	\Psi \Omega - \Omega \Psi) \Big) \\
		&=
		-\tfrac{1}{2} 
		\begin{bsmallmatrix}
			4 \Psi \Omega - B^{\top} C - 4 \Omega \Psi + C^{\top} B - \Psi \Omega + \Omega \Psi \\
			2 B \Omega - 2 C \Psi
		\end{bsmallmatrix} \\
		&=
		- \tfrac{1}{2} 
		\begin{bsmallmatrix}
				3 \Psi \Omega - 3 \Omega \Psi - B^{\top} C  + C^{\top} B   \\
			2 B \Omega - 2 C \Psi
		\end{bsmallmatrix},
	\end{split}
\end{equation}
as well as
\begin{equation}
	\label{equation:comparison_tangential_part_projection_xi_1_V_V_x2}
	\begin{split}
		P_E\Big(  \xi_1 V - V \xi_2 \Big)
		&=
		P_E\Big(
		\begin{bsmallmatrix}
			2 \Psi && - B^{\top} \\
			B && 0
		\end{bsmallmatrix}
		\begin{bsmallmatrix}
			\Omega \\
			C
		\end{bsmallmatrix}
		- 
		\begin{bsmallmatrix}
			\Omega \\
			C
		\end{bsmallmatrix}
		\Psi \Big) \\
		&=
		P_E\Big(
		\begin{bsmallmatrix}
			2 \Psi \Omega - B^{\top} C - \Omega \Psi \\
			B \Omega - C \Psi
		\end{bsmallmatrix}
		\Big) \\
		&=
		\begin{bsmallmatrix}
			(\Psi \Omega - (\Psi \Omega)^{\top}) 
			- \tfrac{1}{2} (B^{\top} C - (B^{\top} C)^{\top})
			- \tfrac{1}{2} (\Omega \Psi - (\Psi \Omega)^{\top}) \\
			B \Omega - C \Psi
		\end{bsmallmatrix} \\
		&=
		\begin{bsmallmatrix}
			\Psi \Omega - \Omega \Psi 
			- \tfrac{1}{2} B^{\top} C  + \tfrac{1}{2}C^{\top} B
			- \tfrac{1}{2} \Omega \Psi + \tfrac{1}{2} \Psi \Omega \\
			B \Omega - C \Psi
		\end{bsmallmatrix} \\
		&=- \tfrac{1}{2}
		\begin{bsmallmatrix}
			3 \Psi \Omega - 3 \Omega \Psi 
			 -B^{\top} C  + C^{\top} B \\
			2B \Omega - 2C \Psi
		\end{bsmallmatrix}.
	\end{split}
\end{equation}
By comparing~\eqref{equation:comparison_tangental_part_identification}
and~\eqref{equation:comparison_tangential_part_projection_xi_1_V_V_x2},
we get 
\begin{equation}
	\begin{split}
	&\big(d_{(I_n, I_k)} (\iota_X \circ \pi) \big) \circ \dot{S}(t) \circ S(t)^{-1} \circ (d_{(I_n, I_k)} (\iota_X \circ \pi)\at{\liealg{p}})^{-1}(V) \\
	&=
		- P_E\big(\xi_1(t) V - V \xi_2(t) \big)
	\end{split}
\end{equation}
Therefore,~\eqref{equation:comparison_hat_S_ode} can be written as 
\begin{equation}
	\dot{\widehat{S}}(t) \circ \widehat{S}(t)^{-1}(V)
	=
	-
	P_E\big(\xi_1(t) V - V \xi_2(t) \big)
\end{equation}
for $V \in T_E \Stiefel{n}{k}$, or equivalently as 
\begin{equation}
	\dot{\widehat{S}}(t) \circ \widehat{S}(t)^{-1}(V)
	=
	-
	P_E \circ f_{(\xi_1(t), \xi_2(t))} (V)
\end{equation}
for $V \in T_E \Stiefel{n}{k}$.
Applying $\vect \colon \matR{n}{k} \to \field{R}^{n k}$ we get
for $\mathcal{S} = \vect \circ \widehat{S}(t) \circ \vect^{-1}$
the ODE
\begin{equation}
	\dot{\mathcal{S}}(t)
	=
	- \big(\vect \circ P_E \circ \vect^{-1} \big)
	\circ
	\big( \vect \circ  f_{(\xi_1(t), \xi_2(t))} \circ \vect^{-1} \big) 
	\circ 
	\mathcal{S}(t) .
\end{equation}
For $W\in \matR{n}{k}$, we have
\begin{equation}\label{FFF}
	\begin{split}
		\vect\big( f_{(\xi_1(t), \xi_2(t))}(W)\big)
		&= 
		\vect( \xi_1(t) W - W \xi_2(t)) \\
		&= 
		(I_k \tensor \xi_1(t) + \xi_2(t) \tensor I_n) \vect(W).
	\end{split}
\end{equation}
Denoting the representation matrix of $\mathcal{S}$ by $\mathcal{S}$, as well, and using the identity \eqref{FFF} with $W$ replaced by $\vect^{-1} \circ \mathcal{S}(t) \circ \vect (V)$, we get
\begin{equation}
	\label{equation:comparison_tangential_S_matrix_first_version}
	\dot{\mathcal{S}}(t) \vect(V)
	=
	-\big(\vect \circ P_E \circ \vect^{-1} \big)
	  (I_k \tensor \xi_1(t) + \xi_2(t) \tensor I_n)
	   \mathcal{S} \vect(V),
\end{equation}
for $V \in T_E \Stiefel{n}{k}$.

Recalling the definition of $P_0 \in SO(nk)$
from~\eqref{equation:stiefel_extrinsic_classical_P_0},
and using~\eqref{equation:stiefel_rolling_classical_projections_matrix_rep},
we can rewrite~\eqref{equation:comparison_tangential_S_matrix_first_version}
for $V \in T_E \Stiefel{n}{k}$ equivalently as 
\begin{equation}
	\begin{split}
		&\dot{\mathcal{S}}(t) \circ \mathcal{S}(t)^{-1}  \vect(V) \\
		&=
		- P_0^{\top} 
		\begin{bsmallmatrix}
			I_{\ell_t} && 0 \\
			0 && 0
		\end{bsmallmatrix}
		P_0
		(I_k \tensor \xi_1(t) + \xi_2(t) \tensor I_n)
		 P_0^{\top} P_0 \vect(V) \\
		&=
		- P_0^{\top} 
		\begin{bsmallmatrix}
		I_{\ell_t} && 0 \\
		0 && 0
		\end{bsmallmatrix}
		P_0
		(I_k \tensor \xi_1(t) + \xi_2(t) \tensor I_n)
		 P_0^{\top} 
	 	\begin{bsmallmatrix}
		 	I_{\ell_t} && 0 \\
		 	0 && 0
		 \end{bsmallmatrix}
	 	P_0 \vect(V) \\
	 	&=
	 	- P_0^{\top} 
	 	\begin{bsmallmatrix}
	 		I_{\ell_t} && 0 \\
	 		0 && 0
	 	\end{bsmallmatrix}
 		P_0
	 	Q(t)^{\top} \dot{Q}(t)
	 	P_0^{\top} 
	 	\begin{bsmallmatrix}
	 		I_{\ell_t} && 0 \\
	 		0 && 0
	 	\end{bsmallmatrix}
	 	P_0 \vect(V) ,
	\end{split}
\end{equation}
where the last equality holds due
to~\eqref{equation:comparison_stiefel_xi_1_xi_2_kronecker}.

Similarly,  for $\exIso \colon I \to O(N_E \Stiefel{n}{k})$,
if we define
$\mathcal{\exIso}(t) = \vect \circ \exIso(t) \circ \vect^{-1}$
and denote
its representation matrix
 by the same symbol, we have, for $V \in N_X \Stiefel{n}{k}$,
\begin{equation}
	\begin{split}
		&\dot{\mathcal{\exIso}}(t) \circ \mathcal{\exIso}(t)^{-1} \vect(V) \\
		&= 
		- 
		P_0^{\top}
		\begin{bsmallmatrix}
			0 && 0 \\
			0 && I_{\ell_n}
		\end{bsmallmatrix}
		P_0
		(I_k \tensor \xi_1(t) + \xi_2(t) \tensor I_n)
		P_0^{\top} 
		\begin{bsmallmatrix}
			0 && 0 \\
			0 && I_{\ell_n}
		\end{bsmallmatrix}
		P_0 \vect(V) \\
		&=
		- 
		P_0^{\top}
		\begin{bsmallmatrix}
			0 && 0 \\
			0 && I_{\ell_n}
		\end{bsmallmatrix}
		P_0
		Q(t)^{\top} \dot{Q}(t)
		P_0^{\top} 
		\begin{bsmallmatrix}
			0 && 0 \\
			0 && I_{\ell_n}
		\end{bsmallmatrix}
		P_0 \vect(V).
	\end{split}
\end{equation}

Next we define  $\widetilde{S} \colon I \ni t \mapsto\widetilde{S}(t)  \in \matR{nk}{nk}$
and show  that this curve $\widetilde{S}(t)$ is exactly the curve $\widetilde{S}(t)$
from Lemma~\ref{lemma:stiefel_rolling_classical}. For that, let $v \in \field{R}^{n k}$ and compute 
\begin{equation}
	\label{equation:comparison_tilde_S_definition}
	\begin{split}
		\dot{\widetilde{S}}(t) \widetilde{S}(t)^{-1} v
		&= 
		\dot{\mathcal{S}}(t) \circ \mathcal{S}(t)^{-1} \circ
		 ( \vect \circ P_E \circ \vect^{-1})(v) \\
		 &\quad
		+
		\dot{\mathcal{\exIso}}(t) \circ \mathcal{\exIso}(t)^{-1} \circ 
		( \vect \circ P_E^{\perp} \circ \vect^{-1})(v)	\\
		&=
			- P_0^{\top} 
		\begin{bsmallmatrix}
			I_{\ell_t} && 0 \\
			0 && 0
		\end{bsmallmatrix}
		P_0
		Q(t)^{\top} \dot{Q}(t)
		P_0^{\top} 
		\begin{bsmallmatrix}
			I_{\ell_t} && 0 \\
			0 && 0
		\end{bsmallmatrix}
		P_0 v \\
		&\quad
		-
		P_0^{\top}
		\begin{bsmallmatrix}
			0 && 0 \\
			0 && I_{\ell_n}
		\end{bsmallmatrix}
		P_0
		Q(t)^{\top} \dot{Q}(t)
		P_0^{\top} 
		\begin{bsmallmatrix}
			0 && 0 \\
			0 && I_{\ell_n}
		\end{bsmallmatrix}
		P_0 v \\
		&=
		- P_0^{\top}
		\Big(
		\begin{bsmallmatrix}
			I_{\ell_t} && 0 \\
			0 && 0
		\end{bsmallmatrix}
		P_0
		Q(t)^{\top} \dot{Q}(t)
		P_0^{\top} 
		\begin{bsmallmatrix}
			I_{\ell_t} && 0 \\
			0 && 0
		\end{bsmallmatrix} \\
		&\quad
		+
		\begin{bsmallmatrix}
			0 && 0 \\
			0 && I_{\ell_n}
		\end{bsmallmatrix}
		P_0
		Q(t)^{\top} \dot{Q}(t)
		P_0^{\top} 
		\begin{bsmallmatrix}
			0 && 0 \\
			0 && I_{\ell_n}
		\end{bsmallmatrix}
		\Big)
		P_0 v.
	\end{split}
\end{equation}
In order to show that $\widetilde{S}(t)$ satisfies indeed the ODE from Lemma~\ref{lemma:stiefel_rolling_classical}
we state the following auxiliar result.

\begin{lemma}
	\label{lemma:comparision_lemma_bl_diag_alternative_formula}
	Let $\ell_t, \ell_n \in \field{N}$ with $\ell_t + \ell_n = n k$,
	and consider the matrix $A\in \liealg{so}(nk)$ partitioned as
	$$A = \begin{bsmallmatrix}
		 A_{11} && A_{12} && 0 \\
		 -A^T_{12} && A_{22} && 0\\
		 0 && 0 && 0
	\end{bsmallmatrix} , \quad \mbox{where}\quad A_{11} \in \liealg{so}(\ell_t), \quad A_{22} \in \liealg{so}(\ell_n - 1).
	$$
	Then, for $v \in \field{R}^{n k}$,
	\begin{equation}
		A_{\mathrm{bl-diag}} v 
		=
		\Big(\begin{bsmallmatrix}
			0 && 0 \\
			0 && I_{\ell_n}
		\end{bsmallmatrix}
		A 
		\begin{bsmallmatrix}
			0 && 0 \\
			0 && I_{\ell_n}
		\end{bsmallmatrix}
		+ 
		\begin{bsmallmatrix}
			I_{\ell_t} && 0 \\
			0 && 0
		\end{bsmallmatrix} 
		A 
			\begin{bsmallmatrix}
			I_{\ell_t} && 0 \\
			0 && 0
		\end{bsmallmatrix} \Big) v
	\end{equation}
	holds.
	\begin{proof}
		Writing $v = \begin{bsmallmatrix}
			v_1 \\
			v_2 \\
			v_3
		\end{bsmallmatrix}$,
		where $v_1 \in \field{R}^{\ell_t}$, $v_2 \in \field{R}^{\ell_n - 1}$ and $v_3 \in \field{R}$, we compute
	\begin{equation}
		\begin{split}
			A_{\mathrm{bl-diag}} v
			=
			\begin{bsmallmatrix}
				A_{11} && 0 && 0 \\
				0 && A_{22} && 0 \\
				0 && 0 && 0
			\end{bsmallmatrix}
			\begin{bsmallmatrix}
				v_1 \\
				v_2 \\
				v_3
			\end{bsmallmatrix}
			=
			\begin{bsmallmatrix}
				A_{11} v_1  \\
				A_{22} v_2 \\
				0
			\end{bsmallmatrix}
		.
		\end{split}
	\end{equation}
	Moreover, we also have
	\begin{equation}
		\begin{split}
			&\Big(\begin{bsmallmatrix}
				0 && 0 \\
				0 && I_{\ell_n}
			\end{bsmallmatrix}
			A 
			\begin{bsmallmatrix}
				0 && 0 \\
				0 && I_{\ell_n}
			\end{bsmallmatrix}
			+ 
			\begin{bsmallmatrix}
				I_{\ell_t} && 0 \\
				0 && 0
			\end{bsmallmatrix} v
			A 
			\begin{bsmallmatrix}
				I_{\ell_t} && 0 \\
				0 && 0
			\end{bsmallmatrix} \Big)
			\begin{bsmallmatrix}
				v_1 \\
				v_2 \\
				v_3
			\end{bsmallmatrix} \\
			&=
			\begin{bsmallmatrix}
				0 && 0 \\
				0 && I_{\ell_n}
			\end{bsmallmatrix}
			\begin{bsmallmatrix}
				A_{11} && A_{12} && 0 \\
				A_{21} && A_{22} && 0\\
				0 && 0 && 0
			\end{bsmallmatrix}
			\begin{bsmallmatrix}
				0 \\
				v_2 \\
				v_3
			\end{bsmallmatrix}
			+
			\begin{bsmallmatrix}
				I_{\ell_t} && 0 \\
				0 && 0
			\end{bsmallmatrix}
			\begin{bsmallmatrix}
				A_{11} && A_{12} && 0 \\
				A_{21} && A_{22} && 0\\
				0 && 0 && 0
			\end{bsmallmatrix}
			\begin{bsmallmatrix}
				v_1 \\
				0 \\
				0
			\end{bsmallmatrix} \\
			&= 
			\begin{bsmallmatrix}
				0 && 0 \\
				0 && I_{\ell_n}
			\end{bsmallmatrix}
			\begin{bsmallmatrix}
				A_{12} v_2 \\
				A_{22} v_2 \\
				0
			\end{bsmallmatrix}
			+
			\begin{bsmallmatrix}
				I_{\ell_t} && 0 \\
				0 && 0
			\end{bsmallmatrix}
			\begin{bsmallmatrix}
				A_{11} v_1  \\
				A_{21} v_1 \\
				A_{31} v_1
			\end{bsmallmatrix} \\
			&= \begin{bsmallmatrix}
				A_{11} v_1  \\
				A_{22} v_2 \\
				0
							\end{bsmallmatrix},
		\end{split}
	\end{equation}
	showing the desired result.
	\end{proof}
\end{lemma}

Applying Lemma~\ref{lemma:comparision_lemma_bl_diag_alternative_formula}
to~\eqref{equation:comparison_tilde_S_definition} 
yields
\begin{equation}
	\dot{\widetilde{S}}(t) {\widetilde{S}}(t)^{-1}
	=
	- P_0^{\top}  \big( P_0 Q(t)^{\top} \dot{Q}(t) P_0^{\top} \big)_{\mathrm{bl-diag}}  P_0 .
\end{equation}
So,
%\\
%
%\noindent\fbox{%
 %   \parbox{\textwidth}{%
 %\begin{quotation}
 $\widetilde{S}$ defined in~\eqref{equation:comparison_tilde_S_definition} fulfills the ODE from Lemma \ref{lemma:stiefel_rolling_classical}.
%\end{quotation}
%}}
%\\

It remains to show that our approach also gives the curve $s \colon I \to \field{R}^{n k}$ from
Lemma~\ref{lemma:stiefel_rolling_classical}.
Recalling that
\begin{equation}
		\big(d_e L_{q(t)} \big)^{-1} \dot{q}(t) = S(t) u(t) =  (\xi_1(t), \xi_2(t)),
\end{equation}
we write $S(t)^{-1}(\xi_1(t), \xi_2(t)) = (u_1(t), u_2(t))$,
and
\begin{equation}
	\dot{\beta}(t) 
	=
	\big(d_{(I_n, I_k)} (\iota_E \circ \pi) \big)\big(u_1(t), u_2(t) \big)
	=
	 u_1(t) E - E u_2(t) ,
\end{equation}
where $\beta$ is the rolling curve for the rolling of $T_E \Stiefel{n}{k}$
over $\Stiefel{n}{k}$.
We now consider the curve $s \colon I \to \field{R}^{n k}$ from
Lemma~\ref{lemma:stiefel_rolling_classical}, and do the following computations:
\begin{equation}
	\label{equation:comparison_dot_small_s_and_control_equal}
	\begin{split}
		\dot{s}(t) 
		&=
		\widetilde{S}(t)^{\top} Q(t)^{\top} \dot{Q}(t) \vect(E) \\
		&=
		\widetilde{S}(t)^{\top} \vect(\xi_1(t) E - E \xi_2(t) ) \\
		&=
		\vect \circ
		\big( d_{(I_n, I_k)} (\iota_E \circ \pi) \big) \circ S(t)^{-1} \circ 
		\big(d_{(I_n, I_k)} (\iota_E \circ \pi\at{\liealg{p}})\big)^{-1} \\
		&\quad \circ
		\vect^{-1} \circ \vect(\xi_1(t) E - E \xi_2(t) ) \\
		&=
		\vect \circ  \big( d_{(I_n, I_k)} (\iota_E \circ \pi) \big) \circ S(t)^{-1} (\xi_1(t), \xi_2(t)) \\
		&=
		\vect\big( d_{(I_n, I_k)} (\iota_E \circ \pi) (u_1(t), u_2(t)) \big) \\
		&=
		\vect(u_1(t) E - E u_2(t)) \\
		&=
		\vect(\dot\beta(t)) .
	\end{split}
\end{equation}
By~\eqref{equation:comparison_dot_small_s_and_control_equal}, $\vect(\beta(t)) + b_0 = s(t)$ holds for $t \in I$,
and some $b_0 \in \field{R}^{n k}$.

Recalling, from Lemma~\ref{lemma:stiefel_rolling_classical},
 that $(R(t)^{\top},s(t))$
defines a rolling of $\Stiefel{n}{k}$ over $T_E \Stiefel{n}{k}$, 
the development curve is given by $Q(t)\vect(E) = \vect(\widehat{\beta}(t))$
and the rolling curve
by $s(t) = \vect(\beta(t))$. 
Thus, $\widehat{\alpha}(t)$, $\alpha(t)$ and $R(t)$ from Proposition~\ref{proposition:relation_extrinsic_rolling}
correspond to $\vect(\beta(t))$, $Q(t) \vect(E)$ and
$(Q(t) \widetilde{S}(t))^{\top}$, respectively.
Therefore, we obtain
\begin{equation}
	\label{equation:comparison_rolling_s_translation}
	\begin{split}
		s(t) 
		&= 
		\vect(\beta(t))- (Q(t) \widetilde{S}(t))^{\top} Q(t) \vect(E) \\
		&= 
		\vect(\beta(t)) - \widetilde{S}(t)^{\top} Q(t)^{\top} Q(t) \vect(E) \\	
		&=
		\vect(\beta(t))  - \widetilde{S}(t)^{\top} Q(t)^{\top} Q(t) \widetilde{S}(t) \vect(E) \\
		&=
		\vect(\beta(t))  - E, 
	\end{split}
\end{equation}
by exploiting  that $\widetilde{S}(t)\vect(E) = \vect(E)$.
Obviously, using \eqref{equation:comparison_dot_small_s_and_control_equal} we may conclude that $s(t)$ from~\eqref{equation:comparison_rolling_s_translation}
fulfills the ODE
\begin{equation}
	\dot{s}(t) = \widetilde{S}(t)^{\top} Q(t)^{\top} \dot{Q}(t) \vect(E)
\end{equation}
from Lemma~\ref{lemma:stiefel_rolling_classical}.

%, i.e.,
%%\\
%%
%%\noindent\fbox{%
%%    \parbox{\textwidth}{%
%% \begin{quotation}
% $s(t)$ defined in~\eqref{equation:comparison_rolling_s_translation}
% fullfills the ODE from Lemma~\ref{lemma:stiefel_rolling_classical}.
%\end{quotation}
%}}
%\vspace*{0,5 cm}

%{\bf In conclusion:}\\

%\noindent\fbox{%
%    \parbox{\textwidth}{%
% \begin{quotation}
 In conclusion, after having developed theoretical results for rolling normal naturally reductive homogeneous spaces over their tangent spaces, we specialized to the Stiefel manifold.  The results  presented here for rolling extrinsically the Stiefel manifold $\Stiefel{n}{k}$ over its tangent space $T_E \Stiefel{n}{k}$ coincide with those obtained previously in~\cite{hueper.kleinsteuber.leite:2008}.
%  \end{quotation}
%}}

\section*{Acknowledgments}
The first two authors have been supported by the 
German Federal Ministry of Education and Research
(BMBF-Projekt 05M20WWA: Verbundprojekt 05M2020 - DyCA).

The third author was partially supported by the project
Pure Mathematics in Norway,
funded by Trond Mohn Foundation and Troms{\o} Research Foundation.

The fourth author  thanks Funda\c{c}\~ao para a Ci\^encia e Tecnologia (FCT) and COMPETE 2020 program for the financial support to the project UIDP/00048/2020.

%%%%%%%%%%%%%%%%%%%%%%%%%%%%%

\end{document}